\documentclass[12pt]{article}
\usepackage{amssymb,amsmath,epsfig,graphics}

\newcounter{sscounter}[section]
\renewcommand{\thesscounter}{\thesection.\arabic{sscounter}}
\newcommand{\ssnnl}[1]{\noindent
\refstepcounter{sscounter}\bf\thesscounter.#1\rm}

\newcommand{\proof}{\noindent\textbf{Proof. }}
\newcommand{\qed}{{\hfill$\blacksquare$}}

\newcommand{\CC}{{\mathbb C }}

\newcommand{\KK}{{\mathbb K }}
\newcommand{\LL}{{\mathbb L }}
\newcommand{\MM}{{\mathbb M }}
\newcommand{\NN}{{\mathbb N }}
\newcommand{\PP}{{\mathbb P }}

\newcommand{\RR}{{\mathbb R }}
\newcommand{\TT}{{\mathbb T }}

\newcommand{\ZZ}{{\mathbb Z }}

\newcommand{\cA}{\mathcal{A}}

\newcommand{\cG}{\mathcal{G}}

\newcommand{\cH}{\mathcal{H}}

\newcommand{\cL}{\mathcal{L}}

\newcommand{\cS}{\mathcal{S}}
\newcommand{\cT}{\mathcal{T}}

\newcommand{\RS}{\mathcal{M}}
\newcommand{\wRS}{\widetilde{\RS}}

\newcommand{\Sb}{\mathcal{S}^\bullet}
\newcommand{\Sc}{\mathcal{S}^\circ}
\newcommand{\Sd}{\mathcal{S}^*}

\newcommand{\sB}{\mathsf{B}}
\newcommand{\sW}{\mathsf{W}}
\newcommand{\sJ}{\mathsf{J}}

\newcommand{\sQ}{\mathsf{Q}}

\newcommand{\sX}{\mathsf{X}}

\newcommand{\va}{\mathbf{a}}
\newcommand{\vb}{\mathbf{b}}
\newcommand{\vc}{\mathbf{c}}

\newcommand{\ve}{\mathbf{e}}

\newcommand{\vg}{\mathbf{g}}
\newcommand{\vm}{\mathbf{m}}
\newcommand{\vp}{\mathbf{p}}
\newcommand{\vr}{\mathbf{r}}
\newcommand{\vs}{\mathbf{s}}

\newcommand{\vu}{\mathbf{u}}
\newcommand{\vv}{\mathbf{v}}
\newcommand{\vw}{\mathbf{w}}
\newcommand{\vx}{\mathbf{x}}

\newcommand{\vz}{\mathbf{z}}

\newcommand{\wt}{\varpi}
\newcommand{\crit}{\mathrm{crit}}
\newcommand{\nv}{\mathbf{0}}

\newcommand{\gf}{\varphi}
\newcommand{\gl}{\lambda}
\newcommand{\gs}{\sigma}

\newcommand{\gG}{\Gamma}

\newcommand{\wQ}{\widehat{\Gamma}^*}
\newcommand{\oQ}{\overline{\Gamma}^*}
\newcommand{\wG}{\widehat{\Gamma}^{\bullet\circ}}
\newcommand{\oG}{\overline{\Gamma}^{\bullet\circ}}
\newcommand{\oLL}{\overline{\LL}}

\newcommand{\Qg}{{\Gamma}^*}
\newcommand{\Gg}{{\Gamma}^{\bullet\circ}}

\newcommand{\rank}{\mathrm{rank}\:}

\newcommand{\modquot}[2]{\mbox{\raisebox{.2ex}{$#1$}\hspace{-.3em}/ \hspace{-.6em} \raisebox{-.2ex}{$#2$}}}

\begin{document}

\title{Hypergeometric Systems in two Variables, Quivers, Dimers and Dessins d'Enfants.}

\author{Jan Stienstra\\
\small Mathematisch Instituut, Universiteit Utrecht, the Netherlands\\ 
\small e-mail: {stien}{`at'}{math.uu.nl} \normalsize}

\date{}

\maketitle

\begin{abstract}This paper presents some parallel developments in 
Quiver/Dimer Models, Hypergeometric Systems
and Dessins d'Enfants. It demonstrates that the setting in which Gelfand, Kapranov and Zelevinsky have formulated the theory of hypergeometric systems,
provides also a natural setting for dimer models. 
The Fast Inverse Algorithm of \cite{HV} and the untwisting procedure of \cite{FHKV} are recasted in this more natural setting and then immediately
produce from the quiver data the Kasteleyn matrix for dimer models, which 
is best viewed 
as the bi-adjacency matrix for the untwisted model.
Some perfect matchings in the dimer models are direct reformulations of the triangulations in GKZ theory and the
rule which maps triangulations to the vertices of the secondary polygon 
extends to a rule for mapping perfect matchings to lattice points in the secondary polygon.
Finally it is observed in many examples and then conjectured to hold in general,
that the determinant of the Kasteleyn matrix
with suitable weights becomes after a simple transformation equal to the principal $\cA$-determinant in GKZ theory.
Illustrative examples are distributed throughout the text.
\end{abstract}

\section{Introduction}
 
In the last decade interesting correspondences were discovered relating Quiver Gauge Theories, lattice polygons and Calabi-Yau singularities. The motivation and evolution of these ideas in physics
are well-documented in many articles; 
e.g. \cite{BFHMS,FFHH,FHKV,FHKVW,FHMSVW,HHV,HK,HV,Ky,Kn,KOS,MS,ORV}. 
In the present paper we want to put some aspects of these correspondences
from physics alongside the 
hypergeometric systems in two variables of
Gelfand, Kapranov, Zelevinsky \cite{gkz1,gkz2,gkz3} and the dessins d'enfants 
of Grothendieck et al. \cite{LZ,Des,SV}. This reveals intriguing connections between these fields.

The beautiful insight of Gelfand, Kapranov and Zelevinsky was that hypergeometric structures greatly simplify if one introduces extra variables and balances this with an appropriate torus action \cite{gkz1,gkz2,gkz3,gkz4,S1}.
In order to profit from the simplication they developed tools
like the \emph{secondary fan, secondary polytope} and
\emph{principal $\cA$-determinant}. This paper demonstrates that these 
are also very practical tools for studying quivers and dimer models.
 
A similar beautiful insight of simplification by going to higher dimensions
appeared in De Bruijn's construction \cite{B}
of Penrose tilings and developed into the well-known projection method in the theory of quasi-crystals; e.g. \cite{Sen}. We apply the same method to 
construct periodic rhombus tilings of the plane, 
which the physicists call \emph{brane tilings} and \emph{dimer models}.

\begin{figure}[h]
\begin{picture}(350,100)(0,-5)
\put(-15,10){
\begin{picture}(150,100)(0,15)
\setlength\epsfxsize{4cm}
\epsfbox{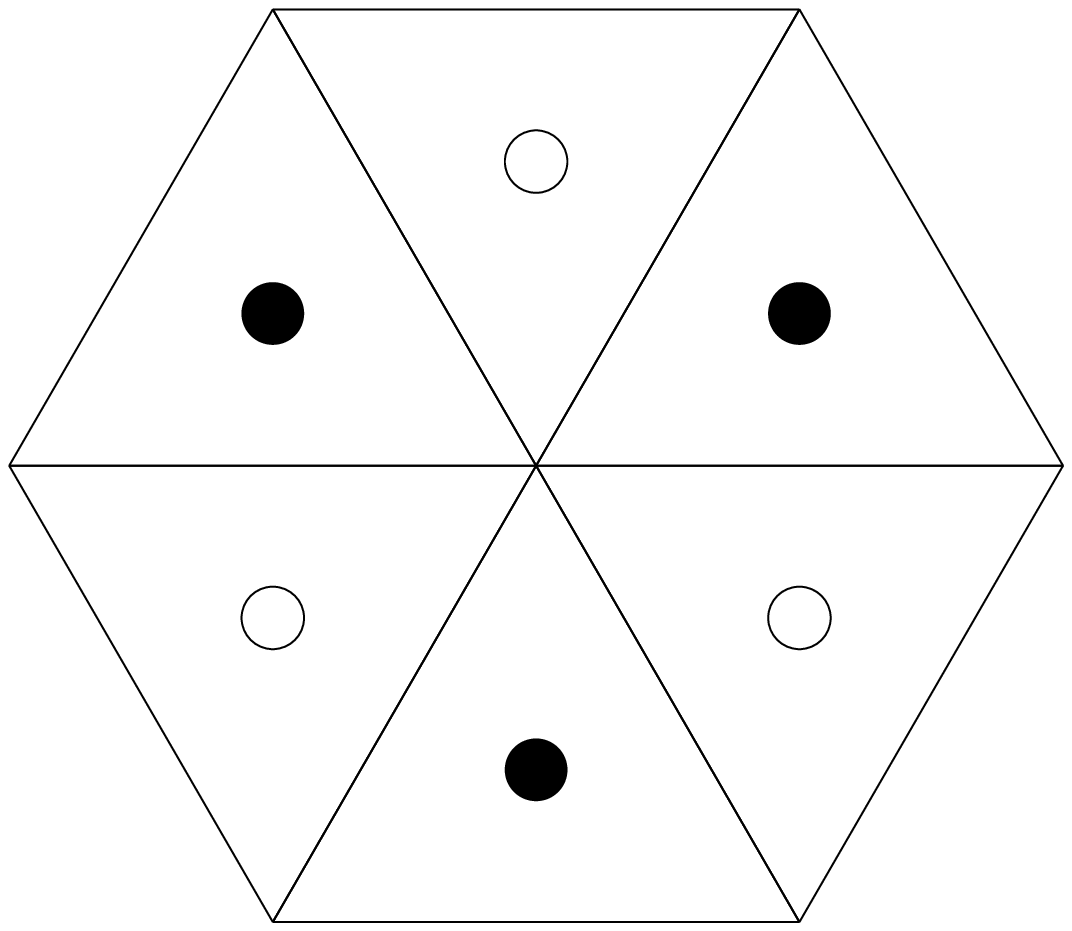}
\end{picture}
}
\put(120,10){
\begin{picture}(150,100)(0,15)
\setlength\epsfxsize{4cm}
\epsfbox{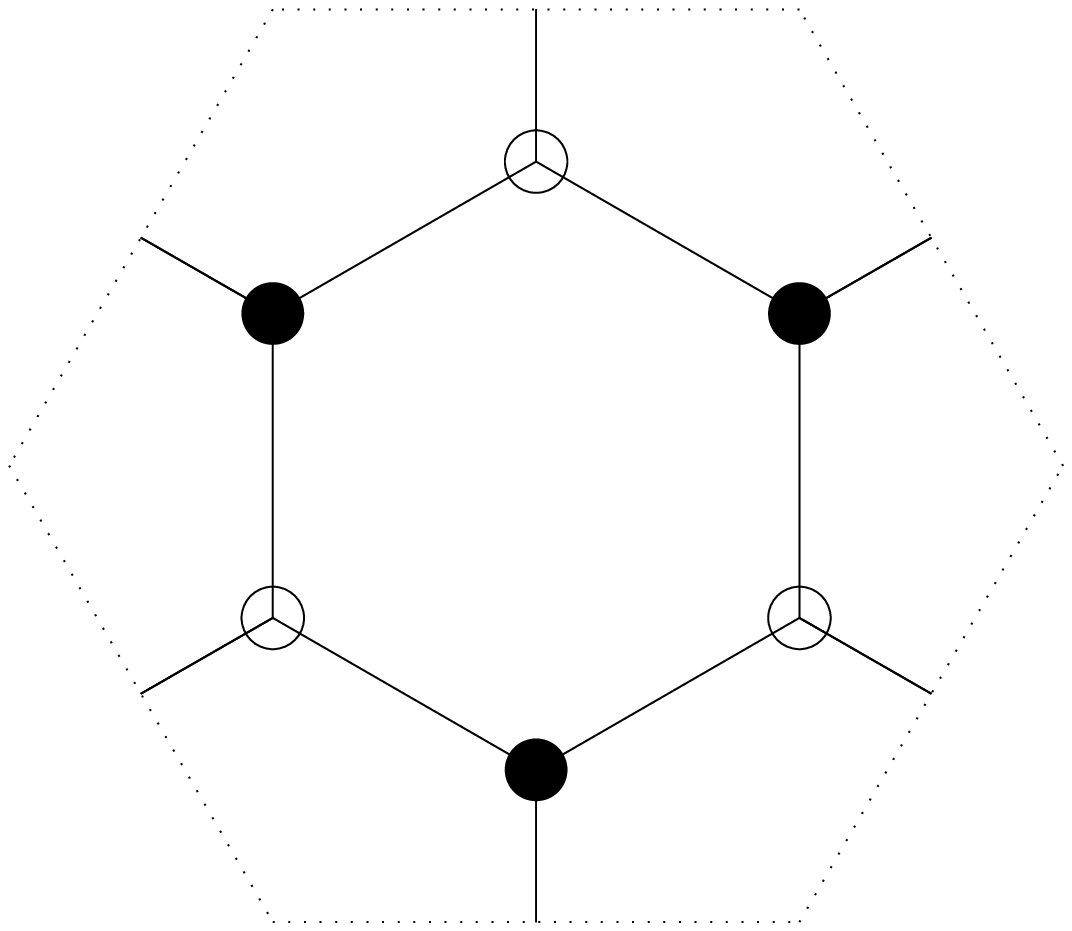}
\end{picture}
}
\put(255,10){
\begin{picture}(150,100)(0,15)
\setlength\epsfxsize{4cm}
\epsfbox{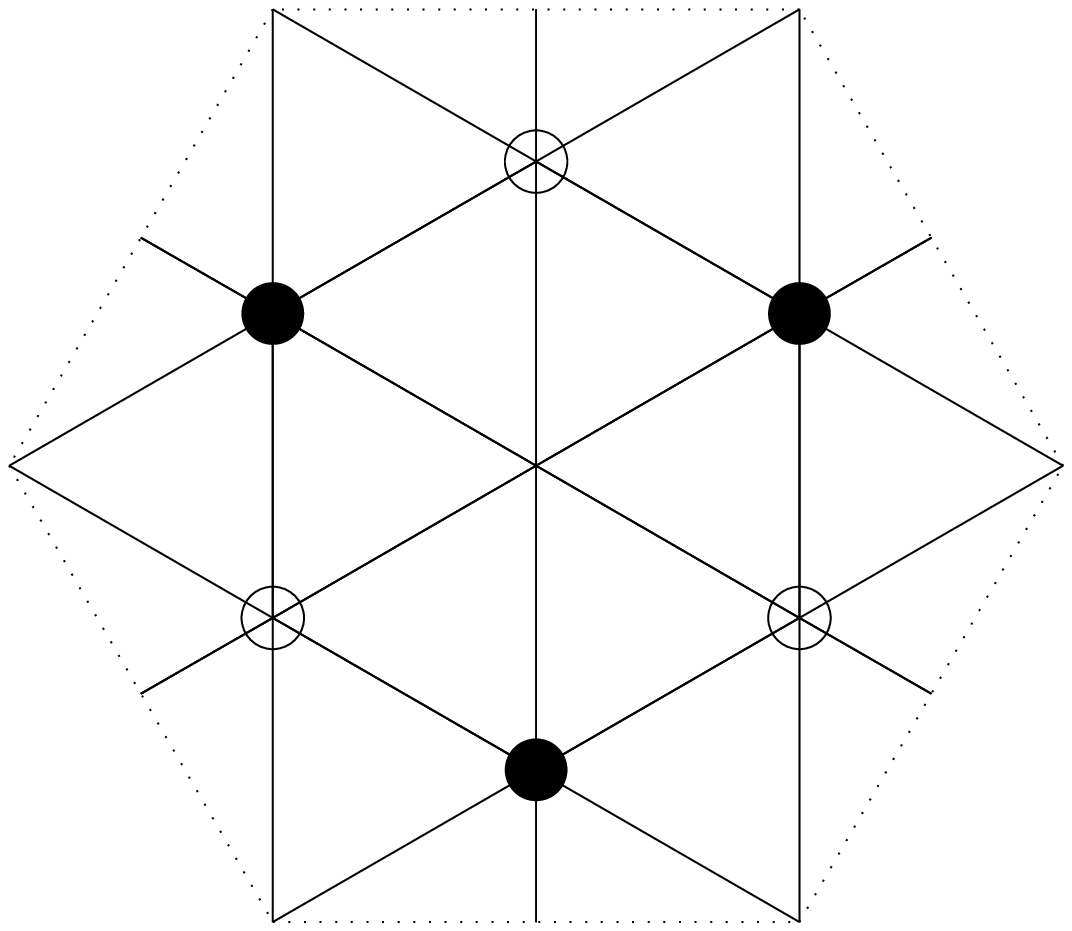}
\end{picture}
}
\end{picture}
\caption{\label{fig:P2combi}
\textit{Three versions of the dessin d'enfants 
for $\PP^2$  (=
case  $B_{1}$ in Fig. \ref{fig: models}). The surface is obtained by identifying opposite sides of the hexagon.}}
\end{figure}

From a geometric perspective this paper deals with 
embeddings of quivers into compact oriented surfaces without boundary.
More specifically, the initial combinatorial data for
a quiver $\sQ$ are two finite sets $E$ (arrows) and $V$ (nodes) 
and two maps $s,t:E\rightarrow V$ (source and target). 
Embedding $\sQ$ into a compact oriented surface without boundary $\RS$ 
means that $V$ becomes a subset of $\RS$ and an arrow $e\in E$ 
becomes a path $p_e$ in $\RS$ from the point $s(e)$ to the point $t(e)$.
It is required that
the boundary of every \emph{face of} $(\RS,\sQ)$ -- i.e. connected
component of $\RS\setminus\bigcup_{e\in E} p_e$ -- is
formed by a sequence of paths $(p_{e_1},\ldots,p_{e_n})$ with 
$p_{e_i}\cap p_{e_{i+1}}=t(e_i)=s(e_{i+1})$ for $1\leq i\leq n$; $e_{n+1}=e_1$.
It is evident from this requirement that
the face boundaries receive an orientation from $\sQ$. When this is combined with the orientation on $\RS$ faces lie either on the positive 
or on the negative side of their boundary.
In pictures we mark faces which lie on the positive (respectively negative) side of their boundary with $\bullet$ (resp. $\circ$). It is obvious that
adjacent faces get different colors. The dual graph of $\sQ$ w.r.t. $\RS$
is thus a bi-partite graph $\Gg$, embedded in $\RS$: the nodes of $\Gg$ correspond with the faces of $(\RS,\sQ)$; the edges of $\Gg$ connect nodes coming from adjacent faces
and correspond bijectively with the edges of $\sQ$. 
A pair $(\RS,\sQ)$ consisting of a compact oriented surface without boundary and an oriented graph embedded in it is called a
\emph{dessin d'enfants} or just \emph{dessin}, in one of the various equivalent definitions
of ``dessin (d'enfants)''; see \cite{Des,LZ,SV}. Other definitions refer to 
$(\RS,\Gg)$ as dessin d'enfants. In case $\RS$ has genus $1$ one often
calls $\Gg$ and its lifting to the plane (i.e. the universal covering of the torus $\RS$) a \emph{dimer model} or
\emph{brane tiling}; see for instance \cite{FHKV,FHKVW,HHV,HK,HV,Ky}. 
Yet another presentation of the same structure gives the dessin as the triangulation of $\RS$ with vertex set the union of the nodes of $\sQ$ and $\Gg$; the triangles are given by the triplets of vertices consisting of
two nodes of $\Gg$ connected by an edge $e^*$ of $\Gg$ and one node $v$ of 
$\sQ$ incident to the edge $e$ of $\sQ$ which is dual to $e^*$.
Figure \ref{fig:P2combi} shows these three get-ups of the dessin d'enfants of the 
$\PP^2$-quiver (case $B_1$ in Figure \ref{fig: models})

From an algebraic perspective this paper deals with 
\emph{superpotentials for quivers}. There is a simple equivalence between
the geometric and algebraic perspectives: the superpotential is
a convenient way of writing the list of oriented boundaries
$(p_{e_1},\ldots,p_{e_n})$ of the faces of 
$(\RS,\sQ)$ with $\pm$-signs indicating how the orientation matches with that of $\RS$. It can also be read as the instruction for building $\RS$ by glueing polygons. For the dessin d'enfants in Figure \ref{fig:P2combi} the
superpotential (for some numbering of the arrows of $\sQ$) is
$$
X_1X_2X_3+X_4X_5X_6+X_7X_8X_9-X_2X_8X_5-X_3X_9X_6-X_4X_1X_7\,.
$$

Another algebraic perspective, equivalent to the previous one, is given by what we want to call
\emph{the bi-adjacency matrix of the dessin d'enfants $(\RS,\sQ)$ with weight 
$\wt$}. It is defined as follows.
The \emph{weight} is a function $\wt:E\rightarrow\CC$. 
The dessins in the present paper have as many black faces as white faces
and
every edge $e\in E$ lies in the boundary of a unique black face, denoted 
$\vb(e)$, and a unique white face, denoted $\vw(e)$.
The bi-adjacency matrix of $(\RS,\sQ)$ with weight $\wt$ is a 
square matrix $\KK^\wt$ with rows 
corresponding with the black faces, columns corresponding with the white faces
and entries in the polynomial ring $\CC[u_v\:|\:v\in V]$
which has one variable for every node $v$ of the quiver $\sQ$:
\emph{the $(\vb,\vw)$-entry of $\KK^\wt$ is} 
\begin{equation}\label{eq:bi-adjacency matrix}
\KK^\wt_{\vb,\vw}\,=\,\sum_{e\in E:\,\vb(e)=\vb,\,\vw(e)=\vw}
\wt (e)\,u_{s(e)}u_{t(e)}\,.
\end{equation}
The bi-adjacency matrix is in fact the \emph{Kasteleyn matrix} of a twist
of the dimer model $(\RS,\Gg)$; see Section \ref{sec:Kasteleyn superpot}.
For the dessin in Figure \ref{fig:P2combi} the
bi-adjacency matrix (for some numbering of the nodes, edges and faces) is
$$
\left[\begin{array}{ccc}
\wt_1u_1u_3&\wt_2u_1u_2&\wt_3u_2u_3\\
\wt_4u_2u_3&\wt_5u_1u_3&\wt_6u_1u_2\\
\wt_7u_1u_2&\wt_8u_2u_3&\wt_9u_1u_3
\end{array}\right]\,.
$$

In Section \ref{sec:quiver} we describe how the quivers for which we can
solve the embedding problem, are associated with certain
rank $2$ subgroups $\LL$ of $\ZZ^N$; here $N$ is the number of nodes 
of the quiver. Such subgroups $\LL\subset\ZZ^N$ are the foundations for
the theory of Gelfand, Kapranov and Zelevinsky.
In Section \ref{sec:Lsecondaries} we describe the 
\emph{secondary fan} and \emph{the secondary polygon} associated
with $\LL\subset\ZZ^N$. These important structures in GKZ-theory
surprisingly turn out to be also quite relevant for the dimer models.
We show for instance in Theorem \ref{thm:Kasteleyn and secondary polygon} that the secondary polygon is the \emph{Newton polygon} of the determinant of the bi-adjacency matrix with non-zero weights.
In Section \ref{sec:A determinant}  we report the observation, based on examples, that for
the \emph{critical weight}
\begin{equation}\label{eq:critical weight}
\crit: E\rightarrow\ZZ_{>0}\,,\qquad \crit(e)=\sharp\{e'\in E\,|\,
s(e')=s(e),\,t(e')=t(e)\}
\end{equation}
the determinant of the bi-adjacency matrix $\KK^\crit$,
after a simple transformation, becomes equal to the 
\emph{principal $\cA$-determinant} of Gelfand, Kapranov, Zelevinsky \cite{gkz4}.
In Sections \ref{sec:hyper} and \ref{sec:solutions}
we recall some topics in the theory of GKZ-hypergeometric systems of differential equations.

Section \ref{sec:dessin} presents the algorithm to solve the quiver embedding problem. It is a combination of the \emph{Fast Inverse Algorithm} of
\cite{HV} and the \emph{untwisting procedure} of \cite{FHKV}.
In the cited papers, however, methods are mainly presented via visual 
inspection of pictures in some concrete examples. So an extrapolation
to general situations was needed. The first impression that the Fast Inverse Algorithm of \cite{HV} is more or less De Bruijn's construction 
\cite{B} of Penrose tilings did not quite yield the sought after embedding of the quiver; the untwisting procedure of \cite{FHKV} is also needed.
The algorithm became a smoothly operating algebraic-combinatorial tool
by consistently working from the philosophy that things look simpler
from a higher dimensional viewpoint (with group action). For doing
the computer experiments behind this paper
I implemented the algorithm in \textsc{matlab}. The reader can find in Section \ref{sec:dessin}
sufficient details for making a computer version of the algorithm.

Sections \ref{sec:perfect sec fan}--\ref{sec:A determinant} 
demonstrate that consistently working from the higher dimensional 
viewpoint leads to new insights in the dimer technology tools \emph{perfect matchings, Kasteleyn matrix and its determinant}
and shows their close relation with the GKZ tools 
\emph{secondary polytope and principal $\cA$-determinant}.

Besides problems such as proving that the algorithm of Section \ref{sec:dessin}
yields at least one superpotential for every quiver which satisfies the
conditions in Theorem \ref{quivergrassmann}, or proving
Conjecture \ref{conj:critical weight} our work raises some interesting questions like:

\textbf{Q1.} The Calabi-Yau singularities in the background of this
work are constructed from toric diagrams, which here are interpreted 
as secondary polygons. In \cite{MS} \S 6 the singularity is 
obtained from the toric diagram by a symplectic quotient construction. 
In \cite{HHV} \S 2 the singularity is given as a toric variety
constructed from a fan with one maximal cone, namely the cone over
the toric diagram (polygon). For this the polygon is put into
the plane of points with first coordinate $1$ in $\RR^3$.
In the work of Gelfand, Kapranov and Zelevinsky the secondary polygon
is put into a $2$-dimensional plane in $\RR^N$, which does not pass through $\nv$; here $N$ is the number of nodes of the quiver. One can perform the
standard toric variety construction with the $3$-dimensional cone over the
secondary polygon in $\RR^N$. 
\emph{Does this toric variety give the appriopriate view
in the GKZ philosophy (higher dimension compensated by group action) on
the singular Calabi-Yau $3$-space?} Since this toric variety
is a natural domain for GKZ hypergeometric functions (see \ref{toric sec fan})
one may wonder: \emph{Do hypergeometric functions provide new useful tools 
for investigating the geometry of $3$-dimensional Calabi-Yau singularities?}
Possibly positive indications are the fact that the dimension of the
solution space of the GKZ system of hypergeometric differential equations
equals the size of the bi-adjacency matrix $\KK^\crit$ (see Theorem
\ref{thm:solution dimension} and Corollary \ref{black volA})
and the observation (Conjecture \ref{conj:critical weight}) 
that the determinant of $\KK^\crit$ equals, up to a simple transformation, 
the principal $\cA$-determinant, which describes the 
singularities of the GKZ system.

\textbf{Q2.} The surface $\RS$ with embedded in it the pair of dual graphs
$\sQ$ and $\Gg$ has been constructed in a purely combinatorial topological
way from the superpotential. It is a general fact (see \cite{LZ} \S 1.2)
that there is then
a continuous map $f$ from $\RS$ to the $2$-sphere $S^2$ which is a ramified
covering with exactly three ramification points such that the set of vertices
of $\sQ$, the set of black vertices of $\Gg$ and the set of white vertices of $\Gg$ are the fibers over the three ramification points.
A highlight in the theory of dessins d'enfants is \emph{Belyi's
Theorem} (see \cite{LZ} Theorem 2.1.1). It states that in the above situation
the surface $\RS$ admits a model $\RS_K$ over a number field $K$ such that 
$f$ becomes a morphism $\RS_K\rightarrow \PP_K^1$ of algebraic curves over $K$ 
which is unramified outside $\{0,1,\infty\}$. The labeling of the ramification points can be taken such that $f^{-1}(\infty)$, $f^{-1}(0)$
and $f^{-1}(1)$ are the sets of white and black vertices of $\Gg$
and the vertices of $\sQ$, respectively.
The bipartite graph $\Gg$ in $\RS$ is then the inverse image of the 
negative real axis $[-\infty,0]$ in the Riemann sphere $\PP^1_\CC\,=\,S^2$ 
and the quiver $\sQ$ is the inverse image of the positively oriented unit 
circle $\{z\in\CC|\:|z|=1\}$. 
It is usually difficult to find explicit algebraic equations for $\RS_K$ and 
the \emph{Belyi function} $f$.\\
On the other hand, the authors of \cite{FHKV} write in footnote ${}^6$:
\emph{.... we have produced a dimer model that is defined on its own spectral
curve $\det\KK^\wt=0$.} Unfortunately their arguments are not (yet) sufficiently
refined to yield the weight $\wt$ that is to be used in this equation.
\\
It seems an interesting challenge to tackle the two problems simultaneously and look
for a weight $\wt$ and a Belyi function $f$ on the algebraic curve
with equation $\det\KK^\wt=0$ that realize $(\RS,\sQ,\Gg)$ as described above.

\section{Rank $2$ subgroups of $\ZZ^N$ and Quivers.}
\label{sec:quiver}

\ssnnl{ The most basic object}\label{introL}
in this paper is a
rank $2$ subgroup $\LL\subset\ZZ^N$
which is not contained in any of the standard 
coordinate hyperplanes of $\ZZ^N$ and is perpendicular to the vector
$(1,\ldots,1)$, with all coordinates $1$.
So, the elements of $\LL$ are vectors $(\ell_1,\ldots,\ell_N)\in\ZZ^N$
with $\ell_1+\ldots+\ell_N\,=\,0$ and for every $i\in\{1,\ldots,N\}$
there is an $(\ell_1,\ldots,\ell_N)\in\LL$ such that $\ell_i\neq 0$.

\

\ssnnl{ Notation.}\label{nota:eJ}
Throughout this paper $\ve_1,\ldots,\ve_N$ is the standard basis of $\ZZ^N$
and $\sJ$ denotes the matrix 
\footnotesize$\left(\hspace{-.8em}
\begin{array}{rr}0&1\\ -1&0\end{array}\hspace{-.5em}\right)$.
\normalsize

\

\ssnnl{ Definition/construction.}\label{Plucker coordinates}
Taking the second exterior powers one finds an inclusion
$\bigwedge^2\LL\hookrightarrow\bigwedge^2\ZZ^N$.
The group $\bigwedge^2\LL$ is a free $\ZZ$-module of rank $1$.
After fixing the orientation on $\LL$,
i.e. choosing one of the two possible isomorphisms $\bigwedge^2\LL\simeq\ZZ$, one gets an inclusion
$\ZZ\hookrightarrow\bigwedge^2\ZZ^N$.
The coordinates of $1\in\ZZ$ with respect to the standard basis $\{\ve_i\wedge\ve_j\}_{1\leq i<j\leq N}$ of $\bigwedge^2\ZZ^N$ are called the
\emph{Pl\"{u}cker coordinates of $\LL\subset\ZZ^N$}.

Dual to the above inclusion one finds the map
$\bigwedge^2\ZZ^N\rightarrow\ZZ$,
i.e. an anti-symmetric bilinear form on $\ZZ^N$. In view of its relation
with the Pl\"{u}cker coordinates we call this the 
\emph{Pl\"{u}cker form of $\LL$}.

To make this look more explicit we take a basis for $\LL$ compatible with the
chosen orientation. One can represent the two basis vectors as the rows
of a $2\times N$-matrix $B$. Let $\vb_1,\ldots,\vb_N$ be the columns of this matrix. Then the Pl\"{u}cker coordinates are 
$\left(\det(\vb_i,\vb_j)\right)_{1\leq i<j\leq N}$ and the (anti-symmetric) 
matrix for the Pl\"{u}cker form is $(\det(\vb_i,\vb_j))_{1\leq i,j\leq N}$.
The Pl\"{u}cker coordinates and the Pl\"{u}cker form do, of course, not change
if one takes another basis for $\LL$ with the same orientation.
Note that the matrix for the Pl\"{u}cker form has rank $2$:
\begin{equation}\label{eq:Plucker form}
\rank\left((\det(\vb_i,\vb_j))_{1\leq i,j\leq N}\right)\,=\,2\,.
\end{equation}

\

\ssnnl{ Definition.}\label{quiver basics} 
A \emph{quiver} $\sQ$ is a finite directed graph. It can be given combinatorially
by two finite sets $E$ and $V$ and two maps $s,t:E\rightarrow V$;
in short hand notation $\sQ=(E,V,s,t)$. 
Each element of $V$ is a node of the graph 
and an element $e\in E$ is an arrow from node $s(e)$ to node $t(e)$.
The \emph{adjacency matrix} of the quiver is the matrix
$(q_{ij})_{i,j\in V}$ with entry $q_{ij}$ equal to the number of arrows from
node $i$ to node $j$. 
If the quiver has no directed loops of length $\leq 2$
(i.e. $q_{ij}q_{ji}=0$ for all $i,j\in V$), it is faithfully described
by the \emph{anti-symmetrized adjacency matrix} $Q=(q_{ij}-q_{ji})_{i,j\in V}$.
If the $(i,j)$-entry of $Q$ is positive it
equals the number of arrows from node $i$ to node $j$. If it is negative it
is the number of arrows from $j$ to $i$. 

\

\ssnnl{ Definition.}\label{pluckerquiver}
 We call the quiver without directed loops of length 
$\leq 2$ with anti-symmetrized adjacency matrix equal to the
matrix of the Pl\"{u}cker form of $\LL$ the 
\emph{Pl\"{u}cker quiver of $\LL$}.

\

\ssnnl{ Examples} of groups $\LL\subset\ZZ^N$ and their Pl\"{u}cker quivers are shown in Figures \ref{fig: models} and \ref{fig:dp3 models}.

\

\ssnnl{ }\label{kirchhof}
Some easily visible geometric properties of a quiver without directed loops of length $\leq 2$
are equivalent to easily seen algebraic properties of its anti-symmetrized adjacency matrix $Q$. For instance, the graph
has no isolated nodes if and only if no column of $Q$ is zero. 
Also, in every node of the quiver the number of incoming 
arrows equals the number of outgoing arrows precisely if in each row of $Q$ 
the sum of the entries is $0$.
For Pl\"{u}cker quivers these properties correspond to the conditions that
$\LL$ is not contained in any of the standard 
coordinate hyperplanes of $\ZZ^N$ and is perpendicular to the vector
$(1,\ldots,1)$.
One property of a Pl\"{u}cker quiver which one can not easily read off from the graph is $\rank Q\,=\,2$.

In the remainder of this section we show that these properties characterize  
Pl\"{u}cker quivers of rank $2$ subgroups $\LL\subset\ZZ^N$ as in \ref{introL}.

\ 

\ssnnl{ Proposition.}\label{prop:rank 2 anti-sym}
\textit{An $N\times N$-matrix $C$ with entries in $\ZZ$ is anti-symmetric
and $\rank C\,=\,2$ if and only if there is a $2\times N$-matrix $B$ 
with entries in $\ZZ$ such that $\rank B\,=\,2$ and
$C\,=\,B^t\,\sJ\, B$.}

\proof
The ``if''-statement is trivial.
So let us consider the ``only if'' and assume that $C$ is anti-symmetric
and $\rank C\,=\,2$. Let $d$ denote the greatest common divisor of the entries of
$C$ and $C'=\frac{1}{d}C$.
Choose an $N\times 2$-matrix $D$ whose columns form a $\ZZ$-basis for the  column space of $C'$. The equality of column spaces of $C'$ and $D$ means that there are an $N\times 2$-matrix $E$ and
a $2\times N$-matrix $F$, both with entries in $\ZZ$ and of rank $2$,
such that
$D=C'E$ and $C'=DF$. Then
$
C'\,=\,DF\,=\,C'EF\,=\,-C'^t EF\,=\,-F^tD^t EF\,.
$ 
So $D^tE$ is an anti-symmetric
$2\times 2$-matrix; say
$D^tE=f\sJ$ with  $f\in\ZZ$. In fact $f=\pm 1$ since it divides all entries of $C'$. If $f=1$ we replace $F$ by $\sJ F$. We then always have $C'\,=\,F^t\sJ F$.
Let $G$ be any $2\times 2$-matrix with entries in $\ZZ$ and $\det G=d$. Let $B=GF$. Then 
$C=B^t\sJ B$ as wanted.
\qed

\

\ssnnl{}\label{quiver to L}
Multiplying in Proposition \ref{prop:rank 2 anti-sym} the matrix $B$
from the left by a matrix from $Sl_2(\ZZ)$ does not change the matrix
$B^t\sJ B$. So it is more natural to interpret $C$ as the matrix of the 
Pl\"{u}cker form of the $\ZZ$-row space of $B$. Note that it follows from the proof
of Proposition \ref{prop:rank 2 anti-sym} that this space is uniquely determined if the greatest common divisor of the entries of $C$ is $1$.
If on the other hand the greatest common divisor of the entries of $C$ is 
$d>1$ this space depends on the additional choice of a $2\times 2$-matrix with entries
in $\ZZ$ and determinant $d$.

\

Let us summarize the above discussion:

\

\ssnnl{ Theorem.}\label{quivergrassmann}
\textit{For every $\LL\subset\ZZ^N$ as in \ref{introL} one has its
Pl\"ucker quiver. This quiver has no isolated nodes
or directed loops of length $\leq 2$ and in every node the number of incoming arrows equals the number of outgoing arrows. Moreover the rank of its anti-symmetrized adjacency matrix is $2$.
Conversely, every quiver with these properties is the Pl\"ucker quiver
of some $\LL\subset\ZZ^N$ as in \ref{introL}.
The correspondence between such quivers and such $\LL\subset\ZZ^N$
is one-to-one for those quivers for which the greatest common divisor of the entries of the anti-symmetrized adjacency matrix is $1$
and those $\LL\subset\ZZ^N$ for which $\modquot{\ZZ^N}{\LL}$ has no torsion.}
\qed

\section{The secondary fan and polygon of $\LL\subset\ZZ^N$.}
\label{sec:Lsecondaries}
\ssnnl{}\label{intro L secondary}
For $\LL\subset\ZZ^N$ as in \ref{introL}
we now define the \emph{secondary fan} and the 
\emph{secondary polytope}. In \ref{A secondary} we will compare this
with the original definitions by Gelfand, Kapranov and Zelevinsky.
The term ``secondary'' refers to the fact that in their theory of 
hypergeometric systems another polytope appears first, which is therefore called the \emph{primary polytope}. Nonetheless
both the secondary fan and the secondary polytope can most conveniently and
directly be described using the lattice $\LL$. 
In case the rank
of $\LL$ is $2$, the constructions become particularly simple.

\

\ssnnl{ Definition.}\label{def:Lsecondary fan}
Let $\LL\subset\ZZ^N$ be as in \ref{introL}.
Let $\LL_\RR^\vee:=\mathrm{Hom} (\LL,\RR)$ denote the real dual space of $\LL$.
Let $\ve_1,\ldots,\ve_N$ be the standard basis of $\ZZ^N$. Let 
$\vb_i\in\LL_\RR^\vee$ be the image of $\ve_i$ under the map 
$\RR^N\rightarrow\LL_\RR^\vee$ dual to the inclusion 
$\LL\hookrightarrow\ZZ^N$. Here and henceforth we identify $\RR^N$ with the real dual space of $\ZZ^N$ by means of the standard dot product.

By definition, the \emph{secondary fan of $\LL$} is the following collection of cones in
$\LL_\RR^\vee$: the $0$-dimensional cone $\{\nv\}$, the $1$-dimensional cones
$\RR_{\geq 0}\vb_i$ ($i=1,\ldots,N$) and the $2$-dimensional cones
which are the closures of the connected components of
$\LL_\RR^\vee\setminus\bigcup_{i=1}^N\RR_{\geq 0}\vb_i$.

\

\ssnnl{ Definition.}\label{def:Lsecondary polytope}
With a $2$-dimensional cone $C$ in the secondary fan
one associates the set 
\begin{equation}\label{eq:LC}
L_C:=\:\left\{ \{i,j\}\subset\{1,\ldots,N\}\;|\;
C\subset (\RR_{\geq 0}\vb_i+\RR_{\geq 0}\vb_j)\;\right\}
\end{equation}
of $2$-element subsets of $\{1,\ldots,N\}$ and the vector
\begin{equation}\label{eq:psiC}
\psi_C\,:=\,\sum_{\{i,j\}\in L_C}|\det (\vb_i,\vb_j)|\,(\ve_i+\ve_j)\,.
\end{equation}
The \emph{secondary polytope} $\Sigma(\LL)$ of $\LL\subset\ZZ^N$ is then 
defined as
$$
\Sigma (\LL)\,:=\,\mathsf{convex\;hull}(\{\psi_C\:|\: C\:\textrm{$2$-dim cone
of secondary fan}\,\})\,.
$$

\

\ssnnl{ Remark.}\label{secondary interpretations}
In \S \ref{A secondary} we give an interpretation of the above $L_C$
as a triangulation of the primary polytope.
In Equation (\ref{eq:PCLC}) we associate with $L_C$ a perfect matching in a bipartite graph.

\

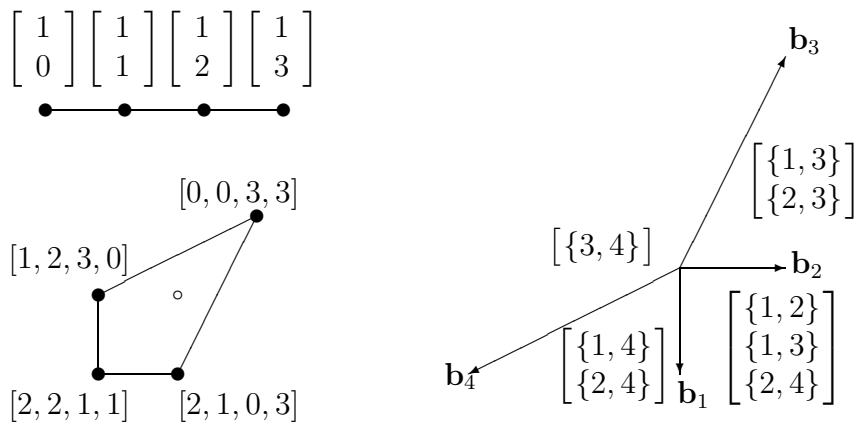
\begin{figure}[h]
\begin{picture}(300,170)(-10,-50)

\put(260,10){\begin{picture}(100,80)(0,0)
\put(0,0){\vector(0,-1){40}}
\put(0,0){\vector(1,0){40}}
\put(0,0){\vector(1,2){40}}
\put(0,0){\vector(-2,-1){80}}
\put(-1,-51){$\vb_1$}
\put(42,-2){$\vb_2$}
\put(41,83){$\vb_3$}
\put(-89,-44){$\vb_4$}
\put(15,-33){$\left[\!\!
\begin{array}{c}\{1,2\}\\ \{1,3\}\\ \{2,4\}\end{array}\!\!\right]$}
\put(25,30){$\left[\!\!
\begin{array}{c}\{1,3\}\\ \{2,3\}\end{array}\!\!\right]$}
\put(-47,-40){$\left[\!\!
\begin{array}{c}\{1,4\}\\ \{2,4\}\end{array}\!\!\right]$}
\put(-50,5){$\left[\!\!
\begin{array}{c}\{3,4\}\end{array}\!\!\right]$}
\end{picture}}

\put(-10,70){\begin{picture}(60,80)(-10,-10)
\put(5,10){$\left[\begin{array}{c}1\\ 0\end{array}\right]$}
\put(35,10){$\left[\begin{array}{c}1\\ 1\end{array}\right]$}
\put(65,10){$\left[\begin{array}{c}1\\ 2\end{array}\right]$}
\put(95,10){$\left[\begin{array}{c}1\\ 3\end{array}\right]$}
\put(20,-10){\line(1,0){90}}
\put(20,-10){\circle*{5}}
\put(50,-10){\circle*{5}}
\put(80,-10){\circle*{5}}
\put(110,-10){\circle*{5}}
\end{picture}}

\put(40,0){\begin{picture}(100,80)(0,0)
\put(0,0){\circle*{5}}
\put(0,-30){\circle*{5}}
\put(30,-30){\circle*{5}}
\put(60,30){\circle*{5}}
\put(30,0){\circle{3}}
\put(0,0){\line(0,-1){30}}
\put(0,-30){\line(1,0){30}}
\put(30,-30){\line(1,2){30}}
\put(60,30){\line(-2,-1){60}}
\put(-34,10){$[1,2,3,0]$}
\put(-34,-45){$[2,2,1,1]$}
\put(30,-45){$[2,1,0,3]$}
\put(30,35){$[0,0,3,3]$}
\end{picture}}
\end{picture}
\caption{\label{fig:secondaries}
\textit{The secondary fan (right) and the primary polytope (top left)
and the secondary polytope (bottom left) for
$\LL\,=\,\ZZ(0,1,1,-2)\oplus\ZZ(-1,0,2,-1)\subset\ZZ^4$.}}
\end{figure}

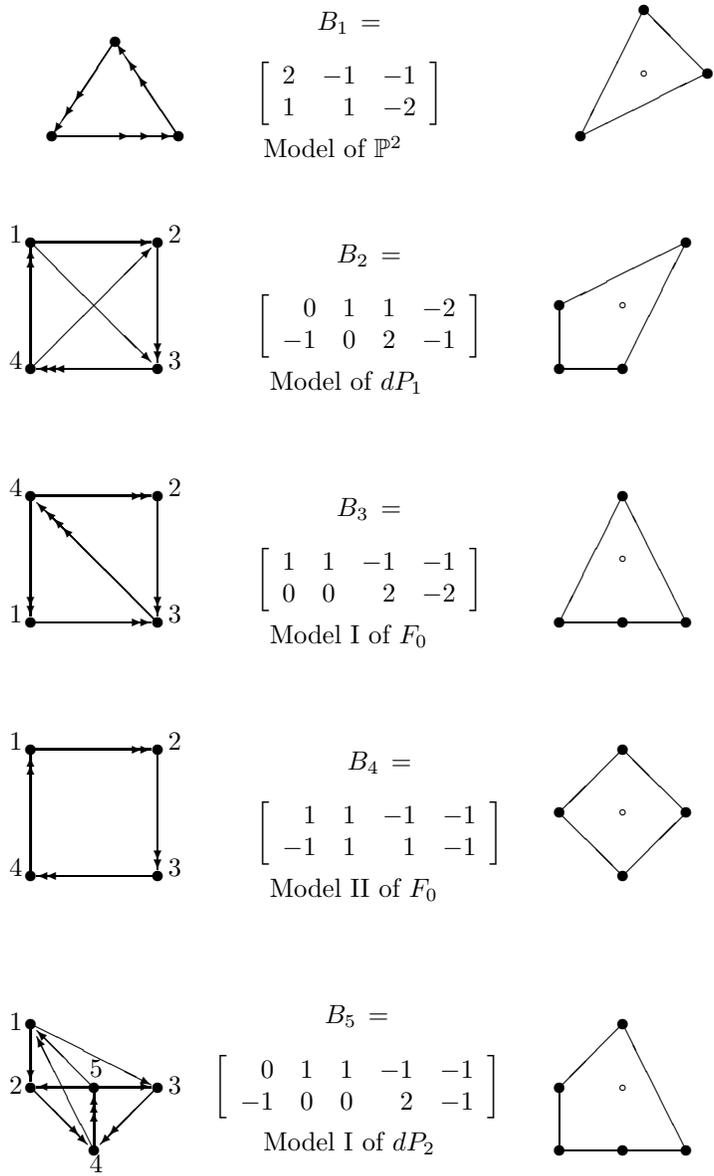
\begin{figure}[t]
\setlength{\unitlength}{0.8pt}
\begin{picture}(340,560)(-30,0)\footnotesize
\put(0,480){
\begin{picture}(300,80)(-20,-10)
\put(0,0){\begin{picture}(100,80)(0,0)
\put(0,0){\circle*{5}}
\put(60,0){\circle*{5}}
\put(30,45){\circle*{5}}
\put(0,0){\vector(1,0){57}}
\put(0,0){\vector(1,0){47}}
\put(0,0){\vector(1,0){37}}
\put(60,0){\vector(-2,3){29}}
\put(60,0){\vector(-2,3){24}}
\put(60,0){\vector(-2,3){19}}
\put(30,45){\vector(-2,-3){29}}
\put(30,45){\vector(-2,-3){24}}
\put(30,45){\vector(-2,-3){19}}
\end{picture}}

\put(70,0){\begin{picture}(60,80)(-10,-10)
\put(10,20){$\begin{array}{c} B_1\,=\\[2ex]
\left[
\begin{array}{rrr}
2&-1&-1\\ 1&1 & -2
\end{array}\right]
\end{array}$}
\put(20,-20){Model of $\PP^2$}
\end{picture}}

\put(250,0){\begin{picture}(100,80)(0,0)
\put(0,0){\circle*{5}}
\put(60,30){\circle*{5}}
\put(30,60){\circle*{5}}
\put(30,30){\circle{3}}

\put(0,0){\line(2,1){60}}
\put(60,30){\line(-1,1){30}}
\put(30,60){\line(-1,-2){30}}
\end{picture}}
\end{picture}}

\put(0,360){
\begin{picture}(300,80)(-10,-10)
\put(0,10){\begin{picture}(100,80)(0,0)
\put(0,0){\circle*{5}}
\put(0,60){\circle*{5}}
\put(60,0){\circle*{5}}
\put(60,60){\circle*{5}}
\put(-10,0){$4$}
\put(-10,60){$1$}
\put(65,0){$3$}
\put(65,60){$2$}
\put(0,0){\vector(0,1){57}}
\put(0,0){\vector(0,1){53}}
\put(0,0){\vector(1,1){57}}
\put(60,0){\vector(-1,0){57}}
\put(60,0){\vector(-1,0){53}}
\put(60,0){\vector(-1,0){49}}
\put(0,60){\vector(1,0){57}}
\put(0,60){\vector(1,-1){57}}
\put(60,60){\vector(0,-1){57}}
\put(60,60){\vector(0,-1){53}}
\end{picture}}

\put(70,10){\begin{picture}(60,80)(-10,-10)
\put(20,20){$\begin{array}{c} B_2\,=\\[2ex]
\left[\begin{array}{rrrr}
0&1&1&-2\\
-1&0&2&-1
\end{array}\right]\end{array}$}
\put(33,-20){Model of $dP_1$}
\end{picture}}

\put(250,40){\begin{picture}(100,80)(0,0)
\put(0,0){\circle*{5}}
\put(0,-30){\circle*{5}}
\put(30,-30){\circle*{5}}
\put(60,30){\circle*{5}}
\put(30,0){\circle{3}}

\put(0,0){\line(0,-1){30}}
\put(0,-30){\line(1,0){30}}
\put(30,-30){\line(1,2){30}}
\put(60,30){\line(-2,-1){60}}
\end{picture}}
\end{picture}}

\put(0,240){
\begin{picture}(300,80)(-10,-10)
\put(0,10){\begin{picture}(100,80)(0,0)
\put(0,0){\circle*{5}}
\put(0,60){\circle*{5}}
\put(60,0){\circle*{5}}
\put(60,60){\circle*{5}}
\put(-10,0){$1$}
\put(-10,60){$4$}
\put(65,0){$3$}
\put(65,60){$2$}
\put(0,60){\vector(1,0){57}}
\put(0,60){\vector(1,0){53}}
\put(0,60){\vector(0,-1){57}}
\put(0,60){\vector(0,-1){53}}
\put(0,0){\vector(1,0){57}}
\put(0,0){\vector(1,0){53}}
\put(60,0){\vector(-1,1){57}}
\put(60,0){\vector(-1,1){53}}
\put(60,0){\vector(-1,1){49}}
\put(60,0){\vector(-1,1){45}}
\put(60,60){\vector(0,-1){57}}
\put(60,60){\vector(0,-1){53}}
\end{picture}}

\put(70,10){\begin{picture}(60,80)(-10,-10)
\put(20,20){$\begin{array}{c} B_3\,=\\[2ex]
\left[\begin{array}{rrrrrr}
1&1&-1&-1\\
0&0&2&-2
\end{array}\right]\end{array}$}
\put(33,-20){Model I of $F_0$}
\end{picture}}

\put(250,10){\begin{picture}(100,80)(0,0)
\put(0,0){\circle*{5}}
\put(30,0){\circle*{5}}
\put(60,0){\circle*{5}}
\put(30,60){\circle*{5}}
\put(30,30){\circle{3}}

\put(0,0){\line(1,0){60}}
\put(60,0){\line(-1,2){30}}
\put(30,60){\line(-1,-2){30}}
\end{picture}}
\end{picture}}

\put(0,120){
\begin{picture}(300,80)(-10,-10)
\put(0,10){\begin{picture}(100,80)(0,0)
\put(0,0){\circle*{5}}
\put(0,60){\circle*{5}}
\put(60,0){\circle*{5}}
\put(60,60){\circle*{5}}
\put(-10,0){$4$}
\put(-10,60){$1$}
\put(65,0){$3$}
\put(65,60){$2$}
\put(0,0){\vector(0,1){57}}
\put(0,0){\vector(0,1){53}}
\put(60,0){\vector(-1,0){57}}
\put(60,0){\vector(-1,0){53}}
\put(0,60){\vector(1,0){57}}
\put(0,60){\vector(1,0){53}}
\put(60,60){\vector(0,-1){57}}
\put(60,60){\vector(0,-1){53}}
\end{picture}}

\put(70,10){\begin{picture}(60,80)(-10,-10)
\put(20,20){$\begin{array}{c} B_4\,=\\[2ex]
\left[\begin{array}{rrrrrr}
1&1&-1&-1\\
-1&1&1&-1
\end{array}\right]\end{array}$}
\put(33,-20){Model II of $F_0$}
\end{picture}}

\put(250,40){\begin{picture}(100,80)(0,0)
\put(0,0){\circle*{5}}
\put(30,-30){\circle*{5}}
\put(60,0){\circle*{5}}
\put(30,30){\circle*{5}}
\put(30,0){\circle{3}}

\put(0,0){\line(1,-1){30}}
\put(30,-30){\line(1,1){30}}
\put(60,0){\line(-1,1){30}}
\put(30,30){\line(-1,-1){30}}
\end{picture}}
\end{picture}}

\put(0,0){
\begin{picture}(300,80)(-10,-10)
\put(0,30){\begin{picture}(100,80)(0,0)
\put(0,0){\circle*{5}}
\put(30,0){\circle*{5}}
\put(30,-30){\circle*{5}}
\put(60,0){\circle*{5}}
\put(0,30){\circle*{5}}
\put(-10,27){$1$}
\put(-10,-3){$2$}
\put(28,-40){$4$}
\put(65,-3){$3$}
\put(28,5){$5$}
\put(0,0){\vector(1,-1){27}}
\put(0,0){\vector(1,-1){23}}
\put(30,-30){\vector(-1,2){27}}
\put(30,-30){\vector(0,1){27}}
\put(30,-30){\vector(0,1){23}}
\put(30,-30){\vector(0,1){19}}
\put(30,0){\vector(-1,1){27}}
\put(30,0){\vector(-1,0){27}}
\put(30,0){\vector(1,0){27}}
\put(0,30){\vector(2,-1){57}}
\put(0,30){\vector(0,-1){27}}
\put(60,0){\vector(-1,-1){27}}
\put(60,0){\vector(-1,-1){23}}
\end{picture}}

\put(70,10){\begin{picture}(60,80)(-10,-10)
\put(0,20){$\begin{array}{c} B_5\,=\\[2ex]
\left[\begin{array}{rrrrrr}
0&1&1&-1&-1\\
-1&0&0&2&-1
\end{array}\right]\end{array}$}
\put(30,-20){Model I of $dP_2$}
\end{picture}}

\put(250,30){\begin{picture}(100,80)(0,0)
\put(0,0){\circle*{5}}
\put(0,-30){\circle*{5}}
\put(30,-30){\circle*{5}}
\put(60,-30){\circle*{5}}
\put(30,30){\circle*{5}}
\put(30,0){\circle{3}}

\put(0,0){\line(0,-1){30}}
\put(0,-30){\line(1,0){30}}
\put(30,-30){\line(1,0){30}}
\put(60,-30){\line(-1,2){30}}
\put(30,30){\line(-1,-1){30}}

\end{picture}}
\end{picture}}
\setlength{\unitlength}{1pt}
\end{picture}
\caption{\label{fig: models}
\textit{Quivers from \cite{FFHH} Figures 10, 11, 4, 12 and 
corresponding polygons.}}
\end{figure}

\ssnnl{ Example.}\label{secondary fan dP1}
Figure \ref{fig:secondaries} shows the secondary fan with the lists $L_C$, 
the secondary polytope $\Sigma (\LL)$ with coordinates for the vertices
and the primary polytope
for $\LL\,=\,\ZZ(0,1,1,-2)\oplus\ZZ(-1,0,2,-1)\subset\ZZ^4$. This is case $B_2$ of Figure \ref{fig: models}; see also Example \ref{cubic roots} and 
\S \ref{A secondary}.

\

\ssnnl{ }\label{secondary polytope from quiver}
The geometric formulation of the construction of the secondary fan and
polytope is quite attractive. Nonetheless in practical algorithms
everything can be most easily obtained from the 
matrix of the Pl\"{u}cker form of $\LL\subset\ZZ^N$. 
Indeed, a $2$-dimensional cone in the secondary fan is bounded by
two half-lines $\RR_{\geq 0}\vb_i$ and $\RR_{\geq 0}\vb_j$
with the property that $\det(\vb_i,\vb_j)>0$ and
$\det(\vb_i,\vb_k)\det(\vb_j,\vb_k)\geq 0$ for $k=1,\ldots,N$.
Thus to find the secondary fan one just needs to find all pairs $\vb_i,\vb_j$ with these properties.

For a cone $C$ bounded by $\RR_{\geq 0}\vb_i$ and $\RR_{\geq 0}\vb_j$
the list $L_C$ is then 
$$
L_C\,=\,\{\,\{k,l\}\:|\: \det(\vb_k,\vb_i)\geq 0\;\textrm{and}\;
\det(\vb_j,\vb_l)\geq 0\,\}\,.
$$

Now
consider two adjacent $2$-dimensional cones in the secondary fan,
say $C$ and $C'$. Let $\vb_{i_1},\ldots,\vb_{i_s}$ be those vectors from the set
$\{\vb_1,\ldots,\vb_N\}$ that lie on the half-line $C\cap C'$. Without loss of generality we may assume $C$ lies to the right of $\vb_{i_1}$ and $C'$ to the left. Then
\begin{eqnarray*}
L_C\setminus(L_C\cap L_{C'})&=&
\{\,\{k,l\}\:|\: l\in\{i_1,\ldots,i_s\}\,,\;\det(\vb_l,\vb_k)<0\:\}\,,
\\
L_{C'}\setminus(L_C\cap L_{C'})&=&
\{\,\{k,l\}\:|\: l\in\{i_1,\ldots,i_s\}\,,\;\det(\vb_l,\vb_k)>0\:\}\,.
\end{eqnarray*}
From this we see, using $\sum_{k=1}^N\det(\vb_i,\vb_k)=0$ for all $i$, that
\begin{equation}\label{eq:secondary edge}
\psi_{C'}-\psi_{C}\,=\,\sum_{r=1}^s\sum_{k=1}^N\det(\vb_{i_r},\vb_k) \ve_k\,.
\end{equation}
In other words, $\psi_{C'}-\psi_C$ is the sum of rows $i_1,\ldots,i_s$
of the matrix of the Pl\"ucker form.

\

\ssnnl{}\label{explict secondary polygon}
To make this even more explicit and simple looking we take a basis for 
$\LL$ compatible with the
chosen orientation. One can represent the two basis vectors as the rows
of a $2\times N$-matrix $B$. Let $\vb_1,\ldots,\vb_N$ be the columns of this matrix. The elements of $\LL$ should now be written as row vectors with two
components and the embedding $\LL\hookrightarrow\ZZ^N$ is given by
$\vv\mapsto \sum_{k=1}^N (\vv\cdot\vb_k)\ve_k$. 
Since $\det(\vb_i,\vb_k)=\vb_{i}^t\sJ\vb_k$,
we can reformulate Equation (\ref{eq:secondary edge}) as
\begin{equation}\label{eq:secondary edge bis}
\psi_{C'}-\psi_C\,=\,\textrm{image of $
\sum_{r=1}^s\vb_{i_r}^t\sJ$ under the embedding $\LL\hookrightarrow\ZZ^N$\,.}
\end{equation}

In order to obtain the simplest formulation for the construction we order
the vectors $\vb_1,\ldots,\vb_N$ 
so that the points $\vp_k=\sum_{i=1}^k \vb_i^t$ for $k=1,\ldots,N$ lie ordered counterclockwise on the 
boundary of the polygon 
\begin{equation}\label{eq:easy secondary}
\Delta\,=\,\mathsf{convex\;hull}\{\vp_1,\ldots,\vp_N\}\,.
\end{equation}

\emph{Then the secondary polygon $\Sigma(\LL)$ is obtained by
first rotating $\Delta$ clockwise over $90^\circ$, next embedding it along with $\LL$ into $\ZZ^N$ and finally translating it over the vector
$\psi_C$, where $C$ is the cone in the secondary fan with left hand boundary
$\RR_{\geq 0}\vb_1$.} 

Note that while the secondary polytope $\Sigma (\LL)$ depends only on the 
embedding $\LL\hookrightarrow\ZZ^N$, the polygon $\Delta$ usually changes when
one puts another vector $\vb_i$ in first position by a cyclic permutation 
or when one multiplies the vectors $\vb_1,\ldots,\vb_N$ by a matrix from $Sl_2(\ZZ)$.

\

\ssnnl{}\label{polygon to L}
For a converse to the above construction we start from a convex
polygon $\Delta$ in $\RR^2$ with vertices in $\ZZ^2$ and non-empty interior. 
Let $\partial\Delta$ denote its boundary. Next we choose a collection of points $\vp_1,\ldots,\vp_N$ in $\partial\Delta\cap\ZZ^2$ which includes all vertices of $\Delta$.
We number these points so that $\vp_i$ and $\vp_{i+1}$ are consecutive points
in the counter-clockwise orientation of $\partial\Delta$. Thinking of the elements 
of $\ZZ^2$ as row vectors we define column vectors
$\vb_1=\vp_1^t$ and $\vb_i=\vp_i^t-\vp_{i-1}^t$ for $i=2,\ldots,N$.
Finally we define $\LL$ to be the $\ZZ$-row space of the $2\times N$-matrix
$B$ with columns $\vb_1,\ldots,\vb_N$.

\emph{Thus every convex
polygon $\Delta$ in $\RR^2$ with vertices in $\ZZ^2$ and non-empty interior
can be viewed as the secondary polygon of some rank $2$ subgroup
$\LL\subset\ZZ^N$ satisfying the conditions in \ref{introL}.}

Note however that in general there are several possible choices for
the points $\vp_1,\ldots,\vp_N$ in $\partial\Delta\cap\ZZ^2$.
The minimal choice takes only the vertices of $\Delta$, while the maximal
choice takes all points of $\partial\Delta\cap\ZZ^2$.

\

\ssnnl{ Theorem.}\label{secondary area}
\textit{ The (Euclidean) area of the polygon $\Delta$ in
(\ref{eq:easy secondary}) is }
\begin{equation}\label{eq:areaD}
\mathrm{area}\,\Delta\,=\,
\textstyle{\frac{1}{2}}\sum_{1\leq i<j\leq N}|\det(\vb_i,\vb_j)|\,-\,
\sum_{\{i,j\}\in L_C}|\det (\vb_i,\vb_j)|
\end{equation}
\textit{for every $2$-dimensional cone $C$ in the secondary fan.}

\proof
After applying, if necessary, a cyclic permutation to $\vb_1,\ldots,\vb_N$ 
we may assume, without loss of generality, that
$\,C=\RR_{\geq 0}\vb_1+\RR_{\geq 0}\vb_N$.

From the triangulation $\Delta$ by the diagonals between the vertex
$\vp_N=\nv$ and the other vertices of $\Delta$ one sees that
$\mathrm{area}\,\Delta\,=\,
\frac{1}{2}\sum_{1\leq i<j\leq M}\det(\vb_i,\vb_j)$
where $M$ is such that $\vb_j\not\in\RR_{\geq 0}\vb_N$ for $1\leq j\leq M$
and $\vb_j\in\RR_{\geq 0}\vb_N$ for $M+1\leq j\leq N$.
The fact that $\sum_{i=1}^N\vb_i=\nv$, implies
$$
\textstyle{\sum_{1\leq i\leq M<j\leq N}\det(\vb_i,\vb_j)\,=\,
\det(\sum_{1\leq i\leq M}\vb_i,\sum_{M<j\leq N}\vb_j)\,=\,0\,.}
$$
Thus we see
$$
\mathrm{area}\,\Delta\,=\,
\textstyle{\frac{1}{2}}\sum_{1\leq i<j\leq N}\det(\vb_i,\vb_j)\,.
$$
Next note that for $1\leq i<j\leq N$ the inequality $\det(\vb_i,\vb_j)<0$ 
holds if and only 
if $C\subset\RR_{\geq 0}\vb_i+\RR_{\geq 0}\vb_j$, if and only if 
$\{i,j\}\in L_C$. Equality (\ref{eq:areaD}) now follows immediately.
\qed

\

\ssnnl{ Corollary.}\label{internal secondary}
\textit{ The number of lattice points in the interior of
the polygon $\Delta$ in (\ref{eq:easy secondary}) is }
$$
\textstyle{\frac{1}{2}}\left(\sum_{1\leq i<j\leq N}|\det(\vb_i,\vb_j)|\,-\,N
\right)-\sum_{\{i,j\}\in L_C}|\det (\vb_i,\vb_j)|\,+\,1\,.
$$
\proof
This follows from Theorem \ref{secondary area} in combination with \emph{Pick's Formula}
(\cite{F} p.113):
$\;
2\mathrm{area}\,\Delta\,=\,
2\sharp\left(\ZZ^2\cap\mathrm{interior}\,\Delta\right)
\,+\,\sharp\left(\ZZ^2\cap\mathrm{boundary}\,\Delta\right)\,-\,2
$
\qed

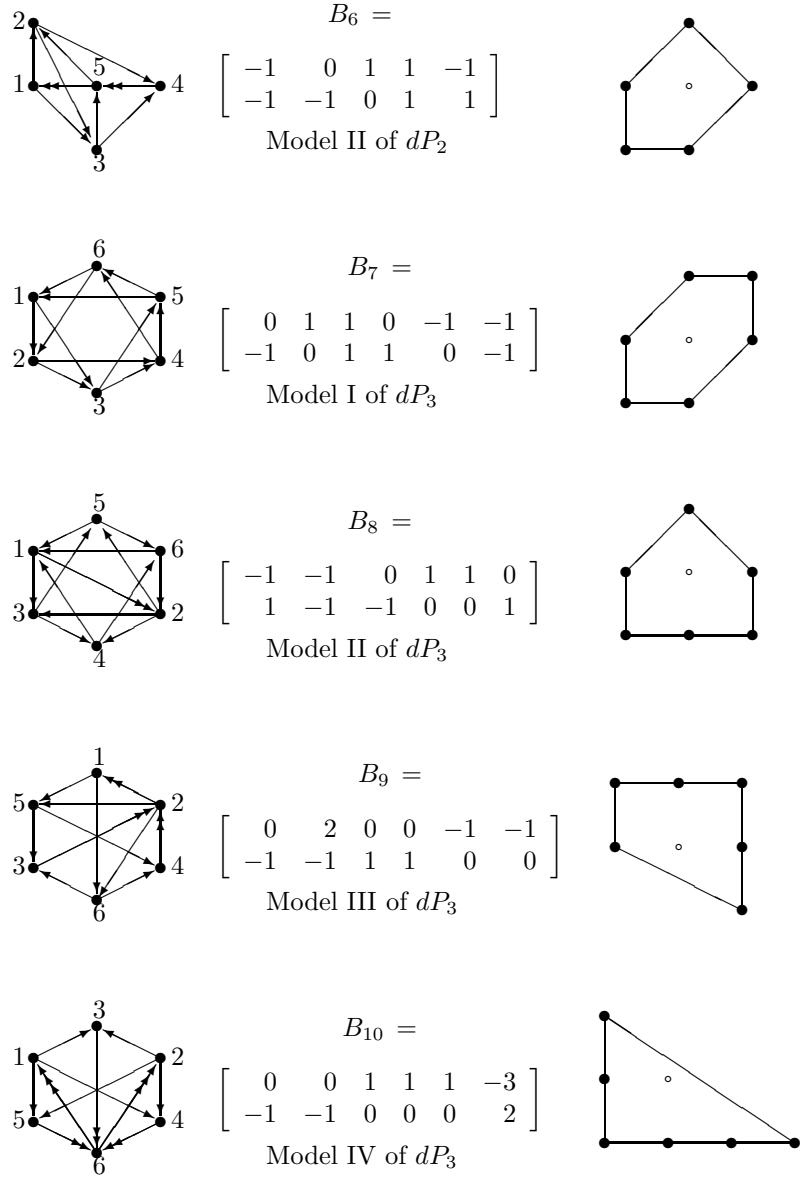
\begin{figure}[t]
\setlength{\unitlength}{0.8pt}
\begin{picture}(340,560)(-30,0)\footnotesize
\put(0,0){
\begin{picture}(300,80)(-10,-10)
\put(0,20){\begin{picture}(100,80)(0,0)
\put(0,0){\circle*{5}}
\put(0,30){\circle*{5}}
\put(30,-15){\circle*{5}}
\put(60,0){\circle*{5}}
\put(60,30){\circle*{5}}
\put(30,45){\circle*{5}}
\put(-10,27){$1$}
\put(-10,-3){$5$}
\put(28,-25){$6$}
\put(65,-3){$4$}
\put(65,27){$2$}
\put(28,49){$3$}
\put(0,30){\vector(0,-1){27}}
\put(0,30){\vector(2,1){27}}
\put(0,30){\vector(2,-1){57}}
\put(60,30){\vector(-2,1){27}}
\put(60,30){\vector(-2,-1){57}}
\put(60,30){\vector(0,-1){27}}
\put(30,45){\vector(0,-1){57}}
\put(30,45){\vector(0,-1){53}}
\put(0,0){\vector(2,-1){27}}
\put(0,0){\vector(2,-1){23}}
\put(60,0){\vector(-2,-1){27}}
\put(60,0){\vector(-2,-1){23}}
\put(30,-15){\vector(2,3){28}}
\put(30,-15){\vector(2,3){25}}
\put(30,-15){\vector(2,3){22}}
\put(30,-15){\vector(-2,3){28}}
\put(30,-15){\vector(-2,3){25}}
\put(30,-15){\vector(-2,3){22}}
\end{picture}}

\put(70,10){\begin{picture}(60,80)(-10,-10)
\put(0,20){$\begin{array}{c} B_{10}\,=\\[2ex]
\left[\begin{array}{rrrrrr}
0&0&1&1&1&-3\\
-1&-1&0&0&0&2
\end{array}\right]\end{array}$}
\put(30,-20){Model IV of $dP_3$}
\end{picture}}

\put(270,70){\begin{picture}(100,80)(0,0)
\put(0,0){\circle*{5}}
\put(0,-30){\circle*{5}}
\put(0,-60){\circle*{5}}
\put(30,-60){\circle*{5}}
\put(60,-60){\circle*{5}}
\put(90,-60){\circle*{5}}
\put(30,-30){\circle{3}}

\put(0,0){\line(0,-1){60}}
\put(0,-60){\line(1,0){90}}
\put(90,-60){\line(-3,2){90}}
\end{picture}}
\end{picture}}

\put(0,120){
\begin{picture}(300,80)(-10,-10)
\put(0,20){\begin{picture}(100,80)(0,0)
\put(0,0){\circle*{5}}
\put(0,30){\circle*{5}}
\put(30,-15){\circle*{5}}
\put(60,0){\circle*{5}}
\put(60,30){\circle*{5}}
\put(30,45){\circle*{5}}
\put(-10,27){$5$}
\put(-10,-3){$3$}
\put(28,-25){$6$}
\put(65,-3){$4$}
\put(65,27){$2$}
\put(28,49){$1$}
\put(0,30){\vector(0,-1){27}}
\put(0,30){\vector(2,-1){57}}
\put(0,0){\vector(2,1){57}}
\put(0,0){\vector(2,1){53}}
\put(30,-15){\vector(2,1){27}}
\put(30,-15){\vector(-2,1){27}}
\put(60,0){\vector(0,1){27}}
\put(60,0){\vector(0,1){23}}
\put(60,30){\vector(-2,-3){29}}
\put(60,30){\vector(-1,0){57}}
\put(60,30){\vector(-2,1){27}}
\put(60,30){\vector(-2,1){23}}
\put(30,45){\vector(-2,-1){27}}
\put(30,45){\vector(0,-1){57}}
\end{picture}}

\put(70,10){\begin{picture}(60,80)(-10,-10)
\put(0,20){$\begin{array}{c} B_9\,=\\[2ex]
\left[\begin{array}{rrrrrr}
0&2&0&0&-1&-1\\
-1&-1&1&1&0&0
\end{array}\right]\end{array}$}
\put(30,-20){Model III of $dP_3$}
\end{picture}}

\put(335,0){\begin{picture}(100,80)(0,0)
\put(0,0){\circle*{5}}
\put(0,30){\circle*{5}}
\put(0,60){\circle*{5}}
\put(-30,60){\circle*{5}}
\put(-60,60){\circle*{5}}
\put(-60,30){\circle*{5}}
\put(-30,30){\circle{3}}

\put(0,0){\line(0,1){60}}
\put(0,60){\line(-1,0){60}}
\put(-60,60){\line(0,-1){30}}
\put(-60,30){\line(2,-1){60}}
\end{picture}}
\end{picture}}

\put(0,240){
\begin{picture}(300,80)(-10,-10)
\put(0,20){\begin{picture}(100,80)(0,0)
\put(0,0){\circle*{5}}
\put(0,30){\circle*{5}}
\put(30,-15){\circle*{5}}
\put(60,0){\circle*{5}}
\put(60,30){\circle*{5}}
\put(30,45){\circle*{5}}
\put(-10,27){$1$}
\put(-10,-3){$3$}
\put(28,-25){$4$}
\put(65,-3){$2$}
\put(65,27){$6$}
\put(28,49){$5$}
\put(0,30){\vector(0,-1){27}}
\put(0,30){\vector(2,-1){57}}
\put(0,30){\vector(2,-1){53}}
\put(0,0){\vector(2,3){27}}
\put(0,0){\vector(2,-1){27}}
\put(30,-15){\vector(-2,3){27}}
\put(30,-15){\vector(2,3){27}}
\put(60,0){\vector(-2,3){27}}
\put(60,0){\vector(-1,0){57}}
\put(60,0){\vector(-2,-1){27}}
\put(60,30){\vector(-1,0){57}}
\put(60,30){\vector(0,-1){27}}
\put(30,45){\vector(2,-1){27}}
\put(30,45){\vector(-2,-1){27}}
\end{picture}}

\put(70,10){\begin{picture}(60,80)(-10,-10)
\put(0,20){$\begin{array}{c} B_8\,=\\[2ex]
\left[\begin{array}{rrrrrr}
-1&-1&0&1&1&0\\
1&-1&-1&0&0&1
\end{array}\right]\end{array}$}
\put(30,-20){Model II of $dP_3$}
\end{picture}}

\put(280,40){\begin{picture}(100,80)(0,0)
\put(0,0){\circle*{5}}
\put(0,-30){\circle*{5}}
\put(30,-30){\circle*{5}}
\put(60,-30){\circle*{5}}
\put(60,0){\circle*{5}}
\put(30,30){\circle*{5}}
\put(30,0){\circle{3}}

\put(0,0){\line(0,-1){30}}
\put(0,-30){\line(1,0){60}}
\put(60,-30){\line(0,1){30}}
\put(60,0){\line(-1,1){30}}
\put(30,30){\line(-1,-1){30}}
\end{picture}}
\end{picture}}

\put(0,360){
\begin{picture}(300,80)(-10,-10)
\put(0,20){\begin{picture}(100,80)(0,0)
\put(0,0){\circle*{5}}
\put(0,30){\circle*{5}}
\put(30,-15){\circle*{5}}
\put(60,0){\circle*{5}}
\put(60,30){\circle*{5}}
\put(30,45){\circle*{5}}
\put(-10,27){$1$}
\put(-10,-3){$2$}
\put(28,-25){$3$}
\put(65,-3){$4$}
\put(65,27){$5$}
\put(28,49){$6$}
\put(0,30){\vector(0,-1){27}}
\put(0,0){\vector(2,-1){27}}
\put(30,-15){\vector(2,1){27}}
\put(60,0){\vector(0,1){27}}
\put(60,30){\vector(-2,1){27}}
\put(30,45){\vector(-2,-1){27}}
\put(0,0){\vector(1,0){57}}
\put(60,0){\vector(-2,3){28}}
\put(30,45){\vector(-2,-3){28}}
\put(60,30){\vector(-1,0){57}}
\put(0,30){\vector(2,-3){28}}
\put(30,-15){\vector(2,3){28}}
\end{picture}}

\put(70,10){\begin{picture}(60,80)(-10,-10)
\put(0,20){$\begin{array}{c} B_7\,=\\[2ex]
\left[\begin{array}{rrrrrr}
0&1&1&0&-1&-1\\
-1&0&1&1&0&-1
\end{array}\right]\end{array}$}
\put(30,-20){Model I of $dP_3$}
\end{picture}}

\put(280,30){\begin{picture}(100,80)(0,0)
\put(0,0){\circle*{5}}
\put(0,-30){\circle*{5}}
\put(30,-30){\circle*{5}}
\put(60,0){\circle*{5}}
\put(60,30){\circle*{5}}
\put(30,30){\circle*{5}}
\put(30,0){\circle{3}}

\put(0,0){\line(0,-1){30}}
\put(0,-30){\line(1,0){30}}
\put(30,-30){\line(1,1){30}}
\put(60,0){\line(0,1){30}}
\put(60,30){\line(-1,0){30}}
\put(30,30){\line(-1,-1){30}}

\end{picture}}
\end{picture}}

\put(0,480){
\begin{picture}(300,80)(-10,-10)
\put(0,30){\begin{picture}(100,80)(0,0)
\put(0,0){\circle*{5}}
\put(30,0){\circle*{5}}
\put(30,-30){\circle*{5}}
\put(60,0){\circle*{5}}
\put(0,30){\circle*{5}}
\put(-10,27){$2$}
\put(-10,-3){$1$}
\put(28,-40){$3$}
\put(65,-3){$4$}
\put(28,5){$5$}
\put(0,0){\vector(0,1){27}}
\put(0,0){\vector(1,-1){27}}
\put(30,-30){\vector(1,1){27}}
\put(30,-30){\vector(0,1){27}}
\put(30,0){\vector(-1,1){27}}
\put(30,0){\vector(-1,0){23}}
\put(30,0){\vector(-1,0){27}}
\put(0,30){\vector(2,-1){57}}
\put(0,30){\vector(1,-2){27}}
\put(60,0){\vector(-1,0){27}}
\put(60,0){\vector(-1,0){23}}
\end{picture}}

\put(70,10){\begin{picture}(60,80)(-10,-10)
\put(0,20){$\begin{array}{c} B_6\,=\\[2ex]
\left[\begin{array}{rrrrrr}
-1&0&1&1&-1\\
-1&-1&0&1&1
\end{array}\right]\end{array}$}
\put(30,-20){Model II of $dP_2$}
\end{picture}}

\put(280,30){\begin{picture}(100,80)(0,0)
\put(0,0){\circle*{5}}
\put(0,-30){\circle*{5}}
\put(30,-30){\circle*{5}}
\put(60,0){\circle*{5}}
\put(30,30){\circle*{5}}
\put(30,0){\circle{3}}

\put(0,0){\line(0,-1){30}}
\put(0,-30){\line(1,0){30}}
\put(30,-30){\line(1,1){30}}
\put(60,0){\line(-1,1){30}}
\put(30,30){\line(-1,-1){30}}

\end{picture}}
\end{picture}}
\setlength{\unitlength}{1pt}
\end{picture}
\caption{\label{fig:dp3 models}
\textit{Quivers from \cite{FFHH} Figures 12, 9 and 
corresponding polygons.}}
\end{figure}

\clearpage
\begin{figure}[t]
\setlength{\unitlength}{0.8pt}
\begin{picture}(340,240)(-30,-10)\footnotesize
\put(0,120){
\begin{picture}(120,80)(-10,-10)
\put(0,0){\line(1,0){60}}
\put(60,0){\line(0,1){60}}
\put(60,60){\line(-1,0){30}}
\put(30,60){\line(-1,-1){30}}
\put(0,30){\line(0,-1){30}}
\put(0,0){\circle*{5}}
\put(30,0){\circle*{5}}
\put(60,0){\circle*{5}}
\put(60,30){\circle*{5}}
\put(60,60){\circle*{5}}
\put(30,60){\circle*{5}}
\put(0,30){\circle*{5}}
\put(30,30){\circle{3}}
\end{picture}
}
\put(140,120){
\begin{picture}(120,80)(-10,-10)
\put(0,0){\line(1,0){90}}
\put(90,0){\line(-1,1){60}}
\put(30,60){\line(-1,-1){30}}
\put(0,30){\line(0,-1){30}}
\put(0,0){\circle*{5}}
\put(30,0){\circle*{5}}
\put(60,0){\circle*{5}}
\put(90,0){\circle*{5}}
\put(60,30){\circle*{5}}
\put(30,60){\circle*{5}}
\put(0,30){\circle*{5}}
\put(30,30){\circle{3}}
\end{picture}
}
\put(280,120){
\begin{picture}(120,80)(-10,-10)
\put(0,0){\line(1,0){120}}
\put(120,0){\line(-1,1){60}}
\put(60,60){\line(-1,-1){60}}
\put(0,0){\circle*{5}}
\put(30,0){\circle*{5}}
\put(60,0){\circle*{5}}
\put(90,0){\circle*{5}}
\put(120,0){\circle*{5}}
\put(90,30){\circle*{5}}
\put(60,60){\circle*{5}}
\put(30,30){\circle*{5}}
\put(60,30){\circle{3}}
\end{picture}
}
\put(0,0){
\begin{picture}(120,80)(-10,-10)
\put(0,0){\line(1,0){60}}
\put(60,0){\line(0,1){60}}
\put(60,60){\line(-1,0){60}}
\put(0,60){\line(0,-1){60}}
\put(0,0){\circle*{5}}
\put(30,0){\circle*{5}}
\put(60,0){\circle*{5}}
\put(60,30){\circle*{5}}
\put(60,60){\circle*{5}}
\put(30,60){\circle*{5}}
\put(0,60){\circle*{5}}
\put(0,30){\circle*{5}}
\put(30,30){\circle{3}}
\end{picture}
}
\put(140,0){
\begin{picture}(120,80)(-10,-10)
\put(0,0){\line(1,0){90}}
\put(90,0){\line(-1,1){60}}
\put(30,60){\line(-1,0){30}}
\put(0,60){\line(0,-1){60}}
\put(0,0){\circle*{5}}
\put(30,0){\circle*{5}}
\put(60,0){\circle*{5}}
\put(90,0){\circle*{5}}
\put(60,30){\circle*{5}}
\put(30,60){\circle*{5}}
\put(0,60){\circle*{5}}
\put(0,30){\circle*{5}}
\put(30,30){\circle{3}}
\end{picture}
}
\put(280,0){
\begin{picture}(120,80)(-10,-10)
\put(0,0){\line(1,0){90}}
\put(90,0){\line(-1,1){90}}
\put(0,90){\line(0,-1){90}}
\put(0,0){\circle*{5}}
\put(30,0){\circle*{5}}
\put(60,0){\circle*{5}}
\put(90,0){\circle*{5}}
\put(60,30){\circle*{5}}
\put(30,60){\circle*{5}}
\put(0,90){\circle*{5}}
\put(0,60){\circle*{5}}
\put(0,30){\circle*{5}}
\put(30,30){\circle{3}}
\end{picture}
}
\put(30,105){$B_{11}$}
\put(180,105){$B_{12}$}
\put(330,105){$B_{13}$}
\put(30,-10){$B_{14}$}
\put(180,-10){$B_{15}$}
\put(330,-10){$B_{16}$}

\end{picture}
\setlength{\unitlength}{1pt}
\caption{\label{fig: reflexive rest}
\textit{Remaining reflexive polygons.}}
\end{figure}
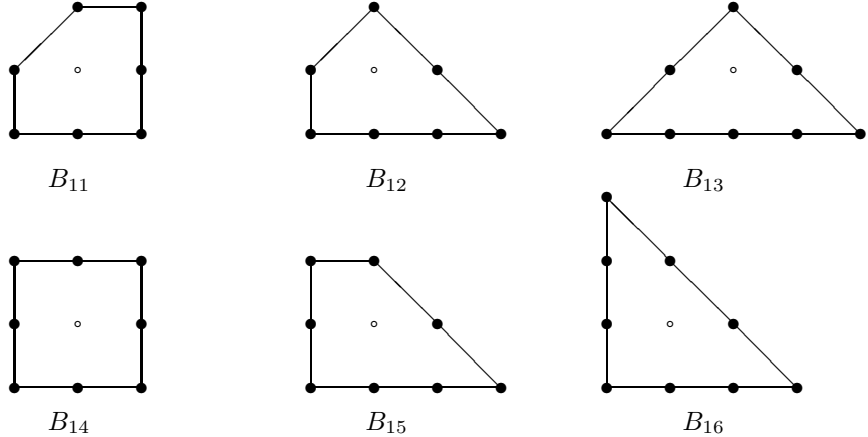

\

\ssnnl{ Example.}\label{examples quiver polygon}
Figures \ref{fig: models} and \ref{fig:dp3 models} show some examples
of the above relation between quivers, rank $2$ subgroups $\LL$ of $\ZZ^N$ as in \ref{introL} and polygons.
The quivers are taken from \cite{FFHH}. 
We name our examples neutrally as $B_1,\ldots,B_{10}$,
but also mention the names given to
these quiver models in \cite{FFHH}. The latter names refer to del Pezzo surfaces
$dP_k$, i.e. the projective plane $\PP^2$ blown up in $k$ points,
some of which, on the present occasion, may happen to be ``infinitely near''
and therefore require repeated blow-ups.
With the latter somewhat liberal use of the name \emph{del Pezzo surface},
$dP_k$ matches well with the toric geometry of the fan whose $1$-dimensional
rays are the half-lines from the interior lattice point $\circ$
through a lattice point $\bullet$ at the boundary.
Note however that this fan is in general not the same as the secondary fan
and that therefore the singularity associated with the polygon
(cf. \cite{HHV,MS}) is not the singularity obtained by contracting to a point
the zero-section of the canonical bundle of the del Pezzo surface.

The polygons in Figures \ref{fig: models} and \ref{fig:dp3 models}, however,
are in general different from, the
polygons associated with the same quiver in \cite{FFHH} Figure 8.
We will clarify this issue in  \ref{dP3IVtwist}, \ref{zigzagquiver}, \ref{twist}.

It is remarkable that the polygons in Figures \ref{fig: models} and 
\ref{fig:dp3 models} 
are exactly the lattice polygons
with one interior lattice point and $\leq 6$ lattice points on the boundary.
There are exactly $16$ lattice polygons with one interior lattice point.
The remaining $6$ are shown in Figure \ref{fig: reflexive rest}.
A nice review of these so-called \emph{reflexive polygons} can be found in \cite{PRV}.

In cases $B_1$, $B_3$ and $B_4$ the greatest common divisor $d$ of the
entries in the quiver's adjacency matrix is $>1$ and the polygon in 
these cases depends on an integer matrix with determinant $d$.
We have chosen this matrix such that the polygon is a reflexive polygon.

\section{GKZ hypergeometric systems}
\label{sec:hyper}
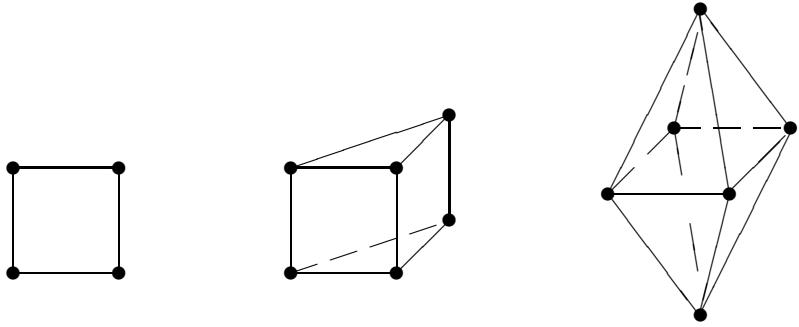
\begin{figure}[t]
\begin{center}
\begin{picture}(300,100)(30,10)
\put(35,20){
\begin{picture}(100,80)(0,0)
\put(0,0){\line(1,0){40}}
\put(0,0){\line(0,1){40}}
\put(40,40){\line(-1,0){40}}
\put(40,40){\line(0,-1){40}}
\put(0,0){\circle*{5}}
\put(40,0){\circle*{5}}
\put(0,40){\circle*{5}}
\put(40,40){\circle*{5}}
\end{picture}}
\put(140,20){
\begin{picture}(100,80)(0,0)
\put(0,0){\line(1,0){40}}
\put(0,0){\line(0,1){40}}
\put(40,40){\line(-1,0){40}}
\put(40,40){\line(0,-1){40}}
\put(60,20){\line(0,1){40}}
\put(40,0){\line(1,1){20}}
\put(40,40){\line(1,1){20}}
\put(0,40){\line(3,1){60}}
\multiput(0,0)(15,5){4}{\line(3,1){10}}
\put(0,0){\circle*{5}}
\put(40,0){\circle*{5}}
\put(0,40){\circle*{5}}
\put(40,40){\circle*{5}}
\put(60,20){\circle*{5}}
\put(60,60){\circle*{5}}
\end{picture}}
\put(260,20){
\begin{picture}(100,80)(0,0)
\put(0,30){\line(1,0){45}}
\put(45,30){\line(1,1){22}}
\multiput(0,30)(15,15){2}{\line(1,1){10}}
\multiput(25,55)(15,0){3}{\line(1,0){10}}
\put(0,30){\line(1,2){35}}
\put(0,30){\line(3,-4){35}}
\put(70,55){\line(-1,-2){35}}
\put(69,55){\line(-3,4){35}}
\put(46,30){\line(-1,6){12}}
\put(46,30){\line(-1,-4){11}}
\multiput(25,55)(6,-36){2}{\line(1,-6){3}}
\multiput(25,55)(5,20){2}{\line(1,4){4}}
\put(0,30){\circle*{5}}
\put(46,30){\circle*{5}}
\put(69,55){\circle*{5}}
\put(25,55){\circle*{5}}
\put(35,-16){\circle*{5}}
\put(35,100){\circle*{5}}
\end{picture}}
\end{picture}
\end{center}
\caption{\label{fig:classical examples}
\textit{The sets $\cA$ for Gauss's hypergeometric functions (left),
Appell's $F_1$ (middle) and Appell's $F_4$ (right).}}
\end{figure}

\ssnnl{ }\label{introGKZ}
In the late 1980's Gelfand, Kapranov and Zelevinsky discovered 
fascinating generalizations of the classical hypergeometric structures of Euler, Gauss, Appell, Lauricella, Horn \cite{gkz1,gkz2,gkz3,gkz4,S1}. 
The main ingredient for these new hypergeometric structures is a finite sequence $\cA\,=\,(\va_1,\ldots,\va_N)$ of vectors in $\ZZ^{k+1}$ which generates $\ZZ^{k+1}$ as an abelian group and for which there exists a group homomorphism 
$h:\ZZ^{k+1}\rightarrow\ZZ$ such that $h(\va_i)=1$ for all $i$.
The latter condition means that $\cA$ lies in a $k$-dimensional affine hyperplane in $\ZZ^{k+1}$.
Figure \ref{fig:classical examples} shows $\cA$ (each black dot represents
one vector)
sitting in this hyperplane for some classical hypergeometric structures.
Gelfand, Kapranov and Zelevinsky called these new structures  \emph{$\cA$-hypergeometric systems}, but nowadays most authors call them \emph{GKZ hypergeometric systems}.

Let $\LL$ denote the lattice (= free abelian group) of linear relations in $\cA$:
\begin{equation}\label{eq:L}
\LL:=\{(\ell_1,\ldots,\ell_N)\in\ZZ^N\;|\;
\ell_1\va_1+\ldots+\ell_N\va_N\,=\,\mathsf{0}\}\,.
\end{equation}
It follows from the above construction and assumptions that the quotient group
$\modquot{\ZZ^N}{\LL}$ is $\ZZ^{k+1}$, and a fortiori, it is torsion free.

In this paper we consider only the case when the rank of $\LL$ is $2$; i.e. $k+1=N-2$. Since we assumed the existence of a linear map $h:\ZZ^{N-2}\rightarrow\ZZ$ such that $h(\va_i)=1$ for all $i$,
$\LL$ lies in the kernel of the map $\vs$ defined by:
\begin{equation}\label{eq:s}
\vs:\RR^N\rightarrow\RR\,,\qquad \vs(z_1,\ldots,z_N)=z_1+\ldots+z_N\,.
\end{equation}

The following lemma shows what the other condition for $\LL$ in \ref{introL}
means in terms of $\cA$.

\

\ssnnl{ Proposition.}\label{lemma:circuits}
\textit{In the situation of \ref{introGKZ} the following statements are equivalent:
\begin{enumerate}
\item
$\LL$ is not contained in any of the standard coordinate hyperplanes of $\ZZ^N$.
\item
No vector in the sequence $\cA$ is linearly independent of the other vectors.
\item
For every $i$ the lattice of linear relations
between the vectors of $\cA\setminus\{\va_i\}$ has rank $1$.
\end{enumerate}
}

\proof Assume $\LL$ is contained in the $i$-th coordinate hyperplane. Then 
the $i$-th component of every element of $\LL$ is $0$ and hence $\va_i$ 
occurs in no linear relation for the vectors in $\cA$. In other words, $\LL$
is the lattice of linear relations
between the vectors of $\cA\setminus\{\va_i\}$. Since $\LL$ has rank $2$, we conclude
that statement 3 implies statement 1.

Conversely, assume that the lattice $\LL_i$ of linear relations
between the vectors of $\cA\setminus\{\va_i\}$ has rank $>1$.
Since $\LL_i$ is contained in $\LL$ and the latter has rank $2$,
we see that $\LL_i\otimes\RR=\LL\otimes\RR$. Therefore, since $\LL_i$ is contained in the $i$-th coordinate hyperplane, so is $\LL$. Hence statement 1 implies statement 3.

It is obvious that statements 1 and 2 are equivalent.
\qed

\

\ssnnl{ Definition.}\label{def:minrel}
A subsequence $\cA'$ of $\cA$ is said to be \emph{minimally dependent} 
if the vectors in $\cA'$ are linearly dependent and the vectors
in every subsequence $\cA''$ of $\cA'$ with $\cA''\neq\cA'$ are linearly independent.
A linear dependence relation $\sum_j \alpha_j\va_j=0$ is a 
\emph{minimal linear dependence relation} in $\cA$ if the subsequence 
$(\va_j\:|\:\alpha_j\neq 0)$ is minimally dependent.

\

\ssnnl{ Proposition.}\label{prop:circuits}
\textit{
\begin{enumerate}
\item
Let $B$ be a $2\times N$-matrix such that its rows are a $\ZZ$-basis for
$\LL$. Let $\vb_1,\ldots,\vb_N$ be the columns of this matrix.
Then one has for every $i$
\begin{equation}\label{eq:minrel}
\sum_{j=1}^N \det(\vb_i,\vb_j)\,\va_j\:=\:\nv\,.
\end{equation}
\item
Assume the conditions in Proposition \ref{lemma:circuits} are satisfied.
Then every linear relation \textup{(\ref{eq:minrel})} is minimal
and every minimal linear dependence relation in $\cA$ is a non-zero scalar 
multiple of some relation \textup{(\ref{eq:minrel})}.
\end{enumerate}
}
\proof
For statement 1,
let $A$ be the matrix with columns $\va_1,\ldots,\va_N$. Then Equation
(\ref{eq:L}) can be rewritten as
$BA^t=0$. This implies $B^t\sJ BA^t=0$, which is just a compressed form of the relations we wanted to prove.

For statement 2 note that condition 3 in Proposition \ref{lemma:circuits} 
implies that every linear relation (\ref{eq:minrel}) is minimal.
Conversely,
consider a minimal linear relation 
$\sum_j \alpha_j\va_j=0$. Then the subsequence 
$(\va_j\:|\:\alpha_j\neq 0)$ has at most $N-1$ elements.
The given minimal relation is therefore a non-zero scalar multiple of
at least one relation (\ref{eq:minrel}).
\qed

\

\ssnnl{}\label{GKZsystem}
Gelfand, Kapranov and Zelevinsky \cite{gkz1,gkz2,gkz3}
associate with a set $\cA$ as in \ref{introGKZ} and a vector $\vc\in\CC^{k+1}$
the following system of partial differential equations
for functions $\Phi$ of $N$ variables $u_1,\ldots,u_N$: 
\begin{itemize}
\item
for every $(\ell_1,\ldots,\ell_N)\in\LL$
one differential equation
\begin{equation}\label{eq:GKZ1}
\prod_{\ell_i<0}\left(\frac{\partial}{\partial u_i}\right)^{-\ell_i}\Phi
\;=\;
\prod_{\ell_i>0}\left(\frac{\partial}{\partial u_i}\right)^{\ell_i}\Phi\,,
\end{equation}
\item
the system of $k+1$ differential equations
\begin{equation}\label{eq:GKZ2}
\va_1\,u_1\frac{\partial  \Phi}{\partial u_1}+
\ldots+\va_N\,u_N\frac{\partial  \Phi}{\partial u_N}\;=\;
\vc\Phi\,.
\end{equation} 
\end{itemize}

\

\ssnnl{ Example.}\label{cubic roots}
The sequence $\cA=(\va_1,\va_2,\va_3,\va_4)\in\ZZ^2$ for Example $B_2$ in 
Figure \ref{fig: models} can be taken to be
$$
\va_1=\left[\begin{array}{r}1\\ 0\end{array}\right]\,,\quad
\va_3=\left[\begin{array}{r}1\\ 1\end{array}\right]\,,\quad
\va_4=\left[\begin{array}{r}1\\ 2\end{array}\right]\,,\quad
\va_2=\left[\begin{array}{r}1\\ 3\end{array}\right]\,.
$$
It is a classical result of K. Mayr that the roots of the $1$-variable
cubic polynomial
$P(x)=u_1+u_3x+u_4x^2+u_2x^3$ as functions of the coefficients $u_1,\ldots,u_4$
satisfy the GKZ system of differential equations for this
$\cA$ and with $\vc=\left[\begin{array}{r}0\\ -1\end{array}\right]$;
see e.g. \cite{S1} \S 2.1.

\

\ssnnl{ Example.}\label{Labc}
An important example in the toric geometry constructions of Sasaki-Einstein
manifolds is known under the name $L^{a,b,c}$; see e.g. \cite{FHMSVW}.
Here $a,b,c$ are integers with $c\leq b$ and $0< a\leq b$
The polygon for this example, displayed in \cite{FHMSVW} Figure 2, is a 
quadrangle with vertices $(0,0)$, $(1,0)$, $(ak,b)$, $(-al,c)$ where $k$ and 
$l$ are integers such that $ck+bl=1$. The method explained in 
\ref{polygon to L} yields $\LL=\ZZ(1,ak-1,-al-ak,al)\oplus\ZZ(0,b,c-b,-c)$.
From this we see that we can take $\cA=(\va_1,\va_2,\va_3,\va_4)\in\ZZ^2$ with
$$
\va_1=\left[\begin{array}{c}1\\ c-a\end{array}\right]\,,\quad
\va_2=\left[\begin{array}{c}1\\ c\end{array}\right]\,,\quad
\va_3=\left[\begin{array}{c}1\\ 0\end{array}\right]\,,\quad
\va_4=\left[\begin{array}{c}1\\ b\end{array}\right]\,.
$$
As in Example \ref{cubic roots} this means that the corresponding GKZ system deals with the roots of the
``four-nomial''  $u_1x^{c-a}+u_2x^c+u_3+u_4x^b$ as functions of
the coefficients $u_1,u_2,u_3,u_4$.

According to \cite{FHMSVW} \S 3 the examples $Y^{p,q}$ of \cite{BFHMS,MS} are special cases of the above:
$Y^{p,q}=L^{p-q,p+q,p}$. So, these correspond to the ``four-nomials''\\
$u_1x^q+u_2x^p+u_3+u_4x^{p+q}$. Thus Example \ref{cubic roots} is in fact $Y^{2,1}=L^{1,3,2}$.

\

\ssnnl{ Example.}\label{F14}
From the pictures of $\cA$ shown in Figure \ref{fig:classical examples}
one easily sees that for the Gauss system the lattice $\LL$ is generated
by the vector $(1,1,-1,-1)$ and thus has rank $1$ and is not of interest here.
For Appell's $F_1$ and $F_4$ the lattice $\LL$ has rank $2$ and the corresponding $B$-matrices are
$$
\;F_1:\; 
\left[\begin{array}{rrrrrr}1&-1&0&-1&1&0\\1&0&-1&-1&0&1\end{array}\right]
\,,\quad
F_4:\;
\left[\begin{array}{rrrrrr}-1&-1&1&1&0&0\\-1&-1&0&0&1&1\end{array}\right].
$$

\

Thus Appell's $F_1$ corresponds to $B_7$ in Figure \ref{fig:dp3 models}.
For Appell's $F_4$ the quiver and the polygon are

\setlength{\unitlength}{0.8pt}
\begin{picture}(300,90)(-10,-10)
\put(80,20){\begin{picture}(100,80)(0,0)
\put(0,0){\circle*{5}}
\put(0,30){\circle*{5}}
\put(30,-15){\circle*{5}}
\put(60,0){\circle*{5}}
\put(60,30){\circle*{5}}
\put(30,45){\circle*{5}}
\put(-10,27){$1$}
\put(-11,-3){$6$}
\put(28,-28){$4$}
\put(65,-3){$2$}
\put(65,27){$5$}
\put(28,49){$3$}
\put(0,30){\vector(0,-1){27}}
\put(0,0){\vector(2,-1){27}}
\put(30,-15){\vector(2,1){27}}
\put(60,0){\vector(0,1){27}}
\put(60,30){\vector(-2,1){27}}
\put(30,45){\vector(-2,-1){27}}
\put(0,30){\vector(1,0){57}}
\put(0,0){\vector(2,3){28}}
\put(30,45){\vector(2,-3){28}}
\put(60,0){\vector(-1,0){57}}
\put(60,30){\vector(-2,-3){28}}
\put(30,-15){\vector(-2,3){28}}
\end{picture}}

\put(250,30){\begin{picture}(100,80)(0,0)
\put(0,-30){\circle*{5}}
\put(30,-30){\circle*{5}}
\put(60,-30){\circle*{5}}
\put(60,0){\circle*{5}}
\put(60,30){\circle*{5}}
\put(30,0){\circle*{5}}

\put(0,-30){\line(1,0){60}}
\put(60,-30){\line(0,1){60}}
\put(60,30){\line(-1,-1){60}}

\end{picture}}
\end{picture}
\setlength{\unitlength}{1pt}

In \cite{gkz1} \S 3.2  Gelfand, Kapranov and Zelevinsky discuss how the
$14$ complete hypergeometric Horn series in two variables fit their theory.
It is an amusing simple exercise to now find the corresponding quivers and polygons. It turns out that, with the exception of $F_4$, all Horn series
give a quiver listed in Example \ref{examples quiver polygon}:
\begin{eqnarray*}
&&\{G_3\}\leftrightarrow B_2\,,\quad \{H_5\}\leftrightarrow B_5\,,\quad
\{G_1,\,H_3,\,H_6\}\leftrightarrow B_6\,,\quad \{F_1,\,G_2\}\leftrightarrow B_7\,,\\ &&\{H_1\}\leftrightarrow B_8\,,\quad
\{H_4,\,H_7\}\leftrightarrow B_9\,,\quad 
\{F_2,\,F_3,\,H_2\}\leftrightarrow B_{11}\,.
\end{eqnarray*}

\

\ssnnl{}\label{L to A}
In \ref{introGKZ} we started from the sequence of vectors $\cA$ and then
defined $\LL$ via Equation (\ref{eq:L}). 
Let us reverse the procedure and start with a rank $2$ subgroup
$\LL\subset\ZZ^N$ as in \ref{introL}. We can now define 
\begin{eqnarray}\label{eq:AfromL}
\cA&=&(\va_1,\ldots,\va_N)\qquad\textrm{with}\\
\nonumber
\va_i&=& \textrm{the class}\;\ve_i\,\bmod\,\LL\; \textrm{in the quotient group}
\;\modquot{\ZZ^N}{\LL}\,.
\end{eqnarray}
If the quotient group $\modquot{\ZZ^N}{\LL}$ is torsion free, it is isomorphic to $\ZZ^{N-2}$
and the sequence of vectors $\cA$ in (\ref{eq:AfromL}) satisfies the 
requirements of \ref{introGKZ}.
\emph{Thus (\ref{eq:L}) and (\ref{eq:AfromL}) give equivalences, inverse 
to each other, between the data $\cA$ and $\LL$.}

In combination with Theorem \ref{quivergrassmann} this gives

\

\ssnnl{ Theorem.}\label{quiverhypergeo}
\emph{There is an equivalence of data between on the one hand
sequences $\cA$ as in \ref{introGKZ} which satisfy also the conditions in 
Proposition \ref{lemma:circuits} and on the other hand
quivers without isolated nodes or directed loops of length $\leq 2$,
such that in every node the number of incoming arrows equals the number of outgoing arrows and such that for its anti-symmetrized adjacency matrix 
the rank is $2$ and the greatest common divisor of the entries
is $1$.}
\qed

\

\ssnnl{ Remark.}\label{GKZwithtorsion}
We refer to \ref{GKZ with torsion} for a discussion of how to
extend the definition of GKZ systems so as to allow
for torsion in the group $\modquot{\ZZ^N}{\LL}$.

\

\ssnnl{}\label{A secondary}
In order to have an efficient method to construct bases for the solution space of the hypergeometric system (\ref{eq:GKZ1})-(\ref{eq:GKZ2})
Gelfand, Kapranov and Zelevinsky developed the theory of the \emph{secondary fan} 
and the \emph{secondary polytope}. 
The term ``secondary'' refers to the habit of considering
the set $\cA$ as primary data and to call the convex hull
of $\cA$ the \emph{primary polytope}. 

In case the lattice
$\LL$ of relations in $\cA$ has rank $2$, the constructions become 
particularly simple. The secondary fan is the one in Definition
\ref{def:Lsecondary fan}. 
With a $2$-dimensional cone $C$ in the secondary fan
one associates a set $L_C$ of $2$-element subsets of $\{1,\ldots,N\}$
as in Equation (\ref{eq:LC}).
One can interpret this $L_C$ is as a triangulation
of the primary polytope $\mathsf{convex\;hull}(\cA)$ as follows.
Recall that $\cA\,=\,(\va_1,\ldots,\va_N)$ and assume for convenience of imagining pictures that all $\va_i$'s are different.
For $\{i,j\}\in L_C$ set $T_{\{i,j\}}\,:=\,
\mathsf{convex\;hull}(\{\va_k\:|\: k\neq i,j\:\})$.
This $T_{\{i,j\}}$ is an $(N-3)$-dimensional simplex and the simplices
$T_{\{i,j\}}$ with $\{i,j\}\in L_C$ together constitute a \emph{triangulation
the primary polytope}.

In \cite{gkz4} p.220 the \emph{secondary polytope} is constructed (and defined)
in terms of these triangulations as follows. For a $2$-element subset  $\{i,j\}\subset\{1,\ldots,N\}$
with $i<j$
let $A_{ij}$ denote the $(N-2)\times(N-2)$-matrix with columns
$\va_k$ for $k\neq i, j$ in the natural order of increasing indices.
To a $2$-dimensional cone $C$ in the secondary fan one then associates
the vector
\begin{equation}\label{eq:phiC}
\gf_C\,:=\,\sum_{\{i,j\}\in L_C}|\det\,A_{ij}|\sum_{k\neq i,j}\ve_k\,.
\end{equation}
Then in \cite{gkz4} p.220 the \emph{secondary polytope} $\Sigma (\cA)$ is defined as
\begin{equation}\label{eq:gkz secondary polytope}
\Sigma (\cA)\,:=\,\mathsf{convex\;hull}(\{\gf_C\:|\: C\:\textrm{$2$-dim cone
of secondary fan}\,\})\,.
\end{equation}
As $|\det\,A_{ij}|$ is $(N-3)!$ times the Euclidean volume of the simplex
$T_{\{i,j\}}$, the number 
$\mathrm{vol}_\cA:=\sum_{\{i,j\}\in L_C}|\det\,A_{ij}|$ is for every $C$ equal to $(N-3)!$ times
the Euclidean volume of the primary polytope
$\mathsf{convex\;hull}(\cA)$
and the point 
$$
\frac{1}{(N-2)\mathrm{vol}_\cA}\sum_{\{i,j\}\in L_C}|\det\,A_{ij}|
\sum_{k\neq i,j}\va_k
$$
coincides for every $C$ with the barycentre of the primary polytope.
The secondary polytope therefore lies in a $2$-dimensional plane
parallel to $\LL_\RR$.

\

\ssnnl{}\label{compare secondary polytopes}
Let us compare $\Sigma (\cA)$ with the polytope $\Sigma (\LL)$
of Definition \ref{def:Lsecondary polytope}. 
It is a well-known and easy to prove fact that
$|\det\,A_{ij}|\,=\,|\det (\vb_i,\vb_j)|$
for all $\{i,j\}\in L_C$ (see e.g. \cite{S1} Eq. (62)).
Thus Equations (\ref{eq:phiC}) and (\ref{eq:psiC}) yield
\begin{equation}\label{eq:phiCpsiC}
\gf_C\,=\,-\psi_C\,+\,\mathrm{vol}_\cA\cdot\sum_{k=1}^N\ve_k\,,
\end{equation}
with
\begin{equation}\label{eq:volA}
\mathrm{vol}_\cA\,=\,\sum_{\{i,j\}\in L_C}|\det\,A_{ij}|
\,=\,\sum_{\{i,j\}\in L_C}|\det (\vb_i,\vb_j)|\,.
\end{equation}
This means that the two secondary polytopes $\Sigma(\cA)$ and $\Sigma(\LL)$ are 
related by a point symmetry with centre 
$\frac{1}{2}\mathrm{vol}_\cA\cdot\sum_{k=1}^N\ve_k$:
\begin{equation}\label{eq:two secondary polytopes}
\Sigma(\cA)\:=\:-\,\Sigma(\LL)\,+\,\mathrm{vol}_\cA\cdot\sum_{k=1}^N\ve_k\,.
\end{equation}

\section{About the solutions to GKZ systems.}
\label{sec:solutions}

\

\ssnnl{}\label{GKZtorusaction}
The variables in GKZ theory (see Section \ref{sec:hyper}) are the natural 
coordinates on the space $\CC^\cA:=\mathrm{Maps}(\cA,\CC)$ of maps 
from $\cA$ to $\CC$. The torus $\TT^{k+1}:=\mathrm{Hom}(\ZZ^{k+1},\CC^*)$ of 
group homomorphisms from $\ZZ^{k+1}$ to $\CC^*$, acts naturally on 
$\CC^\cA$ and on the functions on this space: 
for $\sigma\in\TT^{k+1}$, 
$\vu\in\CC^\cA$, $\va\in\cA$ and $\Phi:\CC^\cA\rightarrow\CC$:
\begin{equation}\label{eq:torus action on functions}
(\sigma\cdot\vu)(\va)\,=\, \sigma(\va)\vu(\va)\,,\qquad
(\Phi\cdot\sigma)(\vu)\,=\,\Phi(\sigma\cdot\vu)\,.
\end{equation}

The GKZ hypergeometric functions associated with $\cA$ and $\vc$
are defined on open domains in $\CC^\cA$. One easily sees that
if a function $\Phi$ on $\CC^\cA$ satisfies the differential equations
(\ref{eq:GKZ1}) then for every $\sigma\in\TT^{k+1}$ the function 
$\Phi\cdot\sigma$ also satisfies these differential equations.
So, the torus $\TT^{k+1}$ acts on the solution space of the system of differential 
equations (\ref{eq:GKZ1}).

On the other hand, if $\Phi_1$ and $\Phi_2$ are two functions which satisfy the
differential equations (\ref{eq:GKZ2}) with the same $\vc$, then their
quotient $\Psi=\frac{\Phi_1}{\Phi_2}$ satisfies
$$
\va_1\,u_1\frac{\partial  \Psi}{\partial u_1}+
\ldots+\va_N\,u_N\frac{\partial  \Psi}{\partial u_N}\;=\;\nv\,.
$$
The latter equation is equivalent to $\Psi$ being $\TT^{k+1}$-invariant.
Thus we find:

\emph{All quotients of pairs of solutions of the system 
(\ref{eq:GKZ1})-(\ref{eq:GKZ2}) are functions on simply connected open subsets of
the orbit space $\modquot{\CC^\cA}{\TT^{k+1}}$.
In the case of interest in the present paper the dimension of this orbit space
is $N-k-1=2$.}

For $\vc\in\ZZ^{k+1}$ the differential equations (\ref{eq:GKZ2})
are equivalent with $\Phi$ transforming under the action of $\TT^{k+1}$
according to the character given by $\vc$:
\begin{equation}\label{eq:character}
\Phi\cdot\sigma=\sigma(\vc)\Phi\,.
\end{equation}
In particular for $\vc=\nv$ all solutions of (\ref{eq:GKZ2}) are $\TT^{k+1}$-invariant.

\

\ssnnl{ }\label{toric sec fan}
The space $(\CC^*)^\cA:=\mathrm{Maps}(\cA,\CC^*)$ of maps 
from $\cA$ to $\CC^*$ is a torus of dimension $N$
which contains $\TT^{k+1}$ as a subtorus. The action 
(\ref{eq:torus action on functions}) of $\sigma\in\TT^{k+1}$
on $\CC^\cA$ obviously restricts to the action of $\TT^{k+1}$ on $(\CC^*)^\cA$ by multiplication.
The quotient space $\modquot{(\CC^*)^\cA}{\TT^{k+1}}$, which is a subspace of
$\modquot{\CC^\cA}{\TT^{k+1}}$, is the torus 
$\mathrm{Hom}(\LL,\CC^*)$ of 
group homomorphisms from $\LL$ to $\CC^*$.
The toric variety associated with the secondary fan (cf. Definition 
\ref{def:Lsecondary fan}) gives a compactification of the torus 
$\mathrm{Hom}(\LL,\CC^*)$; see \cite{F} for the general theory of toric varieties. It is an essential part of the GKZ philosophy that
quotients of hypergeometric functions should be viewed as being defined on
open subsets of this toric variety.

\

\ssnnl{ Theorem.}\label{thm:solution dimension}
(cf. \cite{gkz1} Theorems 2 and 5, $\cite{gkz1}'$, \cite{stu} Prop. 13.5 )
\textit{Let}
\begin{eqnarray*}
\NN\cA&=&\{x_1\va_1+\ldots+x_N\va_N\in\RR^{k+1}\:|\:\forall x_i\in\ZZ_{\geq 0}\}
\,,
\\
\ZZ\cA&=&\{x_1\va_1+\ldots+x_N\va_N\in\RR^{k+1}\:|\:\forall x_i\in\ZZ\}\,,
\\
\mathrm{pos}(\cA)&=&\{x_1\va_1+\ldots+x_N\va_N\in\RR^{k+1}\:|\:\forall x_i\in\RR_{\geq 0}\}\,.
\end{eqnarray*}
\textit{Assume $\NN\cA\,=\,\ZZ\cA\cap\mathrm{pos}(\cA)$,
then the dimension of the space of solutions of the system of differential equations (\ref{eq:GKZ1})-(\ref{eq:GKZ2}) at a general point of $\CC^\cA$
equals the number $\mathrm{vol}_\cA$ in Equation (\ref{eq:volA}).}
\qed

\

\ssnnl{ Remark.}\label{unimodular}
As shown in the proof of \cite{stu} Prop. 13.15 the condition
$\NN\cA\,=\,\ZZ\cA\cap\mathrm{pos}(\cA)$ in the above theorem is satisfied
if $\cA$ admits a \emph{unimodular triangulation}. The latter condition is equivalent to: there is a cone $C$ in the secondary fan such that 
$|\det(\vb_i,\vb_j)|=1$ for all $\{i,j\}\in L_C$ 
(see Definitions \ref{def:Lsecondary fan} and \ref{def:Lsecondary polytope}).

\

\ssnnl{ Remark.}\label{GKZ with torsion} 
From the discussion in \ref{GKZtorusaction} one may get an idea about
extending the GKZ system (\ref{eq:GKZ1})-(\ref{eq:GKZ2}) to the situation
in which $\modquot{\ZZ^N}{\LL}$ has a non-trivial torsion
subgroup $\left(\modquot{\ZZ^N}{\LL}\right)_{\mathrm{tors}}$.
Equation (\ref{eq:GKZ1}) still makes sense in this more general situation,
but Equation (\ref{eq:GKZ2}) must be adapted.

Recall from (\ref{eq:AfromL}) a definition of $\cA$ which also works in the torsion case. Let $\cG_\cA:=
\mathrm{Hom}(\modquot{\ZZ^N}{\LL},\CC^*)$. This is a commutative algebraic
group of which the connected component of the identity is
$\cG_\cA^\circ=\TT^{k+1}=\mathrm{Hom}(\ZZ^{k+1},\CC^*)$ and the group of connected components
is the finite abelian group 
$\modquot{\cG_\cA\:}{\cG_\cA^\circ}=
\mathrm{Hom}(\left(\modquot{\ZZ^N}{\LL}\right)_{\mathrm{tors}},\CC^*)
\simeq \left(\modquot{\ZZ^N}{\LL}\right)_{\mathrm{tors}}$.
Formula (\ref{eq:torus action on functions}) defines an action of
$\cG_\cA$ on $\CC^\cA$ and on the functions on $\CC^\cA$.
It is now clear how to adapt Equation (\ref{eq:GKZ2}):
let $\overline{\va}_i$ denote the projection of $\va_i$ in the free
part $\ZZ^{k+1}$ of $\modquot{\ZZ^N}{\LL}$ and replace (\ref{eq:GKZ2})
by
$$
\overline{\va}_1\,u_1\frac{\partial  \Phi}{\partial u_1}+
\ldots+\overline{\va}_N\,u_N\frac{\partial  \Phi}{\partial u_N}\;=\;
\vc\Phi
$$
plus the requirement that the solution should
transform according to some character of the finite abelian group
$\left(\modquot{\ZZ^N}{\LL}\right)_{\mathrm{tors}}$.

\

\ssnnl{ Example.}\label{example B1}
In case $B_1$ in Figure \ref{fig: models} one has
$\LL=\ZZ(2,-1,-1)\oplus\ZZ(1,1,-2)$ and
$\modquot{\ZZ^3}{\LL}=\ZZ\oplus \modquot{\ZZ\,}{3\ZZ}$. As a generator
for the torsion subgroup we take
$\vg=(1,-1,0)\bmod\LL$. 
Then the polynomials $\Phi_0=\frac{1}{2}u_3^2+u_1u_2$, 
$\Phi_1=\frac{1}{2}u_2^2+u_1u_3$
and $\Phi_2=\frac{1}{2}u_1^2+u_2u_3$ 
satisfy the differential equations (\ref{eq:GKZ1})
for $\LL$, the differential equations (\ref{eq:GKZ2}) for $\vc=2$,
while
$\vg\cdot\Phi_r=e^{2\pi ir/3}\Phi_r$ for $r=0,1,2$.

\section{From rank $2$ subgroups of $\ZZ^N$ to Dessins.}
\label{sec:dessin}
We are going to describe a construction which associates with a rank $2$
subgroup $\LL$ of $\ZZ^N$, as in \ref{introL},
\emph{dessins d'enfants}, i.e. bipartite graphs embedded in oriented
Riemann surfaces. 

The construction is given in \cite{HV} \S 5 and \S6
as \emph{Fast Inverse Algorithm}. In op. cit., however, this algorithm
is only presented via explicit visual inspection of pictures in 
some concrete examples.
In this section we want to present a general
principle behind the Fast Inverse Algorithm of \cite{HV} which 
uses only linear algebra and can be performed by computer.

The main part of our construction produces rhombus tilings of
the plane and
is the same as N.G. de Bruijn's \cite{B}
construction of Penrose tilings. His work also led to effective methods
for making \emph{quasi-crystals} \cite{Sen}. Since $\LL$ is 
defined over $\ZZ$ the construction yields in our situation periodic 
tilings of the plane $\LL_\RR^\vee=\RR^2$, not just quasi-periodic ones.
We can therefore pass to $\RR^2$ modulo the period lattice and find a tiling of the two-dimensional torus.

\

\ssnnl{ }\label{grids}
In this section $\LL$ is as in \ref{introL}. The quotient group $\modquot{\ZZ^N}{\LL}$ is allowed to have torsion.
Let $\LL_\RR$ denote the real $2$-plane in $\RR^N$ which
contains $\LL$. We write $\gl+\LL_\RR$ for the real $2$-plane in $\RR^N$ obtained by translating $\LL_\RR$ over a vector $\gl\in\RR^N$.
On the other hand one has in $\RR^N$ the \emph{standard $N$-grid}
consisting of the hyperplanes
$$
\cH_{i,k}:=\{(z_1,\ldots,z_N)\in\RR^N\:|\:z_i=k\}\qquad
\textrm{for}\; i=1,\ldots,N\;\textrm{and}\;k\in\ZZ\,.
$$
By intersecting with this standard $N$-grid
we obtain in the $2$-plane $\gl+\LL_\RR$ an $N$-grid of lines:
\begin{equation}\label{eq:grid L}
\cL^\gl_{i,k}:=\cH_{i,k}\cap(\gl+\LL_\RR)\qquad
\textrm{for}\; i=1,\ldots,N\;\textrm{and}\;k\in\ZZ\,.
\end{equation}
Note that this crucially uses the assumption that $\LL$ is not contained in any of the standard coordinate hyperplanes in $\ZZ^N$.

\emph{The grid lines are naturally oriented}: the orientation on $\cL^\gl_{i,k}$
is such that the points in $\gl+\LL_\RR$ which are to the left of 
$\cL^\gl_{i,k}$, have a larger $i$-th coordinate than the points to the right
of $\cL^\gl_{i,k}$.

\

\ssnnl{ Definition.}\label{def:nonresonant}
We say that $\gl$ is \emph{non-resonant} if
no three grid lines pass through one point.

\

\ssnnl{ First step.}\label{first step} 
Let $\LL$ and $\gl$ be as in \ref{grids} with $\gl$ non-resonant
as in \ref{def:nonresonant}.
Consider the map 
$$
F:\RR^N\rightarrow\ZZ^N\,,\qquad 
F(z_1,\ldots,z_N)=(\lfloor z_1\rfloor,\ldots\lfloor z_N\rfloor)\,,
$$
where $\lfloor z\rfloor$ for a real number $z$ denotes the largest integer
$\leq z$.
The map $F$ contracts an open $N$-cube in the $\ZZ^N$ structure on $\RR^N$ 
(and part of its boundary) onto one of its corners. 
$F$ commutes with the translation action of $\ZZ^N$ on $\RR^N$ and $\ZZ^N$.
Suppressing $\LL$ and $\gl$ from the notation 
we define
\begin{equation}\label{eq:S0}
\cS^0:=F(\gl+\LL_\RR)\,.
\end{equation}
The map $F$ is constant on each $2$-cell ( = connected component) 
of the grid complement in $\gl+\LL_\RR$. Every point of $\cS^0$ is in fact the image of a unique such $2$-cell.
The distance between two points of 
$\cS^0$ is $1$ if and only if the corresponding $2$-cells are separated by exactly one grid line.

Because $\gl$ is non-resonant, an intersection point of 
grid lines is in the closure of exactly four $2$-cells. 
If $\vx= \cL^\gl_{i,k}\cap \cL^\gl_{j,m}$ the four points of $\cS^0$
corresponding to these cells, i.e. 
$F(\vx)$, $F(\vx)-\ve_i$, $F(\vx)-\ve_j$, $F(\vx)-\ve_i-\ve_j$,
are the vertices of a square $\square_{\vx}$.
When $\vx$ runs through the set of all intersection points of grid lines,
these squares fit together to a connected surface 
$\cS:=\bigcup_\vx \square_{\vx}$ embedded in $\RR^N$.
One can characterize the surface $\cS$ also as
\begin{equation}\label{eq:squares}
\cS:=\textit{union of all unit squares in $\RR^N$ with vertices in the set $\cS^0$\,.}
\end{equation}
Here one may define 
\emph{unit square} as the convex hull of four points 
$\vp_1,\vp_2,\vp_3,\vp_4$ in $\RR^N$
with Euclidean distances
$\parallel\!\!\vp_1-\vp_2\!\!\parallel=\parallel\!\!\vp_2-\vp_3\!\!\parallel=
\parallel\!\!\vp_3-\vp_4\!\!\parallel=$ $\parallel\!\!\vp_4-\vp_1\!\!\parallel=1$
and $\parallel\!\!\vp_1-\vp_3\!\!\parallel=\parallel\!\!\vp_2-\vp_4\!\!\parallel=
\sqrt{2}$.
A unit square with vertices in $\ZZ^N$ is necessarily  of the form
$\vp+\mathrm{convex\, hull}(\nv,\ve_i,\ve_j,\ve_i+\ve_j)$ 
for some $\vp\in\ZZ^N$ and some $i,j$.
We say that a unit square in $\RR^N$ has 
\emph{type $(i,j)$} if its sides are parallel to the vectors $\ve_i$ or $\ve_j$.
Similarly, a side of a unit square in $\RR^N$ has 
\emph{type $i$} if it is parallel to the vector $\ve_i$.

\

\ssnnl{ Proposition.}\label{prop:Sprojection}
\emph{Let $\LL^\vee_\RR=\mathrm{Hom}(\LL,\RR)$ denote the real dual
of $\LL$. Identify the real dual $\mathrm{Hom}(\ZZ^N,\RR)$ of $\ZZ^N$ with
$\RR^N$ by means of the standard dot product.
Then the linear map $\RR^N\rightarrow\LL^\vee_\RR$ which is dual to the inclusion $\LL\rightarrow\ZZ^N$, restricts to a homeomorphism
$\cS\stackrel{\sim}{\longrightarrow}\LL^\vee_\RR$.
}

\proof Let $\vb_1,\ldots,\vb_N\in\LL^\vee_\RR$ denote the images of the standard basis vectors of $\RR^N$. After identifying in the obvious way,
the real plane $\LL_\RR$ with the real $2$-plane $\gl+\LL_\RR$ in $\RR^N$
one obtains in $\LL_\RR$ an $N$-grid in which the grid lines $\cL^\gl_{i,k}$
are perpendicular to $\vb_i$.
After choosing a basis for $\LL$ and taking the dual basis for $\LL_\RR^\vee$
we may identify $\LL_\RR$ and $\LL_\RR^\vee$ with $\RR^2$. The term
``perpendicular'' then means perpendicular with respect to the standard inner product on $\RR^2$.

\begin{figure}
\begin{picture}(350,130)(0,0)
\put(-20,0){
\begin{picture}(100,100)(0,15)
\setlength\epsfxsize{8cm}
\epsfbox{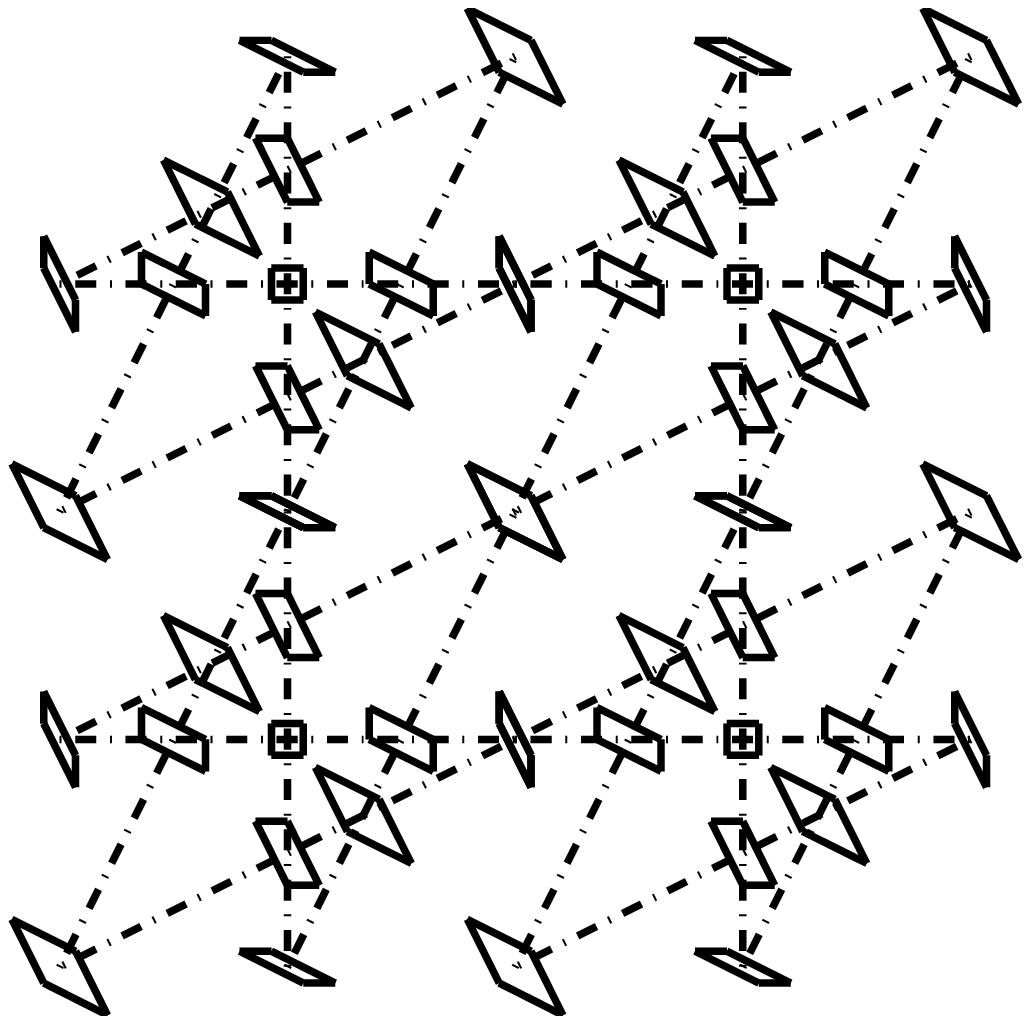}
\end{picture}
}
\put(180,0){
\begin{picture}(150,100)(0,0)
\setlength\epsfxsize{7cm}
\epsfbox{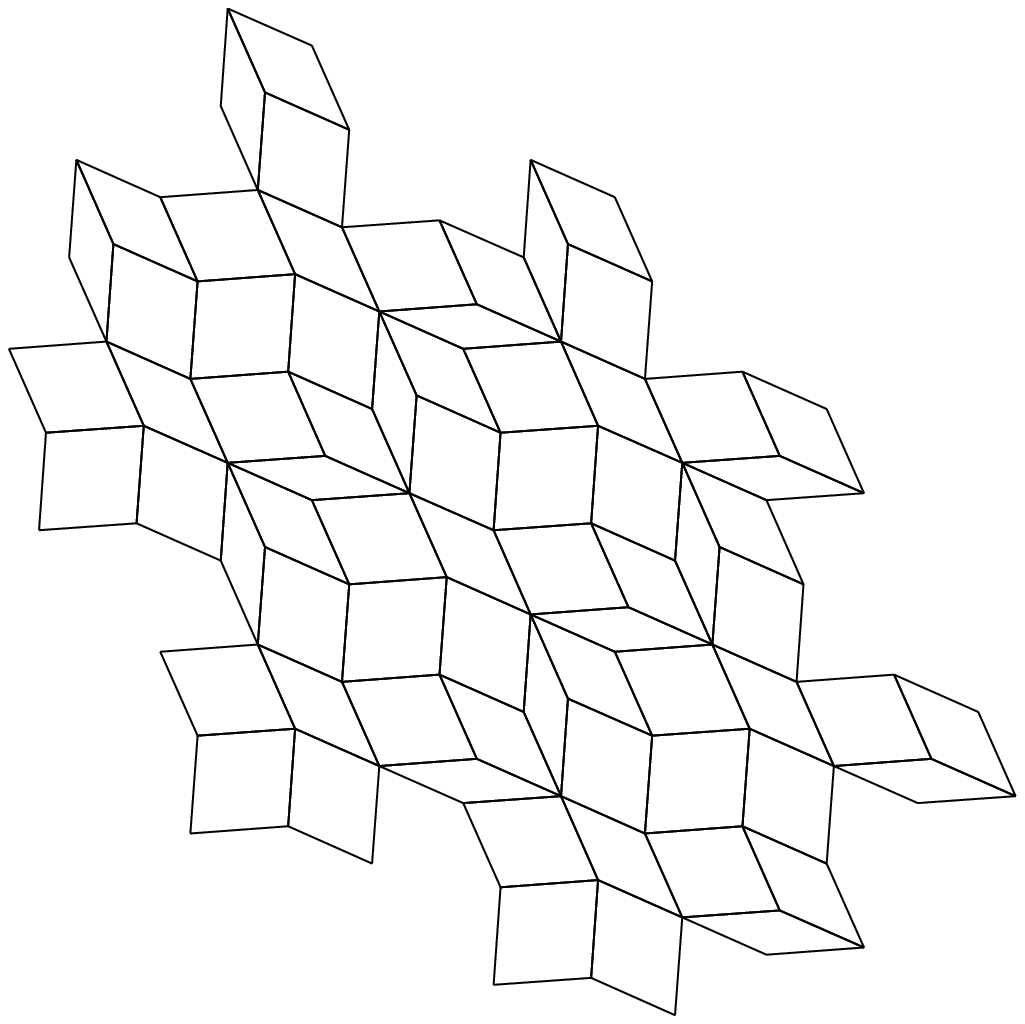}
\end{picture}
}
\end{picture}
\caption{\label{fig:grid2tiling}
\textit{converting grid to tiling for $B_2$ from Figure \ref{fig: models}}}
\end{figure}

Next draw for each intersection point of grid lines, say
$\vx= \cL^\gl_{i,k}\cap \cL^\gl_{j,m}$, a parallelogram with centre $\vx$
and sides $\epsilon\vb_i$ and $\epsilon\vb_j$. Here the positive
real number $\epsilon$ is so small that the parallelograms obtained from all grid 
intersection points are disjoint; see Figure \ref{fig:grid2tiling}
for an example. It is clear from this figure how one can glue the parallelograms of two consecutive intersection points on a grid line
along their sides perpendicular to the grid line.
The result of this glueing is then the same as the image (scaled with a factor 
$\epsilon$) of the surface $\cS$ under the projection $\RR^N\rightarrow\LL_\RR^\vee$.
\qed

\

It is obvious from the construction that the group $\LL$ acts by translations 
on the point set $\cS^0$ and on the surface $\cS$, preserving the types
of the squares and their sides. The cell structure on $\cS$ given by the squares, their sides and vertices induces therefore on $\modquot{\cS\,}{\LL}$ a cell structure. 
The following proposition counts the cells on
$\modquot{\cS\,}{\LL}$.

\

\ssnnl{ Proposition.}\label{count cells}
\textit{
Let $B$ be a $2\times N$-matrix such that its rows are a $\ZZ$-basis for
$\LL$. Let $\vb_1,\ldots,\vb_N$ be the columns of this matrix. Then
on $\modquot{\cS\,}{\LL}$
\begin{enumerate}
\item
the number of $2$-cells of type $(i,j)$ is $|\det(\vb_i,\vb_j)|$.
\item
the number of $1$-cells of type $i$ is $\sum_{j=1}^N|\det(\vb_i,\vb_j)|$.
\item
the number of $0$-cells is
$\sum_{1\leq i<j\leq N}|\det(\vb_i,\vb_j)|$.
\end{enumerate}
}
\proof
The counts can be made on the dual structure, which is given by the intersecting grid lines on $\LL_\RR$.
Let $B_{ij}$ denote the $2\times 2$-matrix with columns $\vb_i$ and $\vb_j$.
Then 1. amounts to counting for $\gl_i,\gl_j\in\RR$
the number of elements in the set
$\{\vx=(x_1,\,x_2)\in\RR^2\;|\; 
0\leq x_1,\,x_2<1\,,\;\vx\,B_{ij}+(\gl_i,\gl_j)\in\ZZ^2\}$.
Since there is a bijection between this set and the coset space $\modquot{\ZZ^2}{\ZZ^2B_{ij}}$, the number of elements is
$|\det B_{ij}|$.
The number of $1$-cells of type $i$ is equal
to the number of intersection points of the $i$-th grid with the other grids.
In view of 1. this number is therefore as stated in 2.
Knowing the numbers of $0$- and $1$-cells on the $2$-torus
$\modquot{\LL_\RR\,}{\LL}$ one computes the number of $2$-cells
from the fact that the Euler characteristic of a $2$-torus is $0$.
\qed

\

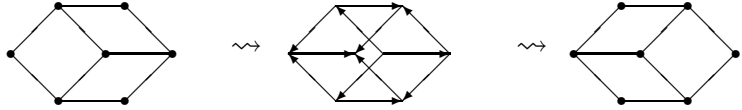
\begin{figure}[h]
\setlength{\unitlength}{.6pt}
\begin{picture}(300,90)(-20,-10)
\put(35,20){
\begin{picture}(100,80)(0,0)
\put(0,0){\line(1,1){30}}
\put(0,0){\line(1,-1){30}}
\put(60,0){\line(-1,1){30}}
\put(60,0){\line(-1,-1){30}}
\put(60,0){\line(1,0){42}}
\put(30,30){\line(1,0){42}}
\put(30,-30){\line(1,0){42}}
\put(102,0){\line(-1,1){30}}
\put(102,0){\line(-1,-1){30}}
\put(0,0){\circle*{5}}
\put(60,0){\circle*{5}}
\put(102,0){\circle*{5}}
\put(30,30){\circle*{5}}
\put(30,-30){\circle*{5}}
\put(72,30){\circle*{5}}
\put(72,-30){\circle*{5}}
\end{picture}}
\put(180,20){$\rightsquigarrow$}
\put(210,20){
\begin{picture}(100,80)(0,0)

\put(30,30){\vector(-1,-1){30}}
\put(30,-30){\vector(-1,1){30}}
\put(60,0){\vector(-1,-1){30}}
\put(60,0){\vector(-1,1){30}}
\put(60,0){\vector(1,0){42}}
\put(30,30){\vector(1,0){42}}
\put(30,-30){\vector(1,0){42}}
\put(102,0){\vector(-1,1){30}}
\put(102,0){\vector(-1,-1){30}}
\put(72,30){\vector(-1,-1){30}}
\put(72,-30){\vector(-1,1){30}}
\put(0,0){\vector(1,0){42}}

\end{picture}}
\put(360,20){$\rightsquigarrow$}
\put(390,20){
\begin{picture}(100,80)(0,0)
\put(0,0){\line(1,1){30}}
\put(0,0){\line(1,-1){30}}
\put(42,0){\line(1,1){30}}
\put(42,0){\line(1,-1){30}}
\put(42,0){\line(-1,0){42}}
\put(30,30){\line(1,0){42}}
\put(30,-30){\line(1,0){42}}
\put(102,0){\line(-1,1){30}}
\put(102,0){\line(-1,-1){30}}
\put(0,0){\circle*{5}}
\put(42,0){\circle*{5}}
\put(102,0){\circle*{5}}
\put(30,30){\circle*{5}}
\put(30,-30){\circle*{5}}
\put(72,30){\circle*{5}}
\put(72,-30){\circle*{5}}
\end{picture}}\end{picture}
\setlength{\unitlength}{1pt}
\caption{\label{fig:elementary transformation}
\textit{An elementary transformation}}
\end{figure}

\ssnnl{ Second step.}\label{deformS}
For the constructions in \ref{grids} and \ref{first step} we may randomly
choose a non-resonant $\gl\in\RR^N$. Once we have constructed 
the set $\cS^0$ and surface $\cS$ we modify it by \emph{elementary 
transformations}. For such an elementary transformation we need a point
$\vp_0\in\cS^0$ such that there are exactly three points 
$\vp_1,\,\vp_2,\,\vp_3$ in $\cS^0$ at distance $1$ from
$\vp_0$ and such that the three points $\vp_1+\vp_2-\vp_0$,
$\vp_1+\vp_3-\vp_0$ and $\vp_2+\vp_3-\vp_0$ also lie in $\cS^0$.
These seven points are the vertices of three squares in $\cS$.
The elementary transformation now replaces $\vp_0$ by 
$\vp'_0:=\vp_1+\vp_2+\vp_3-2\vp_0$ and replaces the three squares
by the three squares in $\RR^N$ with vertex configuration 
$
\{\vp_1,\,\vp_2,\,\vp_3,\,\vp_1+\vp_2-\vp_0,$ $
\vp_1+\vp_3-\vp_0,\,\vp_2+\vp_3-\vp_0,\,\vp'_0\}\,;
$
see Figure \ref{fig:elementary transformation}.
In order to preserve $\LL$-periodicity we perform this 
transformation simultaneously at all configurations
$\{\vp_0+\ell,$ $\vp_1+\ell,\,\vp_2+\ell,\,\vp_3+\ell\}$
with $\ell\in\LL$. 

After this elementary transformation we obtain a surface $\cS'$ in $\RR^N$
which is the union of all unit squares with vertices in the set 
$\cS'^0$, obtained from $\cS^0$ by replacing the points $\vp_0+\ell$
by $\vp_1+\vp_2+\vp_3-2\vp_0+\ell$ for all $\ell\in\LL$.

The surface $\cS'$, in turn, can be further modified by elementary transformations: 
just replace in the above construction
$\cS^0$ by $\cS'^0$ and $\vp_0$ by an appropriate point of $\cS'^0$.

\

\ssnnl{ Third step.}\label{specialS}
After starting the construction in \ref{first step} with a randomly chosen non-resonant $\gl\in\RR^N$ one can create by repeated elementary transformations
many surfaces $\cS$ each of which is the union of unit squares with vertices
in an $\LL$-invariant subset $\cS^0$ of $\ZZ^N$ and which is mapped homeomorphically onto $\LL_\RR^\vee$ by the projection
$\RR^N\rightarrow\LL_\RR^\vee$.
Moreover the numbers of cells on each of these surfaces are still the same as 
in Proposition \ref{count cells}. 

The number of surfaces one can make in this way is finite and it is
possible to generate (by computer) a complete list.

\

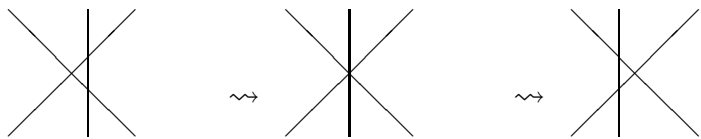
\begin{figure}[h]
\setlength{\unitlength}{.6pt}
\begin{picture}(300,90)(-50,10)
\put(35,20){
\begin{picture}(100,80)(0,0)
\put(0,0){\line(1,1){80}}
\put(80,0){\line(-1,1){80}}
\put(50,0){\line(0,1){80}}
\end{picture}}
\put(180,40){$\rightsquigarrow$}
\put(210,20){
\begin{picture}(100,80)(0,0)
\put(0,0){\line(1,1){80}}
\put(80,0){\line(-1,1){80}}
\put(40,0){\line(0,1){80}}
\end{picture}}
\put(360,40){$\rightsquigarrow$}
\put(390,20){
\begin{picture}(100,80)(0,0)
\put(0,0){\line(1,1){80}}
\put(80,0){\line(-1,1){80}}
\put(30,0){\line(0,1){80}}
\end{picture}}\end{picture}
\setlength{\unitlength}{1pt}
\caption{\label{fig:grid change}
\textit{Local change in grid at resonance}}
\end{figure}

\ssnnl{ Remark.}\label{start independent}
The list of surfaces produced in \ref{specialS} is, up to permutation of
its entries, independent of the choice of $\gl\in\RR^N$ at the start of the algorithm.
This can be seen as follows.
Fix $i$, $1\leq i\leq N$. Take $\gl\in\RR^N$ and let $\gl_t=\gl+t\ve_i$
for $t\in\RR$. Set 
$t_1=\min\{t\in\RR_{>0}\:|\:\gl_t\;\textrm{is not non-resonant}\}$
and $t_2=\min\{t\in\RR_{>t_1}\:|\:\gl_t\;\textrm{is not non-resonant}\}$.
Then the surface $\cS_t$ constructed from $\gl_t$ in 
\ref{first step} is equal to the surface $\cS_0$ for $0<t<t_1$.
Figure \ref{fig:grid change} shows how the grid locally changes as  $t$ passes through $t_1$.
Comparing Figures \ref{fig:grid change} and 
\ref{fig:elementary transformation} one sees that such a local change of the grid corresponds to an elementary transformation
of the surface $\cS_0$. Thus the surfaces
$\cS_t$ with $t_1<t<t_2$ are obtained from $\cS_0$ by a number of
elementary transformations.

Since any two non-resonant $\gl,\,\gl'\in\RR^N$ can be moved to a
common value by coordinatewise changes as above we conclude that
the surfaces $\cS$ and $\cS'$ obtained in \ref{first step}
from $\gl$ and $\gl'$, respectively, are related by 
elementary transformations, and that, hence, the
lists of surfaces which the algorithm produces from start values $\gl$ and
$\gl'$ are the same, up to possibly a permutation of the entries.

\

\ssnnl{ Definition.}\label{perfect surface}
We say that a surface $\cS$ on this list is 
\emph{perfect} if the function 
$\vs:\RR^N\rightarrow\RR$, $\vs(z_1,\ldots,z_N)=z_1+\ldots+z_N$,
takes only three values on the
set $\cS^0$ of vertices in $\cS$.
In that case these three values 
are consecutive integers, say $a+1,a,a-1$. We denote by
$\Sb$ (resp. $\Sd$ resp. $\Sc$) the set of vertices
where $\vs$ takes the value $a+1$ (resp. $a$ resp. $a-1$)
and say that the vertices in $\Sb$ are \emph{black}, those in $\Sc$ are 
\emph{white} and those in $\Sd$ are \emph{grey}.
Each of these sets is invariant under the translation action of $\LL$.

\

\ssnnl{ Remark.}\label{more perfect surfaces}
I have at present no proof that for every rank $2$ subgroup $\LL\subset\ZZ^N$
as in \ref{introL} the above construction indeed yields at least one perfect surface.
On the other hand, in all examples I investigated the computer produced
at least one perfect surface. It would also be very interesting to know
which perfect surfaces can arise already in the first step \ref{first step}
of the construction, by choosing $\gl$ appropriately (instead of at random).

\

\ssnnl{}\label{diagonals}
For every surface $\cS$ in the list of \ref{specialS} one has the graph 
$\widehat{\Gamma}$ with 
set of vertices $\cS^0$ and arrows given by the diagonals in the squares,
oriented such that in a square of type $(i,j)$
the oriented diagonals are $\ve_i+\ve_j$ and $\mathrm{sign}(\det(\vb_i,\vb_j))(\ve_i-\ve_j)$.
The advantage of having $\cS$ embedded in $\RR^N$ is that the vertices
and the arrows in these graphs are actual points and vectors in $\RR^N$.
If the surface $\cS$ is perfect, every square in $\cS$ has
one vertex in $\Sb$, one in $\Sc$ and
two vertices in $\Sd$. One of its diagonals
goes from the white to the black vertex and the other diagonal connects the
two  grey vertices. 
The graph $\widehat{\Gamma}$ therefore is the disjoint union of two 
oriented graphs
$\wG$ and $\wQ$, with vertex sets $\Sb\cup\Sc$ and $\Sd$, respectively.
The graph $\wG$ is a \emph{bi-partite graph}, i.e. its vertex set is the 
disjoint union of two sets (the black resp. white vertices) and with edges connecting only vertices of different colors. 

There is a natural duality between $\wG$ and $\wQ$: every vertex of
$\wG$ lies in a unique connected component of $\cS\setminus\wQ$
and vice versa. Moreover the boundaries of these connected components
are polygons and an arrow between two nodes in one graph is a common
side of the corresponding polygons for the other graph.

\

\ssnnl{}\label{periodic models 1}
The perfect surface $\cS$
is mapped homeomorphically onto $\LL_\RR^\vee$ by the projection
$\RR^N\rightarrow\LL_\RR^\vee$, which is dual to the inclusion 
$\LL\hookrightarrow\ZZ^N$. This projection maps $\ve_i$ to $\vb_i$, for
$i=1,\ldots,N$. One may however pre-compose this projection with
the linear map $\RR^N\rightarrow\RR^N,\,\ve_i\mapsto
\frac{1}{\parallel\vb_i\parallel}\ve_i,\,1\leq i\leq N$. The composite linear map
$\RR^N\rightarrow\LL_\RR^\vee$ maps the tiling by squares on the surface $\cS$
piecewise linearly and homeomorphically onto a tiling of $\LL_\RR^\vee$
by rhombi. It maps $\LL$ isomorphically onto a lattice $\oLL$ in $\LL_\RR^\vee$
and it maps the graphs $\wG$ and $\wQ$  `isomorphically'
onto $\oLL$-periodic graphs $\oG$ and $\oQ$ in the plane $\LL_\RR^\vee$.
In the literature  $\oG\subset\LL_\RR^\vee$ is called  
a \emph{periodic dimer model} or
\emph{periodic brane tiling}; see for instance \cite{FHKV,FHKVW,HHV,HK,HV,Ky}.

\

\ssnnl{ Definition.}\label{graphs}
All structures in \ref{diagonals} are invariant under translation by vectors in $\LL$. Passing to the orbit space we obtain from a perfect surface $\cS$
the $2$-dimensional torus $\cT=\modquot{\cS\,}{\LL}$
(i.e. compact oriented surface of genus $1$ without boundary). Embedded in this torus are the two graphs
$\Gg=\modquot{\wG}{\LL}$ and $\Qg=\modquot{\wQ}{\LL}$,
which are dual to each other.
Note that \ref{periodic models 1} also yields
$\cT=\modquot{\LL_\RR^\vee\,}{\oLL}$, $\Gg=\modquot{\oG}{\oLL}$ and $\Qg=\modquot{\oQ}{\oLL}$.

After Grothendieck one calls $(\cT,\Qg,\Gg)$  a \emph{dessin d'enfants} \cite{LZ, Des}. 

\

\begin{figure}[h]
\begin{picture}(350,110)(-40,10)
\put(-40,10){
\begin{picture}(150,150)(0,0)
\setlength\epsfxsize{4cm}
\epsfbox{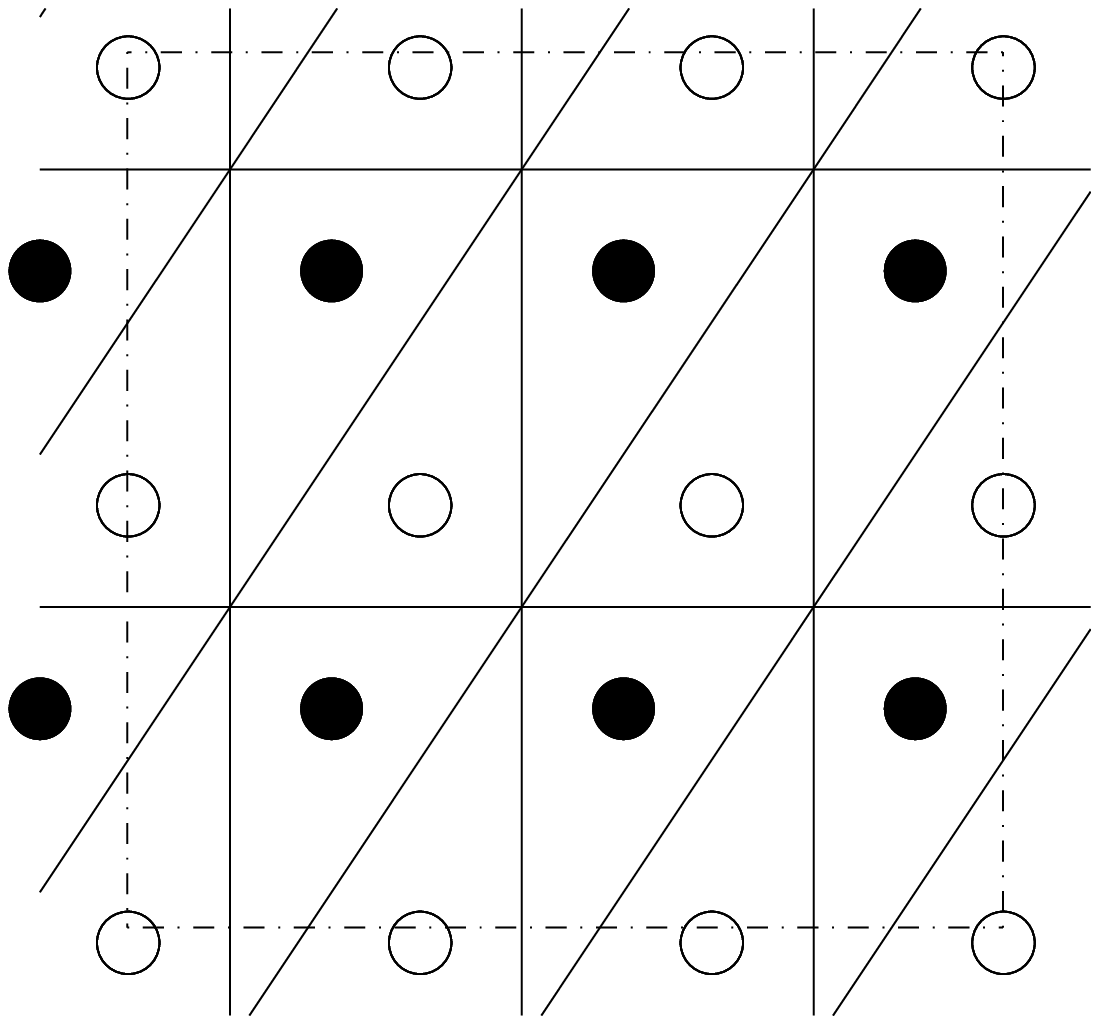}
\end{picture}
}
\put(100,10){
\begin{picture}(150,150)(0,0)
\setlength\epsfxsize{4cm}
\epsfbox{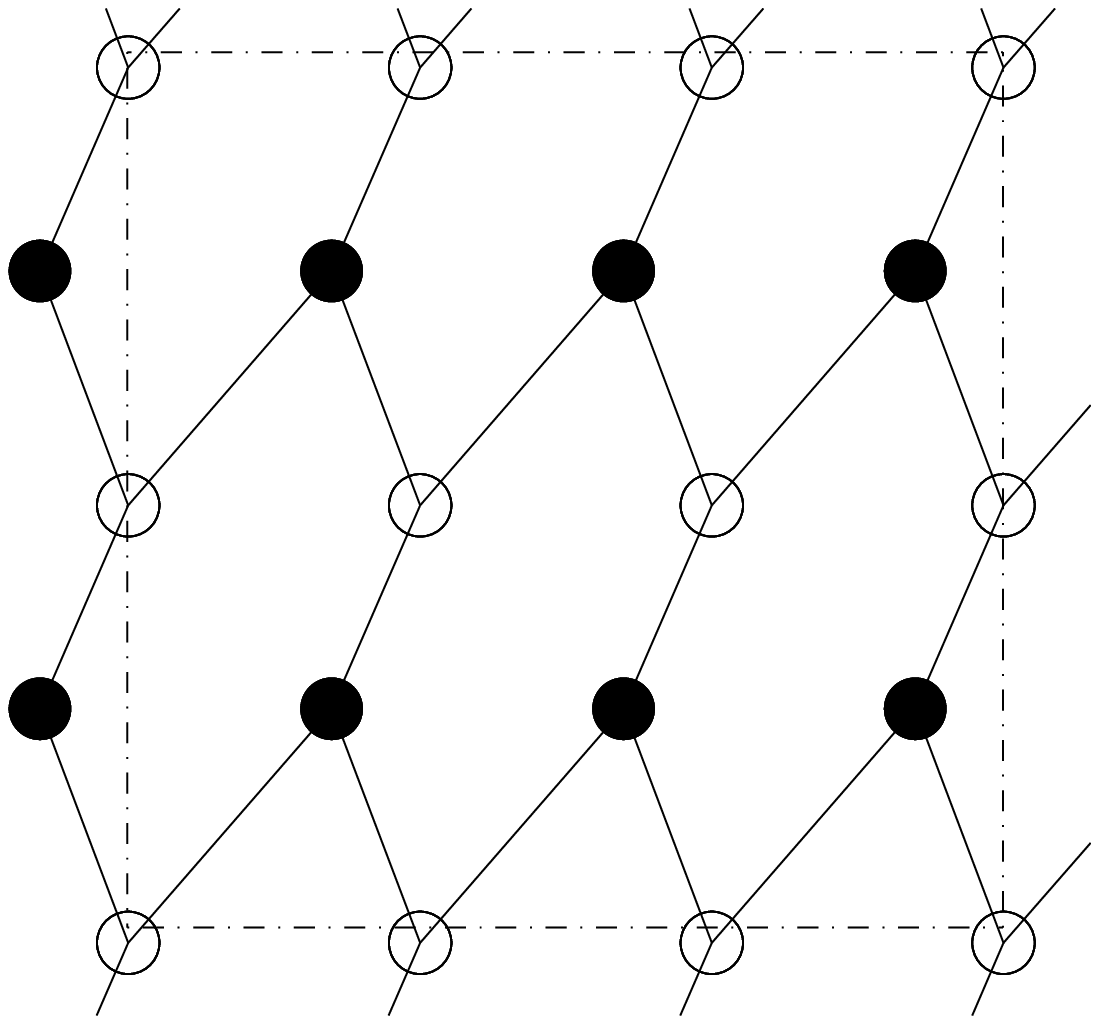}
\end{picture}
}
\put(250,50){
\begin{picture}(100,80)(0,0)
\put(0,0){\circle*{5}}
\put(0,30){\circle*{5}}
\put(30,-15){\circle*{5}}
\put(60,0){\circle*{5}}
\put(60,30){\circle*{5}}
\put(30,45){\circle*{5}}
\put(0,30){\vector(0,-1){27}}
\put(0,0){\vector(2,-1){27}}
\put(30,-15){\vector(2,1){27}}
\put(60,0){\vector(0,1){27}}
\put(60,30){\vector(-2,1){27}}
\put(30,45){\vector(-2,-1){27}}
\put(0,0){\vector(1,0){57}}
\put(60,0){\vector(-2,3){28}}
\put(30,45){\vector(-2,-3){28}}
\put(60,30){\vector(-1,0){57}}
\put(0,30){\vector(2,-3){28}}
\put(30,-15){\vector(2,3){28}}
\put(30,43){\vector(0,-1){55}}
\put(30,-13){\vector(0,1){55}}
\put(2,29){\vector(2,-1){55}}
\put(58,1){\vector(-2,1){55}}
\put(2,1){\vector(2,1){55}}
\put(58,29){\vector(-2,-1){55}}
\end{picture}
}
\end{picture}
\caption{\label{fig:dP3IVtwist}
\textit{Two versions of the dessin $(\cT,\Qg,\Gg)$ for $B_{10}$ from Figure \ref{fig: models}. The dashed square is the period parallelogram.
The picture on the left shows how the quiver on the right is embedded in a $2$-torus.}}
\end{figure}

\ssnnl{ Example.}\label{dP3IVtwist}
Figure \ref{fig:dP3IVtwist} shows two drawings of the dessin $(\cT,\Qg,\Gg)$ for $B_{10}$ from Figure \ref{fig: models}; or, rather, it shows a lifting of this dessin to the plane $\LL_\RR^\vee$ with basis so that the period lattice $\oLL$ becomes just $\ZZ^2$.
It also shows $\Qg$ as an unembedded quiver.
\emph{Note however that this is not the same as the quiver $\sQ$
in case $B_{10}$ in Figure \ref{fig: models}.}

The quiver in Figure \ref{fig:dP3IVtwist} is in fact the quiver for the
singularity $\modquot{\CC^3}{\ZZ_6}$ where a generator of the cyclic 
group  $\ZZ_6$ acts on $\CC^3$ as multiplication by the diagonal matrix
$\mathrm{diag}(e^{2\pi i/6},\,e^{2\pi i/3},\,e^{2\pi i/2})$.
Also the brane tiling picture in Figure \ref{fig:dP3IVtwist} appears
in \cite{HK} \S 3 in connection with the
singularity $\modquot{\CC^3}{\ZZ_6}$.

This difference between \cite{HK} and our approach comes, because
we still have to perform the untwist, as is explained below.

\

\ssnnl{}\label{zigzagquiver}
The quiver $\Qg$ in \ref{graphs} is not the Pl\"{u}cker quiver
of $\LL$, defined in \ref{pluckerquiver}. The latter appears in the current setting through the \emph{zigzag loops}, as follows.
Let $\cS$ be a perfect surface. Fix $i\in\{1,\ldots, N\}$. In every square
on $\cS$ which has two sides parallel to $\ve_i$ draw the line segment
connecting the midpoints of these two sides. If the square has type $(i,j)$
this line segment is oriented so that it is a translate of the vector 
$\mathrm{sign}(\det(\vb_i,\vb_j))\ve_j$.
The union of these line segments projects to an oriented closed
curve on $\cT$, called the 
\emph{$i$-th zigzag loop}. Note that this zigzag loop may consist of
several connected components.

According to Proposition \ref{count cells}
the $i$-th and $j$-th zigzag loops on $\cT$ intersect in exactly
$|\det(\vb_i,\vb_j)|$ points.
\emph{Thus we see that the nodes of the Pl\"{u}cker quiver can be identified with
the zigzag loops on $\cT$ and the arrows between two nodes are 
the intersection points of the corresponding zigzag loops.}
The orientation of the arrows follows from the orientation of the zigzag loops
and can in pictures be indicated as an over/undercrossing of the zigzag loops.

\begin{figure}[h]
\begin{picture}(350,100)(-40,20)
\put(40,10){
\begin{picture}(150,150)(0,0)
\setlength\epsfxsize{7cm}
\epsfbox{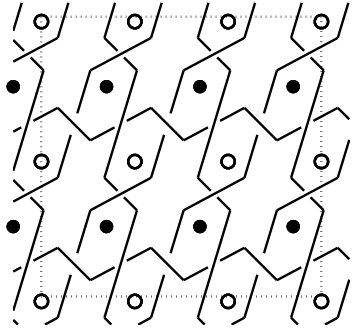}
\end{picture}
}
\end{picture}
\caption{\label{fig:dP3IVzigzag}
\textit{The zigzag loops for the dessin in Figure \ref{fig:dP3IVtwist}}}
\end{figure}

Figure \ref{fig:dP3IVzigzag} shows the zigzag loops for the dessin in Figure \ref{fig:dP3IVtwist}.

\

\begin{figure}[h]
\setlength{\unitlength}{0.5pt}
\begin{picture}(350,200)(-100,-20)
\thicklines
\put(0,0){
\begin{picture}(150,100)(0,0)
\put(0,0){\line(1,1){80}}
\put(0,0){\line(-1,1){80}}
\put(0,160){\line(-1,-1){80}}
\put(0,160){\line(1,-1){80}}
\put(0,0){\circle{15}}
\put(0,160){\circle*{15}}
\put(0,0){\vector(1,1){30}}
\put(0,0){\vector(-1,1){30}}
\put(-80,80){\vector(1,1){30}}
\put(80,80){\vector(-1,1){30}}
\put(-10,20){\vector(1,0){20}}
\put(-10,140){\vector(1,0){20}}
\put(-12,175){$\vb$}
\put(-12,-25){$\vw$}
\put(45,120){$\vr$}
\put(-60,120){$\vr'$}
\put(45,34){$\vr'$}
\put(-60,34){$\vr$}

\end{picture}
}
\put(200,0){
\begin{picture}(150,100)(0,0)
\put(0,0){\line(1,1){40}}
\put(0,0){\line(-1,1){40}}
\put(0,160){\line(-1,-1){40}}
\put(0,160){\line(1,-1){40}}
\put(-40,40){\line(1,1){80}}
\put(40,40){\line(-1,1){35}}
\put(-40,120){\vector(1,-1){35}}
\put(-40,40){\vector(1,1){35}}
\put(0,0){\vector(1,1){30}}
\put(0,0){\vector(-1,1){30}}
\put(-40,120){\vector(1,1){30}}
\put(40,120){\vector(-1,1){30}}
\put(0,0){\circle{15}}
\put(0,160){\circle*{15}}
\put(-10,20){\vector(1,0){20}}
\put(-10,140){\vector(1,0){20}}
\put(-12,175){$\vb$}
\put(-12,-25){$\vw$}
\put(45,120){$\vr$}
\put(-60,120){$\vr'$}
\put(45,34){$\vr'$}
\put(-60,34){$\vr$}

\end{picture}
}

\put(360,0){
\begin{picture}(150,100)(0,0)
\put(0,0){\line(1,1){40}}
\put(0,0){\line(-1,1){40}}
\put(0,160){\line(-1,-1){40}}
\put(0,160){\line(1,-1){40}}
\put(40,40){\line(0,1){80}}
\put(-40,40){\line(0,1){80}}
\put(40,40){\vector(0,1){40}}
\put(-40,120){\vector(0,-1){60}}
\put(0,0){\vector(1,1){30}}
\put(0,0){\vector(-1,1){30}}
\put(-40,120){\vector(1,1){30}}
\put(40,120){\vector(-1,1){30}}
\put(0,0){\circle{15}}
\put(0,160){\circle*{15}}
\put(10,20){\vector(-1,0){20}}
\put(-10,140){\vector(1,0){20}}
\put(-12,175){$\vb$}
\put(-12,-25){$\vw$}
\put(45,80){$\vr$}
\put(-60,80){$\vr'$}
\end{picture}
}
\put(560,0){
\begin{picture}(150,100)(0,0)
\put(0,0){\line(1,1){80}}
\put(0,0){\line(-1,1){80}}
\put(0,160){\line(-1,-1){80}}
\put(0,160){\line(1,-1){80}}
\put(0,0){\circle{15}}
\put(0,160){\circle*{15}}
\put(10,20){\vector(-1,0){20}}
\put(-10,140){\vector(1,0){20}}
\put(-60,90){\vector(0,-1){20}}
\put(60,70){\vector(0,1){20}}

\put(-40,40){\line(1,1){80}}
\put(40,40){\line(-1,1){80}}
\put(0,80){\vector(1,1){20}}
\put(0,80){\vector(-1,-1){20}}
\put(40,40){\vector(-1,1){20}}
\put(-40,120){\vector(1,-1){20}}
\put(-12,175){$\vb$}
\put(-12,-25){$\vw$}
\put(85,80){$\vr$}
\put(-100,80){$\vr'$}
\end{picture}
}
\end{picture}
\setlength{\unitlength}{1pt}
\caption{\label{fig:untwist helix}
Untwisting:  \textit{from left to right:
original square tile, helix tile, untwisted helix tile, untwisted helix 
completed to square tile.}}
\end{figure}
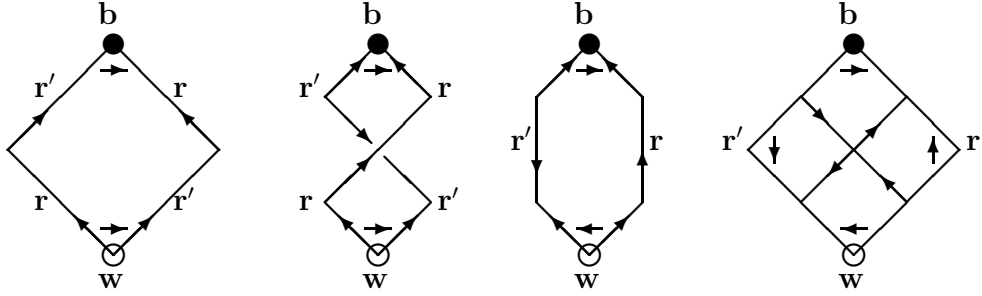

\ssnnl{ Final untwist.}\label{twist}
Pictures like Figure \ref{fig:dP3IVzigzag} can be viewed as showing an oriented
surface with boundary, in which the connected components of the zigzag
loops are the connected components of the boundary and which is embedded in
three space in a twisted way so that the black dots are on one side
of the surface and the white dots are on the other side. 
This surface comes with a tiling by helices as shown
in Figure \ref{fig:untwist helix}.

One can `untwist' the helix tiles and complete the untwisted helix tiles to squares as indicated in Figure \ref{fig:untwist helix}. 
The effect is that
the boundary cycles are being capped off by discs. The result is an oriented surface $\RS$ without boundary which is tiled by squares with one black, one white and two red vertices. The black-white diagonals give an embeding of the bipartite graph $\Gg$ into $\RS$. The red-red diagonals give an embeding of the quiver $\sQ$ into $\RS$ (if all zigzag loops have just one connected component; otherwise some points on the surface must be pinched
together). 

The untwisting procedure is described in \cite{FHKV} \S 5 via examples, 
pictures and visual inspection. Here we want to build the `untwisted
surface' by simply reinterpreting the combinatorial data of the tiling
of the torus $\cT$, i.e. the
$\LL$-periodic tiling of the perfect surface $\cS$ by unit squares
taken modulo $\LL$
(see \ref{perfect surface}, \ref{graphs}).
Let $\sB$ (resp. $\sW$) denote the set of black (resp. white) vertices.
There is a bijection between the set of tiles and the set $E$ of arrows
of the quiver $\sQ$ and we can refer to tiles as $e\in E$.
A tile $e$ has one black vertex $\vb(e)$ and one white vertex $\vw(e)$.
A tile $e$ is the image of a unit square with sides parallel to two of the
basis vectors $\ve_1,\ldots,\ve_N$; let us denote the indices of this (unordered) pair of basis vectors as a $2$-element subset 
$\{\vr(e), \,\vr'(e)\}$ of $\{1,\ldots,N\}$.
Thus the tiling on $\cT$ yields the list of quadruples 
\begin{equation}\label{eq:quadruples}
\MM\,=\,
\{\:(\,\vb(e),\,\vw(e),\,\vr(e),\,\vr'(e)\,)\:\}_{e\in E}\,.
\end{equation}
Now, for every $e\in E$ we take a unit square $\square_e$ and attach labels
$\vb(e)$, $\vw(e)$, $\vr(e)$, $\vr'(e)$ and colors black, white, red, red,
respectively, to its vertices such that the vertices labeled
$\vb(e)$ and $\vw(e)$ are not adjacent. If two such squares $\square_e$ and
$\square_{e'}$ have equal labels on two adjacent vertices, we glue
$\square_e$ and $\square_{e'}$ along the corresponding sides.
The result of this glueing is a surface $\wRS$ tiled with squares.

On the surface $\wRS$ there is for every $\vb\in\sB$ one black point 
with label $\vb$ and for every $\vw\in\sW$ there is one white point 
with label $\vw$. However for $i\in\{1,\ldots,N\}$
there may be several red points with label $i$.
We finally obtain the desired surface $\RS$ by identifying for
every $i\in\{1,\ldots,N\}$ the red points with label $i$.

\

\ssnnl{ Proposition.}\label{oriented M}
\textit{The surface $\wRS$, and hence also $\RS$, is oriented.}

\proof Cut the square $\square_e$ in four pieces like in the right-hand picture 
of Figure \ref{fig:untwist helix}.
Then for fixed $\vb\in\sB$ the small squares with black vertex
$\vb$ are glued together to a polygon isomorphic
to the polygon in $\cT$ formed by the unit squares with one vertex $\vb$.
Similarly, for fixed $\vw\in\sW$ the small squares with white vertex
$\vw$ are glued together to a polygon isomorphic
to the polygon in $\cT$ formed by the unit squares with one vertex $\vw$
but with the orientation reversed.
For $i\in\{1,\ldots,N\}$ the small squares with red vertex
$i$ are glued together to a number 
of disjoint polygons, one for every connected component of the $i$-th zigzag loop and with oriented boundary `equal to' that component.
It is now clear from Figure \ref{fig:untwist helix} that these polygons are oriented consistently so that $\wRS$ is oriented.
\qed 

\

\ssnnl{ Conclusion.}\label{final M}
\textit{The surface $\RS$ is tiled with squares so that the black-white diagonals
form a bi-partite graph isomorphic to the bi-partite graph $\Gg$
in \ref{graphs}, while the red-red diagonals
form a graph isomorphic to the quiver $\sQ$.}

\

\begin{figure}[h]
\begin{picture}(350,150)(-40,-10)
\put(0,0){
\begin{picture}(150,100)(0,15)
\setlength\epsfxsize{5cm}
\epsfbox{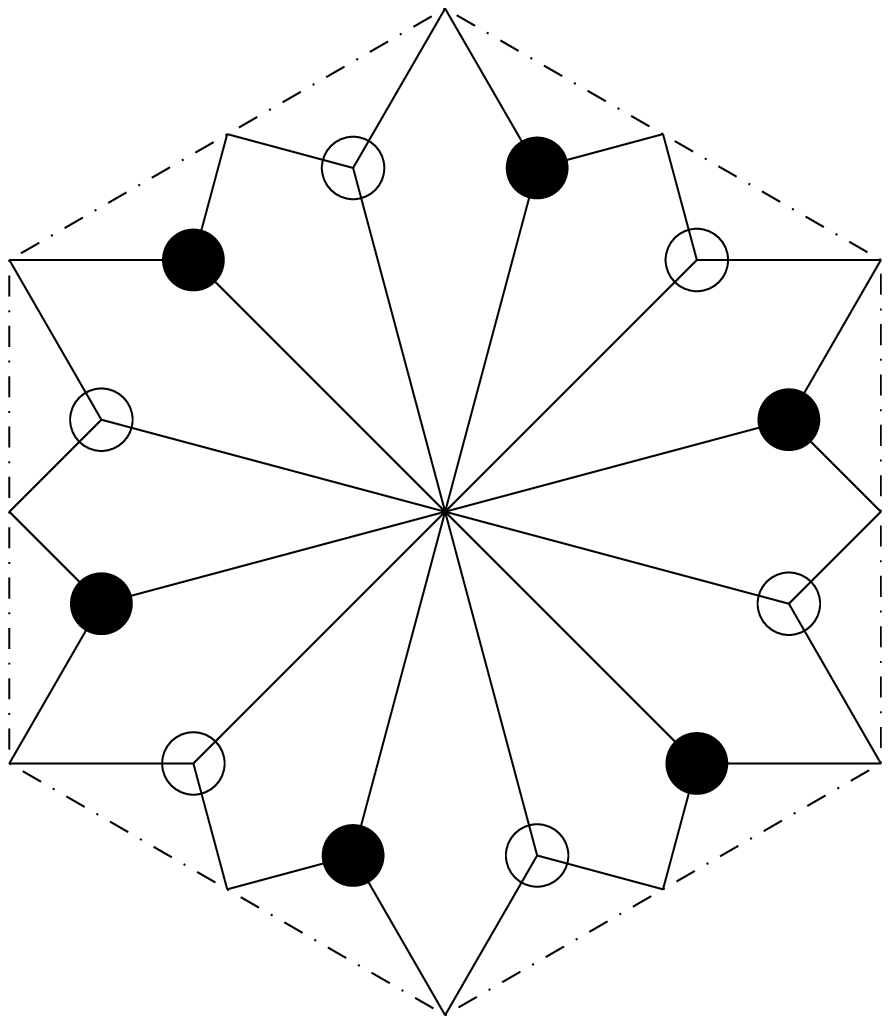}
\end{picture}
}

\put(180,0){
\begin{picture}(150,100)(0,15)
\setlength\epsfxsize{5cm}
\epsfbox{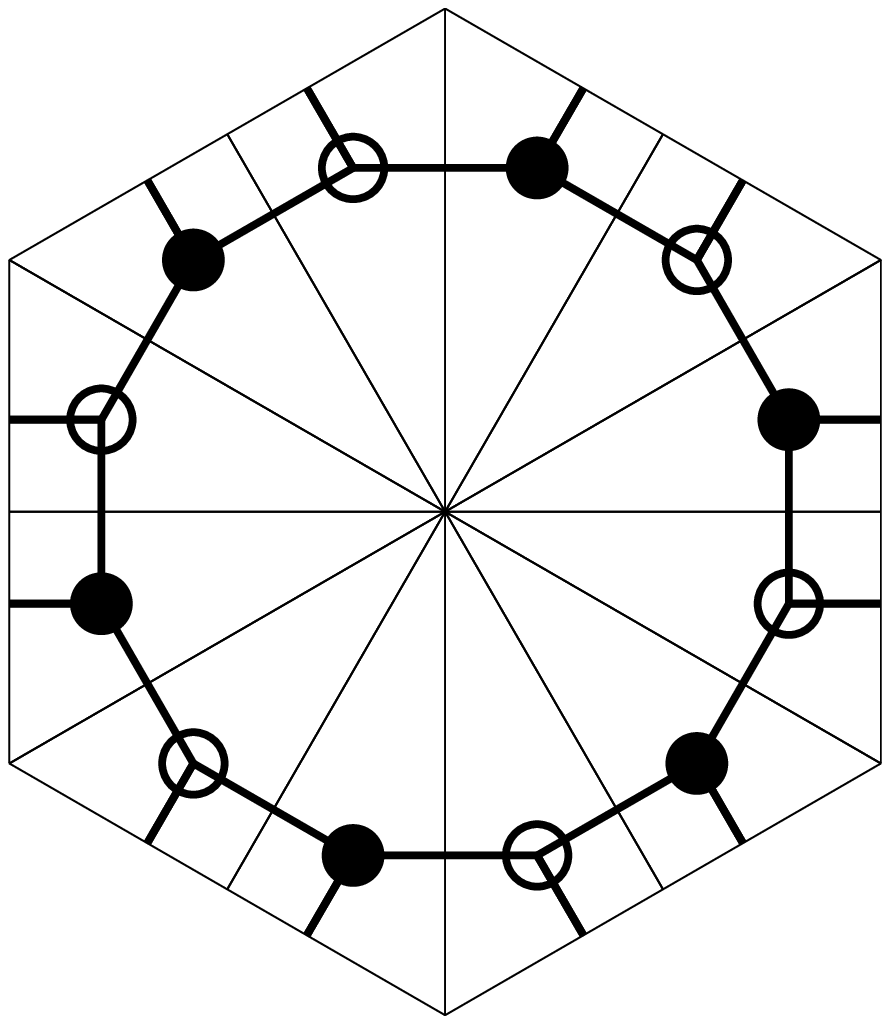}
\end{picture}
}

\end{picture}
\caption{\label{fig:dP3IV}
\textit{Surface $\RS$ with kite tiling (left) and quiver $\sQ$ and bipartite graph $\Gg$ (right)
for case  $B_{10}$ in Figure \ref{fig:dp3 models}
(= model IV of $dP_3$). The surface $\RS$ is obtained by identifying opposite
sides of the hexagon.}}
\end{figure}

\ssnnl{ Example.}\label{correct dP3IV}
The list $\MM$ (\ref{eq:quadruples}) which our algorithm produces for 
model IV of $dP_3$, i.e.
example $B_{10}$ from Figure \ref{fig: models}, is
$$
\begin{array}{rcrrrrrrrrrrrrrrrrrr}
      e&:&1&2&3&4&5&6&7&8&9&10&11&12&13&14&15&16&17&18\\
 \vb(e)&:&6&5&1&2&4&3&1&5&3&6&4&2&5&6&1&2&3&4\\
 \vw(e)&:&5&6&2&1&3&4&3&1&5&2&6&4&5&6&1&2&3&4\\
 \vr(e)&:&2&2&3&3&4&4&1&1&1&1&1&1&5&6&5&6&5&6\\
\vr'(e)&:&1&1&1&1&1&1&5&5&5&6&6&6&2&2&3&3&4&4
\end{array}
$$
The surface $\RS$ with the corresponding tiling by kites 
(instead of squares) is shown in the left-hand picture in
Figure \ref{fig:dP3IV}. The right-hand picture in
Figure \ref{fig:dP3IV} shows the bi-partite graph $\Gg$ and the
quiver $\sQ$ embedded in $\RS$. Note that this quiver is indeed the same
as the one for example $B_{10}$ in Figure \ref{fig:dP3IV}.

\section{Perfect matchings and the secondary fan.}
\label{sec:perfect sec fan}

\ssnnl{ Definition.}\label{def:perfect matching}
(see e.g. \cite{KOS, Kn, ORV})
A \emph{perfect matching} $P$ on a bi-partite graph  $\Gg$ (also known as a 
\emph{dimer configuration})
is a subset of the edges of 
$\Gg$ such that every node of $\Gg$ is incident to exactly one
edge in $P$. 

\

\ssnnl{ Theorem.}\label{thm:perfect and secondary fan}
\textit{Consider a subgroup $\LL\subset\ZZ^N$ as in \ref{introL}.
Let $\Gg$ be a bi-partite graph which 
is obtained from $\LL\subset\ZZ^N$ as explained in Section \ref{sec:dessin}, and in Definition \ref{graphs} in particular.
Let $C$ be a $2$-dimensional cone in the secondary fan of
$\,\LL\subset\ZZ^N$ and let $L_C$ be the corresponding collection of
$2$-element subsets of $\{1,\ldots,N\}$ as in 
\ref{def:Lsecondary fan} and \ref{def:Lsecondary polytope}.
Identify $E$ with the set of edges of $\Gg$.}
\\
\textit{Then, with the notation as in (\ref{eq:quadruples}), the set
\begin{equation}\label{eq:PCLC}
P_C\,=\,\{e\in E\:|\: \{\vr(e),\vr'(e)\}\,\in\,L_C\:\}
\end{equation}
is a perfect matching on $\Gg$.}

\proof
Consider a perfect surface $\cS\subset\RR^N$ as in \ref{perfect surface}.
The map $\pi:\RR^N\rightarrow\LL_\RR^\vee$ which is dual to the inclusion
$\LL\hookrightarrow\ZZ^N$ projects $\cS$ homeomorphically onto the plane
$\LL_\RR^\vee$. The tiling of $\cS$ by unit squares gives a tiling
of $\LL_\RR^\vee$ by parallelograms. A unit square with edges
$\ve_i$ and $\ve_j$ projects onto a parallelogram with edges
$\vb_i$ and $\vb_j$. Now let $C$ be a $2$-dimensional cone of the secondary 
fan (see \ref{def:Lsecondary fan}) and let $\vw$ be a white vertex of the 
tiling. By translating the origin of the secondary fan to $\vw$ one sees that 
there is exactly one parallelogram with vertex $\vw$
and with edges $\vb_i$ and $\vb_j$ such that
$C\,\subset\,\RR_{\geq 0}\vb_i\,+\,\RR_{\geq 0}\vb_j$, i.e. such that 
$\{i,j\}\in L_C$. This shows that there is exactly one element of $P_C$ incident to $\vw$, namely the black-white diagonal of this parallelogram.
A similar argument works for the black vertices.
\qed

\

\ssnnl{ Corollary.}\label{black volA}
\textit{In the situation of Theorem \ref{thm:perfect and secondary fan}
the number of black vertices and the number of white vertices
are both equal to the 
number $\mathrm{vol}_\cA$ in (\ref{eq:volA}).}

\proof
The number of black (resp. white) vertices of $\Gg$ is equal to the 
number of edges in any perfect matching. This holds in particular for
the perfect matching $P_C$ corresponding to a $2$-dimensional cone $C$ in the secondary fan. The result now follows from
(\ref{eq:PCLC}), (\ref{eq:LC}) and (\ref{eq:volA}).
\qed

\

The following result was also derived in \cite{FHKV} \S 5.1.

\

\ssnnl{ Corollary.}\label{genus}
\textit{The genus of the surface $\RS$ in \ref{twist} is equal to the number of lattice points in the interior of the secondary polytope $\Delta$
in (\ref{eq:easy secondary}).}

\proof The number of $1$-cells in the tiling of $\RS$ is twice the number of $2$-cells. According to Proposition \ref{count cells} the number of $2$-cells is $\sum_{1\leq i<j\leq N}|\det(\vb_i,\vb_j)|$. The number of $0$-cells is,
according to Corollary \ref{black volA} and Formula (\ref{eq:volA}) equal to
$N+2\sum_{\{i,j\}\in L_C}|\det (\vb_i,\vb_j)|$. The desired result now follows from Corollary \ref{internal secondary} and a calculation of the
Euler characteristic  $2-2\,\mathrm{genus}(\RS)$.
\qed

\section{Kasteleyn matrix, bi-adjacency matrix,\\
$3$-constellation and superpotential.}
\label{sec:Kasteleyn superpot}
\ssnnl{}\label{traditional Kasteleyn matrix}
In \ref{periodic models 1} we found the
$\oLL$-periodic bi-partite graph $\oG$ in the plane $\LL_\RR^\vee$.
Periodic bi-partite graphs in the plane are also known as  
\emph{periodic dimer models}.
One of the main tools for studying periodic dimer models is the \emph{Kasteleyn 
matrix}. It is usually defined as an adjacency 
matrix with extra factors $x,x^{-1},y,y^{-1}$ for edges crossing 
the sides of a fundamental parallelogram
and determined by visual inspection of
a picture of the periodic dimer model; see e.g. \cite {KOS,HK,Ky}. 
The purpose of the extra factors
is to keep track of the homology class of a closed path
on the bi-partite graph $\Gg$ in the torus $\cT$ (see \ref{graphs}).
We obtained $\oG\,\subset\,\LL_\RR^\vee$ as the projection of the
black-white diagonals in the tiling of the perfect surface $\cS$
by unit squares. The advantage of having $\cS$ embedded in $\RR^N$ is that 
these diagonals are actual vectors in $\ZZ^N$ and that points lying on paths
formed by these diagonals are simply given by vector addition.
Thus all edges become ``equally responsible'' for the periodicity.
An edge which is the diagonal in a square of type $(i,j)$ contributes
$\ve_i+\ve_j$.

Thus we are led to a reformulation of the Kasteleyn matrix of the dimer model
$(\cT,\Gg)$ in terms of the combinatorial data given as the list
$\MM$ in Equation (\ref{eq:quadruples}). As this same list
is used in \ref{twist} for building the surface $\RS$ we call the same matrix
also the bi-adjacency matrix of the dessin $(\RS,\sQ)$.

\
 
\ssnnl{ Definition.}\label{def:Kasteleyn matrix} 
From the list  
$\MM=\{\:(\,\vb(e),\,\vw(e),\,\vr(e),\,\vr'(e)\,)\:\}_{e\in E}$ in
(\ref{eq:quadruples}) and a function $\wt:E\rightarrow\CC$
on the set of edges of the quiver $\sQ$
we construct a matrix $\KK^\wt$ with rows indexed by the set of black vertices 
$\sB$, columns indexed by the set of white vertices $\sW$
and with entries in the polynomial ring $\CC[u_1,\ldots,u_N]$:
\begin{equation}\label{eq:Kasteleyn matrix}
(\vb,\vw)\textrm{-entry of}\;\KK^\wt\quad=
\sum_{e\in E:\,\vb(e)=\vb,\,\vw(e)=\vw}
\wt (e)\,u_{\vr(e)}u_{\vr'(e)}\,.
\end{equation}
We use for $\KK^\wt$ the name \emph{Kasteleyn matrix of the dimer model $(\cT,\Gg)$}
as well as the name \emph{bi-adjacency matrix of the dessin $(\RS,\sQ)$}.

\

\ssnnl{ Example.}\label{dP3IV Kasteleyn}
Although the relation between the list $\MM$ and the matrix $\KK^\wt$ is so straightforward, that there seems no point in writing here the
matrix $\KK^\wt$ for Example \ref{correct dP3IV}, we do present
this $\KK^\wt$ here to compare it with results in the literature.
$$
\textrm{bi-adjacency matrix }\;\KK^\wt\quad\textrm{for model IV of}\; dP_3\;=\;
$$
$$
\left[
\begin{array}{cccccc}
\wt_{15}u_3u_5&\wt_3u_1u_3&\wt_7u_1u_5&0&0&0\\
\wt_4u_1u_3&\wt_{16}u_3u_6&0&\wt_{12}u_1u_6&0&0\\
0&0&\wt_{17}u_4u_5&\wt_6u_1u_4&\wt_9u_1u_5&0\\
0&0&\wt_5u_1u_4&\wt_{18}u_4u_6&0&\wt_{11}u_1u_6\\
\wt_8u_1u_5&0&0&0&\wt_{13}u_2u_5&\wt_2u_1u_2\\
0&\wt_{10}u_1u_6&0&0&\wt_1u_1u_2&\wt_{14}u_2u_6
\end{array}
\right]
$$
(with $\wt_e$ denoting $\wt(e)$).
In view of the discussion in Example \ref{dP3IVtwist} this
bi-adjacency matrix for model IV of $dP_3$ should be somehow the same as the Kasteleyn matrix for the singularity $\modquot{\CC^3}{\ZZ_6}$ in
\cite{HK} Equation (3.5):
$$
K(z,w)=
\left[
\begin{array}{cccccc}
w&-1&0&-w&0&0\\
0&1&-1&0&-1&0\\
-1&0&1&0&0&-1\\
-zw&0&0&w&-1&0\\
0&-z&0&0&1&-1\\
0&0&-z&-1&0&1
\end{array}
\right]\,.
$$
Indeed, by dividing the first, third and fifth column of $\KK^\wt$ by $u_5$
and subsequently setting $u_5=z^{-1}$, $u_3=w$, $u_1=u_2=u_4=u_6=1$,
and appropriately setting $\wt_e=\pm 1$ one obtains a matrix which up to a permutation of rows and columns is $K(z,w)$. This way of getting $K(z,w)$
from $\KK^\wt$ is nothing but passing from homogeneous to
inhomogeneous coordinates.

\

\ssnnl{ Definition.}\label{def:constellation}(cf. \cite{LZ} Definition 1.1.1.) 
A \emph{$3$-constellation} $(E,\gs_0,\gs_1)$ consists of a finite set $E$ 
and two permutations $\gs_0,\gs_1$ of $E$ such that the permutation group 
generated by $\gs_0,\gs_1$ acts transitively on $E$.

\

\ssnnl{ Notations.}\label{cycle notations}
The cycle notation $(i_1\;i_2\ldots i_{k-1}\;i_k)$ denotes the permutation $\rho$ such that $\rho(i_j)=i_{j+1}$ for $j=1,\ldots,k-1$ and
$\rho(i_k)=i_1$.
By a \emph{cycle} of a permutation $\gs$ we mean an orbit of the group
generated by $\gs$. We denote the set of cycles of $\gs$  
by $E_{\gs}$.
In the product $\gs\tau$ of two permutations
$\gs$ and $\tau$ one first applies
$\tau$ ; so $(\gs\tau) (i)=\gs(\tau (i))$.

\

\ssnnl{ Definition.}\label{def:superpotential}
For our purposes we do not need the most general notion of a superpotential
for a quiver. The superpotentials we need are just notational reformulations of 
$3$-constellations. The \emph{superpotential attached to a
$3$-constellation} $(E,\gs_0,\gs_1)$ is the following polynomial $W$
in non-commuting variables $X_e$ ($e\in E$):
$$
W=\sum_{\gamma\in E_{\gs_0}} \sX_\gamma\:-\:
\sum_{\gamma\in E_{\gs_1}} \sX_\gamma\,,
$$
where $\sX_\gamma\,:=\, X_{i_1}X_{i_2}\cdot\ldots\cdot X_{i_k}$
for a cyclic permutation $\gamma=(i_1\;i_2\ldots i_k)$.

\

\ssnnl{ Definition.}\label{def:dessin constellation}
The \emph{$3$-constellation $(E,\gs_0,\gs_1)$ associated
with the dessin $(\RS,\sQ)$} in \ref{final M} is
the following. $E$ is the set of arrows of the quiver $\sQ$.
Every $e\in E$ is a path $p_e$ on $\RS$.
Every connected component of $\RS\setminus \cup_{e\in E}\: p_e$
has an oriented boundary which can be viewed as a cyclic permutation of elements of $E$. Then $\gs_0$ (resp. $\gs_1$) is the composition of the cyclic permutations which are boundaries of connected components containing
a black (resp. white) point.

\

\ssnnl{ Remark.}\label{def:twist constellation}
In the same way one associates a $3$-constellation with the dessin 
$(\cT,\Qg)$ in \ref{graphs}. It is obvious from the untwisting procedure 
in \ref{twist} that, if $(E,\gs_0,\gs_1)$ is the $3$-constellation associated
with the dessin $(\RS,\sQ)$, then $(E,\gs_0,\gs_1^{-1})$ is the 
$3$-constellation associated with the dessin $(\cT,\Qg)$.

\

\ssnnl{}\label{M from constellation}
It is obvious that the $3$-constellation $(E,\gs_0,\gs_1)$ associated
with the dessin $(\RS,\sQ)$ contains the complete instructions for
building $\RS$: for every cycle of $\gs_0$ and every cycle of $\gs_1$
take a (convex planar) polygon with sides labeled by the elements in the cycle
in their cyclic order.
Next glue these polygons by identifying sides with the same label.
The sides of these polygons correspond to arrows of the quiver $\sQ$
and the vertices of these polygons correspond to vertices of $\sQ$.
After the above procedure of glueing polygons along their sides
one must still identify points which correspond to the same vertex of $\sQ$.
The result is then $\RS$ with $\sQ$ embedded in it.

\

\ssnnl{}\label{bi-adjacency constellation}
The $3$-constellation $(E,\gs_0,\gs_1)$ in \ref{M from constellation}
and the list $\MM$ in \ref{twist} both are completely determined by 
and do completely determine the dessin $(\RS,\sQ)$.
One can, however, also describe the relation between $(E,\gs_0,\gs_1)$ 
and $\MM$ in a direct algebraic/combinatorial way.

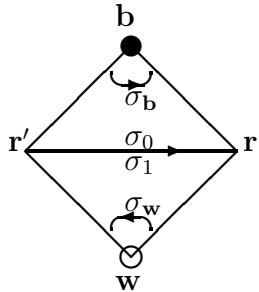
\begin{figure}[h]
\setlength{\unitlength}{0.5pt}
\begin{picture}(350,200)(-100,-20)
\thicklines
\put(200,0){
\begin{picture}(150,100)(0,0)
\put(0,0){\line(1,1){80}}
\put(0,0){\line(-1,1){80}}
\put(0,160){\line(-1,-1){80}}
\put(0,160){\line(1,-1){80}}
\put(0,0){\circle{15}}
\put(0,160){\circle*{15}}
\put(-80,80){\vector(1,0){120}}
\put(-80,80){\line(1,0){160}}
\put(-12,175){$\vb$}
\put(-12,-25){$\vw$}
\put(85,80){$\vr$}
\put(-93,80){$\vr'$}
\put(0,140){\oval(30,20)[b]}
\put(0,20){\oval(30,20)[t]}
\put(5,30){\vector(-1,0){15}}
\put(-5,130){\vector(1,0){15}}
\put(-5,85){$\gs_0$}
\put(-5,66){$\gs_1$}
\put(-5,115){$\gs_\vb$}
\put(-5,35){$\gs_\vw$}
\end{picture}
}
\end{picture}
\setlength{\unitlength}{1pt}
\caption{\label{fig:constellation}
\textit{orientation of cycles in $3$-constellation
and orientation of $\RS$.}}
\end{figure}

For the construction of $(E,\gs_0,\gs_1)$ from $\MM$ one first forms for $\vb\in\sB$, $\vw\in\sW$ and $i\in\{1,\ldots,N\}$
the sets
\begin{eqnarray*}
|\gs_\vb|&=&\{e\in E\:|\: \vb(e)=\vb\}\,,\qquad
|\gs_\vw|\:=\:\{e\in E\:|\: \vw(e)=\vw\}\,,\\
|\vz_i|&=&\{e\in E\:|\: i\in\{\vr(e),\vr'(e)\}\,\}\,.
\end{eqnarray*}
These sets have an unoriented cyclic structure:
$e,e'\in|\gs_\vb|$ (resp. $e,e'\in|\gs_\vw|$)  are neighbors if and only if 
$\{\vr(e),\vr'(e)\}\cap\{\vr(e'),\vr'(e')\}\neq\emptyset$, while 
$e,e'\in|\vz_i|$  are neighbors if and only if 
$\{\vb(e),\vw(e)\}\cap\{\vb(e'),\vw(e')\}\neq\emptyset$.
In order to put a consistent orientation on these cyclic sets we choose
one of the two possible orientations of $|\vz_1|$.
For every $\vb\in\sB$ we have $\sharp (|\vz_1|\cap|\gs_\vb|)=0$ or $2$.
In the latter case  $|\vz_1|\cap|\gs_\vb|$ consists of two elements,
say $e$ and $e'$, which are neighbors both in $|\vz_1|$ and in $|\gs_\vb|$.
The orientation on $|\vz_1|$ then induces an orientation on $|\gs_\vb|$ 
such that $e$ is the successor of $e'$ in $|\gs_\vb|$
if $e$ is the successor of $e'$ in $|\vz_1|$.
In this way the orientation on $|\vz_1|$ induces an orientation on 
every $|\gs_\vb|$ and every $|\gs_\vw|$ which has a non-empty intersection with
$|\vz_1|$. Next this induces an orientation on 
every $|\vz_i|$ which has a non-empty intersection with any of the already
oriented sets $|\gs_\vb|$ or $|\gs_\vw|$. And so on. It is because of the geometric background of $\MM$ (i.e. the orientations shown in the right-hand picture in
Figure \ref{fig:untwist helix}) that we can indeed go on.
In the end all sets $|\gs_\vb|$ and $|\gs_\vw|$ have an oriented cyclic 
structure and can be identified with cyclic permutations $\gs_\vb$ and 
$\gs_\vw$. 
The construction of the $3$-constellation $(E,\gs_0,\gs_1)$ is finished by setting
$$
\gs_0\,=\,\prod_{\vb\in\sB} \gs_\vb\,,\qquad
\gs_1\,=\,\prod_{\vw\in\sW} \gs_\vw^{-1}\,;
$$
the reason for inverting the white cycles is shown in Figure
\ref{fig:constellation}

To go from $(E,\gs_0,\gs_1)$ to $\MM$ one views the cycles of
the permutations $\gs_0$, $\gs_1$ and $\gs_0\gs_1^{-1}$  as elements of the
sets $E_{\gs_0}$, $E_{\gs_1}$ $E_{\gs_0\gs_1^{-1}}$, respectively, and then sets $\MM=\{(\vb(e),\,\vw(e),\,\vr(e),\,\vr'(e))\}_{e\in E}$
with
\begin{eqnarray*}
\vb(e)&=&\textrm{the cycle of $\gs_0$ which contains $e$}\,,
\\
\vw(e)&=&\textrm{the cycle of $\gs_1$ which contains $e$}\,,
\\
\vr(e)&=&\textrm{the cycle of $\gs_0\gs_1^{-1}$ which contains $e$}\,,
\\
\vr'(e)&=&\vr(\gs_1(e))\,.
\end{eqnarray*}

\

\ssnnl{ Example.}\label{dp3IVconstellation}
The procedure of \ref{bi-adjacency constellation} converts the list $\MM$ 
in \ref{correct dP3IV} into the $3$-constellation $(E,\gs_0,\gs_1)$
with 
\begin{eqnarray*}
E&=&\{1,2,3,4,5,6,7,8,9,10,11,12,13,14,15,16,17,18\}\,,
\\
\gs_0&=&(3,7,15)(4,12,16)(6,9,17)(5,11,18)(2,8,13)(1,10,14)\,,
\\
\gs_1&=&(4,8,15)(3,10,16)(5,7,17)(6,12,18)(1,9,13)(2,11,14)\,.
\end{eqnarray*}
Since in this example the cycle
$\vz_1=(1,10,3,7,5,11,2,8,4,12,6,9)$ meets each of the above $3$-cycles
the conversion went especially fast.

The corresponding superpotential $W$ is
\begin{eqnarray*}
&&\hspace{-1.5em}
X_{3}X_{7}X_{15}+X_{4}X_{12}X_{16}+X_{6}X_{9}X_{17}+X_{5}X_{11}X_{18}+X_{2}X_{8}X_{13}+X_{1}X_{10}X_{14}
\\ 
&&\hspace{-1.5em}
-X_{4}X_{8}X_{15}-X_{3}X_{10}X_{16}-X_{5}X_{7}X_{17}-X_{6}X_{12}X_{18}-
X_{1}X_{9}X_{13}-X_{2}X_{11}X_{14}.
\end{eqnarray*}

\section{The determinant of the bi-adjacency matrix and the secondary polygon}
\label{sec:bi-adjacency determinant}

\ssnnl{}\label{kasteleyn determinant}
Corollary \ref{black volA} implies that the bi-adjacency matrix $\KK^\wt=(\kappa_{ij})$ of the dessin $(\RS,\sQ)$ is a square matrix.
Its determinant is by definition
$$
\det\KK^\wt\,=\,\sum_{\tau}\mathrm{sign}(\tau) \kappa_{i\,\tau(i)}
$$
where $\tau$ runs over the set of all bijections
$\sB\stackrel{\simeq}{\rightarrow}\sW$ and $\mathrm{sign}(\tau)$ 
is defined as the sign of the permutation 
$\tau_0^{-1}\tau$ of $\sB$
for some reference bijection $\tau_0:\sB\stackrel{\simeq}{\rightarrow}\sW$.
Of course, if we want to actually write $\KK^\wt$ as a matrix, we must
choose bijections between $\sB$, $\sW$ and the set of numbers
$\{1,\ldots,\mathrm{vol}_\cA\}$ and that fixes $\tau_0$.
Changing the reference bijection multiplies
$\det\KK^\wt$ by $\pm 1$, but that is for our purpose unimportant. Next note that the $(\vb,\vw)$-entry of $\KK^\wt$
is non-zero
if and only if there is an edge of $\Gg$ connecting the nodes $\vb$ and $\vw$.
So the only bijections $\tau$ that contribute to the determinant of $\KK^\wt$
are the \emph{perfect matchings}; see Definition \ref{def:perfect matching}.
Thus we find
\begin{eqnarray}\nonumber
\det\KK^\wt&=&\pm \sum_{\mathrm{perfect\,matchings}\, P}\mathrm{sign}(P) 
\prod_{e\in P} w(e)u_{\vr(e)}u_{\vr'(e)}
\\
\label{eq:Kasteleyn det}
&=&\pm \sum_{\mathrm{perfect\,matchings}\, P}\mathrm{sign}(P) 
\left(\prod_{e\in P} \wt(e)\right) \vu^{\widehat{P}}
\end{eqnarray}
where 
\begin{equation}\label{eq:lift perfect}
\widehat{P}=\sum_{e\in P}(\ve_{\vr(e)}+\ve_{\vr'(e)})\,\in\,\ZZ^N\,,
\end{equation}
and $\vu^\vp\,=\,u_1^{p_1}\cdot\ldots\cdot u_N^{p_N}$ for 
$\vp=(p_1,\ldots,p_N)\in\ZZ^N$.

\

\ssnnl{ Definition.}\label{def:Newton polytope}
The \emph{Newton polytope} of a Laurent polynomial 
$$
f\,=\,\sum_{(k_1,\ldots,k_N)\in\ZZ^N} c_{(k_1,\ldots,k_N)}
u_1^{k_1}\cdot\ldots\cdot u_N^{k_N}\,\in\CC[u_1^\pm,\ldots,u_N^\pm]
$$
is the polytope
$$
\mathsf{Newton}(f)\,:=\,\mathsf{convex\;hull}\left(
\{(k_1,\ldots,k_N)\:|\:c_{(k_1,\ldots,k_N)}\neq 0\,\}\right)\,.
$$

\

\ssnnl{ Theorem.}\label{thm:Kasteleyn and secondary polygon}
\textit{The Newton polytope of the determinant of the bi-adjacency matrix $\KK^\wt$ of the dessin $(\RS,\sQ)$ is the same as the secondary polygon 
of $\LL\subset\ZZ^N$ (see \ref{def:Lsecondary polytope}):
\begin{equation}\label{eq:Newton is secondary}
\mathsf{Newton}(\det\KK_\wt)\,=\,\Sigma(\LL)\,.
\end{equation}
The vertices of $\mathsf{Newton}(\det\KK_\wt)$ are the points
$\widehat{P_C}$ (see (\ref{eq:lift perfect})) given by the 
perfect matchings $P_C$ with $C$ a $2$-dimensional cone in the secondary fan.
}

\proof
If $P$ and $P'$ are two perfect matchings then the edges in
$P$ oriented black-to-white together with the edges in
$P'$ oriented white-to-black form a collection of closed loops on the torus $\cT$ in 
\ref{graphs}. The first homology group of this torus is $\LL$.
Therefore $\widehat{P}-\widehat{P'}$ lies in $\LL\subset\ZZ^N$. 
This shows that $\mathsf{Newton}(\det\KK^\wt)$ lies in a $2$-dimensional plane
in $\RR^N$ parallel to $\LL_\RR$.

From Equation (\ref{eq:lift perfect}) one sees that for every
perfect matching $P$
$$
\widehat{P}=\sum_{i=1}^N\sharp\{e\in P\:|\:i\in\{\vr(e),\vr'(e)\}\}\:\ve_i\,,
$$
i.e. the $i$-th coordinate of the vector $\widehat{P}$ equals the number
of edges in the perfect matching $P$ which intersect the $i$-th zigzag loop
(cf. \ref{zigzagquiver}).
It follows from Proposition \ref{count cells}
that the number of edges of the $i$-th zigzag loop is
$\sum_{j=1}^N|\det(\vb_i,\vb_j)|$.
Since consecutive edges in a zigzag loop can never intersect the same perfect matching we see that for every perfect matching $P$ and every $i$
\begin{equation}\label{eq:secondary bound 1}
\sharp\{e\in P\:|\:i\in\{\vr(e),\vr'(e)\}\}
\leq\textstyle{\frac{1}{2}}\,\sum_{j=1}^N|\det(\vb_i,\vb_j)|\,.
\end{equation}
Fix $i$ and let $C$ and $C'$ be the two $2$-dimensional cones in the secondary fan
containing the half-line $\RR_{\geq 0}\vb_i$. Then 
$\{i,j\}\in L_C\Leftrightarrow \det(\vb_i,\vb_j)\geq 0$ and
$\{i,j\}\in L_{C'}\Leftrightarrow \det(\vb_i,\vb_j)\leq 0$.
Since $\sum_{j=1}^N\det(\vb_i,\vb_j)=0$ we find
$$
\sum_{j,\,\{i,j\}\in L_C}|\det(\vb_i,\vb_j)|\,=\,
\sum_{j,\,\{i,j\}\in L_{C'}}|\det(\vb_i,\vb_j)|\,=\,
\textstyle{\frac{1}{2}}\,\sum_{j=1}^N|\det(\vb_i,\vb_j)|\,.
$$
Thus we see that for $i$ and the perfect matchings $P_C$ and $P_{C'}$
Equation (\ref{eq:secondary bound 1}) is in fact an equality. 

Fix a $2$-dimensional cone $C$ in the secondary fan.
It is bounded by two half-lines
$\RR_{\geq 0}\vb_i$ and $\RR_{\geq 0}\vb_{i'}$.
The above argument now shows that 
Equation (\ref{eq:secondary bound 1}) is in fact an equality for
the perfect matching $P_C$ and $i$ and for $P_C$ and $i'$ in place of $i$.
Thus the point $\widehat{P_C}$ is a vertex of
$\mathsf{Newton}(\det\KK_w)$. This together with (\ref{eq:secondary bound 1})
proves:
$$
\mathsf{Newton}(\det\KK_\wt)\,=\,
\mathsf{convex\;hull}\left(\!\!
\{\widehat{P_C}\,|\,C\;\textrm{$2$-dim. cone of secondary fan of $\LL$}\}
\!\!\right).
$$
To complete the proof of (\ref{eq:Newton is secondary}) we note that
Equations (\ref{eq:LC}), (\ref{eq:psiC}), (\ref{eq:PCLC}) and 
(\ref{eq:lift perfect}) 
together with Proposition \ref{count cells} show
$\widehat{P_C}\,=\,\psi_C$.
\qed

\

\ssnnl{ Example.}\label{dP3 model II}
According to \cite{FHKVW} \S 8 the quiver for model II of $dP_3$
(i.e. case $B_8$ in Figure \ref{fig:dp3 models}) affords two 
different superpotentials. This is confirmed by our algorithm,
which yields the following two bi-adjacency matrices
$$
\KK^\wt_1=\left[\begin{array}{cccc}
\wt_2u_1u_2&\wt_8u_2u_4&0&\wt_4u_1u_4\\
\wt_{14}u_5u_6&\wt_{10}u_2u_6&\wt_9u_2u_5&0\\
\wt_6u_1u_6+\wt_7u_2u_3&\wt_{13}u_4u_6&\wt_1u_1u_2&\wt_{11}u_3u_4\\
\wt_{12}u_3u_5&0&\wt_5u_1u_5&\wt_3u_1u_3
\end{array}\right]
$$
$$
\KK^\wt_2=\left[\begin{array}{cccc}
\wt_2u_1u_2&\wt_9u_2u_5&0&\wt_5u_1u_5\\
\wt_4u_1u_4&\wt_1u_1u_2&\wt_8u_2u_4&0\\
\wt_7u_2u_3&\wt_{14}u_5u_6&\wt_{10}u_2u_6&\wt_{12}u_3u_5\\
\wt_{11}u_3u_4&\wt_6u_1u_6&\wt_{13}u_4u_6&\wt_3u_1u_3
\end{array}\right]
$$
The corresponding superpotentials are:

\begin{eqnarray*}
W_1&=&
X_{2}X_{8}X_{4}+X_{9}X_{14}X_{10}+X_{1}X_{7}X_{11}X_{13}X_{6}+X_{3}X_{12}X_{5}
\\ &&
-X_{2}X_{7}X_{12}X_{14}X_{6}-X_{8}X_{13}X_{10}-X_{1}X_{9}X_{5}-X_{3}X_{11}X_{4}
\end{eqnarray*}

\begin{eqnarray*}
W_2&=&
X_{2}X_{9}X_{5}+X_{1}X_{8}X_{4}+X_{7}X_{12}X_{14}X_{10}+X_{3}X_{11}X_{13}X_{6}
\\ &&
-X_{2}X_{7}X_{11}X_{4}-X_{1}X_{9}X_{14}X_{6}-X_{8}X_{13}X_{10}-X_{3}X_{12}X_{5}
\end{eqnarray*}
After figuring out the rule for translating the edge labels
one finds that the prepotentials in Eqs. (8.1) and (8.2) of \cite{FHKVW} 
are $W_A=W_1$ and $W_B=W_2$.

For both superpotentials the surface $\RS$ has genus $1$.
Figures 18 and 19 in \cite{FHKVW} show pictures of the planar periodic bi-partite graph
which is the lifting of $\Gg$ to the simply connected cover of $\RS$.

One can not expect the Kasteleyn matrices $K_A$ and $K_B$ in
\cite{FHKVW} Eq. (8.3) to be the same as the bi-adjacency matrices
$\KK^\wt_1$ and $\KK^\wt_2$, respectively, because of the ``untwist''.
$\KK^\wt_1$ and $\KK^\wt_2$ are Kasteleyn matrices of the dimer
models before the untwisting (see \ref{def:Kasteleyn matrix}).
Moreover Eqs. (8.4) and (8.5) of \cite{FHKVW} show that the Newton polygons of the determinants of the Kasteleyn matrices $K_A$ and $K_B$ are different.
On the contrary, for both $\KK^\wt_1$ and $\KK^\wt_2$ the determinant
is a linear combination of the monomials
$$
\vu^{[3,3,1,0,0,1]}\!,\,\vu^{[2,2,1,1,1,1]}\!,\,\vu^{[1,1,1,2,2,1]}
\!,\,\vu^{[2,0,0,2,2,2]}\!,\,
\vu^{[3,1,0,1,1,2]}\!,\,\vu^{[0,2,2,2,2,0]}\!,\,\vu^{[1,3,2,1,1,0]}\,.
$$
So,
in agreement with Theorem \ref{thm:Kasteleyn and secondary polygon},
the Newton polygons of $\det\KK^\wt_1$ and $\det\KK^\wt_2$ coincide with
the polygon in Figure \ref{fig:dp3 models} case $B_8$
and \ref{explict secondary polygon}.
In spite of having the same Newton polygon $\det\KK^\wt_1$ and $\det\KK^\wt_2$ 
are not equal; indeed when all weights are $1$ they differ in the coefficient
of $\vu^{[2,2,1,1,1,1]}$.

\section{Bi-adjacency matrix with critical weights
and the principal $\cA$-determinant}
\label{sec:A determinant}

\ssnnl{}\label{intro Adet}
In \cite{gkz4} p.297 Eq. (1.1) Gelfand, Kapranov and Zelevinsky define for a set $\cA\,=\,(\va_1,\ldots,\va_N)\subset\ZZ^{k+1}$ as in \ref{introGKZ} the
\emph{principal $\cA$-determinant $E_\cA(f_\cA)$}. The definition uses the
Laurent polynomial
\begin{equation}\label{eq:fA}
f_\cA\,=\,\sum_{i=1}^N u_i\vx^{\va_i}\,,
\end{equation}
where $\vx^\vm=x_1^{m_1}x_2^{m_2}\cdot\ldots\cdot x_{k+1}^{m_{k+1}}$
for $\vm=(m_1,\ldots,m_{k+1})\in\ZZ^{k+1}$ and where the coefficients $u_1,\ldots,u_N$ are variables.  
The name ``principal $\cA$-determinant'' refers to the fact 
(\cite{gkz4} p.298 Prop. 1.1) that
$E_\cA(f_\cA)$ can be written as the 
determinant of some exact complex, i.e. 
$$
E_\cA(f_\cA)=\prod_j(\det\,M_j)^{(-1)^j}
$$
for certain matrices $M_j$.
Another useful description of $E_\cA(f_\cA)$ is given in
\cite{gkz4} p.299 Prop. 1.2:
\begin{equation}\label{eq:Adetdisc}
E_\cA(f_\cA)=\pm
\prod_{\gG\subset\mathsf{convex\;hull}(\cA)}
\Delta_{\cA\cap\gG}(f_{\cA\cap\gG})^{m(\gG)}\,
\end{equation}
where the product runs over all faces $\gG$ of the
primary polytope $\mathsf{convex\;hull}(\cA)$ (cf. \ref{A secondary});
$m(\gG)$ is a multiplicity and $\Delta_{\cA\cap\gG}(f_{\cA\cap\gG})$ is 
the \emph{$(\cA\cap\gG)$-discriminant} of the Laurent polynomial
$f_{\cA\cap\gG}=\sum_{i:\,\va_i\in\gG} u_i\vx^{\va_i}$.
The latter discriminant is a polynomial in the variables 
$u_i$ with $\va_i\in\gG$. To define it (see \cite{gkz4} p.271) one needs
the algebraic set $\nabla_{\cA\cap\gG}$ which is the closure in 
$\CC^{\cA\cap\gG}$ of:
$$
\{\vu\in\CC^{\cA\cap\gG}\;|\;\exists
\vx_0\in(\CC^*)^{k+1}\;s.t.\; 
f_{\cA\cap\gG}(\vx_0)=\frac{\partial f_{\cA\cap\gG}}{\partial x_i}(\vx_0)= 0,
\:\forall i\}\,.
$$
Then, by definition, $\Delta_{\cA\cap\gG}(f_{\cA\cap\gG})\,=\,1$ if $\mathrm{codim}_{\CC^{\cA\cap\gG}}(\nabla_{\cA\cap\gG})>1$
and
$$
\mathsf{zero\; locus\; of\; } \Delta_{\cA\cap\gG}(f_{\cA\cap\gG})\,=\,
\nabla_{\cA\cap\gG}\qquad\textrm{if}\quad
\mathrm{codim}_{\CC^{\cA\cap\gG}}(\nabla_{\cA\cap\gG})=1\,.
$$
So, $E_\cA(f_\cA)$ gives the locus of the points $(u_1,\ldots,u_N)\in\CC^N$
for which at least one of the Laurent polynomials $f_{\cA\cap\gG}$
has a critical point with critical value $0$.

\

\ssnnl{ Theorem.}\label{thm:secAdet}(\cite{gkz4}p.302 Thm.1.4;
cf.(\ref{eq:gkz secondary polytope}))\\
\textit{The Newton polytope of $E_\cA(f_\cA)$ coincides with the secondary
polytope $\Sigma(\cA)$}.

\qed

\

\ssnnl{}\label{motivate conjecture}
Theorems \ref{thm:secAdet} and \ref{thm:Kasteleyn and secondary polygon}
in combination with Formula (\ref{eq:two secondary polytopes})
make one wonder whether for an appropriate choice of the weight $\wt$
the determinant of the bi-adjacency matrix $\KK^\wt$ is equal
to the principal $\cA$-determinant $E_\cA(f_\cA)$, up to the simple
transformation necessitated by (\ref{eq:two secondary polytopes}).
In order to formulate this transformation we must make the
dependence on $u_1,\ldots,u_N$ visible by writing
$\KK^\wt(u_1,\ldots,u_N)$ and $E_\cA(f_\cA(u_1,\ldots,u_N))$.

After some experimenting with examples I found
a very natural and simple choice for the weight that does the job.

\

\ssnnl{ Definition.}\label{def:critical weight}
The \emph{critical weight} for the arrows of the quiver $\sQ$
is the function
\begin{equation}\label{eq:def critical weight}
\crit: E\rightarrow\ZZ_{>0}\,,\qquad \crit(e)=\sharp\{e'\in E\,|\,
s(e')=s(e),\,t(e')=t(e)\}
\end{equation}

\

\ssnnl{ Conjecture.}\label{conj:critical weight}
\textit{For every set $\cA$ as in \ref{introGKZ} and every
dessin $(\RS,\sQ)$ (see \ref{final M}) constructed from $\cA$ by the algorithm
in Section \ref{sec:dessin} the determinant of the bi-adjacency matrix with critical weight and the principal $\cA$-determinant satisfy:}
\begin{equation}\label{eq:crit det}
(u_1\cdot\ldots\cdot u_N)^{\mathrm{vol}_\cA}\det\,\KK^\crit(u_1^{-1},\ldots,u_N^{-1})
\,=\,E_\cA(f_\cA(u_1,\ldots,u_N))\,.
\end{equation}

\

\ssnnl{ Remark.}\label{secvertices}
\emph{In support of the above conjecture we can point out that the
coefficients of the monomials corresponding to the vertices of the secondary polygons are, up to sign, the same for the two sides of Equation 
(\ref{eq:crit det}).} Indeed,
Theorem 1.4 of \cite{gkz4}p.302 gives the coefficient of the monomial in
the principal $\cA$-determinant with corresponds to a vertex of the secondary polytope. Such a vertex corresponds to a maximal cone $C$ 
in the secondary fan. In the notations of
Definition \ref{def:Lsecondary polytope} the formula in loc. cit. for the
coefficient of the monomial corresponding to $C$ reads:
\begin{equation}\label{eq:vertexcoeff}
\pm \prod_{\{i,j\}\in L_C}|\det(\vb_i,\vb_j)|^{|\det(\vb_i,\vb_j)|}\,.
\end{equation}
On the other hand, Theorems \ref{thm:perfect and secondary fan}
and \ref{thm:Kasteleyn and secondary polygon} show that the same vertex 
corresponds to a perfect matching $P_C$. From the role of perfect matchings in
the computation of the determinant of the bi-adjacency matrix
(see \ref{kasteleyn determinant}) one now easily checks that the coefficient
of the monomial in the determinant, which corresponds to $C$, is
(possibly up to sign) the same as (\ref{eq:vertexcoeff}).

\

In the remainder of this section we explicitly verify Conjecture \ref{conj:critical weight} in some examples.

\

\ssnnl{ Example.}\label{dP1 Kasteleyn}
For case $B_2\,=
\left[\begin{array}{rrrr}
0&1&1&-2\\
-1&0&2&-1
\end{array}\right]$ 
in Figure \ref{fig: models} the algorithm in Section \ref{sec:dessin} yields
in  the following bi-adjacency matrix with critical weights:
$$
\KK^\crit\,=\,
\left[
\begin{array}{ccc}
2u_1u_4&\quad 3u_3u_4&\quad u_1u_3\\
2u_2u_3&\quad u_2u_4&\quad 3u_3u_4\\
u_1u_2+3u_3u_4&\quad 2u_2u_3&\quad 2u_1u_4
\end{array}
\right]\,.
$$
One easily computes
$$
det\,\KK^\crit\,=\,27\vu^{[0,0,3,3]}+4\vu^{[1,2,3,0]}+4\vu^{[2,1,0,3]}
-18\vu^{[1,1,2,2]}-\vu^{[2,2,1,1]}
$$
Note that the exponents are the same as the coordinates of the vertices
and the interior point
of the secondary polygon in \ref{secondary fan dP1}.
For the computation of the principal $\cA$-determinant we note that
in this case (see also Example \ref{cubic roots})
$$
\cA\,=\,(\va_1,\,\va_2,\,\va_3,\,\va_4)\,=\,
\left(
\left[\begin{array}{c}1\\ 0\end{array}\right],\, 
\left[\begin{array}{c}1\\ 3\end{array}\right],\,
\left[\begin{array}{c}1\\ 1\end{array}\right],\,
\left[\begin{array}{c}1\\ 2\end{array}\right]\right)
$$
and hence
$f_\cA=x_1(u_1\,+\,u_2x_2^3\,+\,u_3x_2\,+\,u_4x_2^2)$.
The primary polytope is shown in Figure \ref{fig:secondaries}.
The two boundary points of this primary polytope contribute factors
$u_1$ and $u_2$ to the principal $\cA$-determinant (cf. Equation 
(\ref{eq:Adetdisc})). The contribution from the full primary polytope is
the very classical discriminant of the cubic polynomial
$u_1+u_3x_2+u_4x_2^2+u_2x_2^3$, which is 
(e.g. \cite{gkz4} p.405)
$$
27 u_1^2u_2^2+4u_1u_4^3+4u_3^3u_2-u_3^2u_4^2-18u_1u_3u_4u_2\,.
$$
Thus we find
$$
E_\cA(f_\cA)\,=\,
27\vu^{[3,3,0,0]} +4\vu^{[2,1,0,3]}+4\vu^{[1,2,3,0]}
-\vu^{[1,1,2,2]}-18\vu^{[2,2,1,1]}\,.
$$
\emph{This shows that conjecture \ref{conj:critical weight} holds in this case.}

\

\ssnnl{ Example.}\label{dP3 determinants}
When critical weights are used the two bi-adjacency matrices in 
Example \ref{dP3 model II} become
$$
\KK^\crit_1=\left[\begin{array}{cccc}
2u_1u_2&u_2u_4&0&u_1u_4\\
u_5u_6&u_2u_6&u_2u_5&0\\
u_1u_6+u_2u_3&u_4u_6&2u_1u_2&u_3u_4\\
u_3u_5&0&u_1u_5&u_1u_3
\end{array}\right]
$$
$$
\KK^\crit_2=\left[\begin{array}{cccc}
2u_1u_2&u_2u_5&0&u_1u_5\\
u_1u_4&2u_1u_2&u_2u_4&0\\
u_2u_3&u_5u_6&u_2u_6&u_3u_5\\
u_3u_4&u_1u_6&u_4u_6&u_1u_3
\end{array}\right]
$$
One easily computes $\det\,\KK^\crit_1$ and $\det\,\KK^\crit_2$
and finds that both are equal to
\begin{eqnarray*}
&&
4\vu^{[3,3,1,0,0,1]}\,+2\vu^{[1,1,1,2,2,1]}\,+\vu^{[1,3,2,1,1,0]}
\,+\vu^{[3,1,0,1,1,2]}\,-
\\ &&\hspace{7em}
\,-\vu^{[2,0,0,2,2,2]}\,-\vu^{[0,2,2,2,2,0]}\,-6\vu^{[2,2,1,1,1,1]}.
\end{eqnarray*}
This gives a remarkable contrast with the last line of Example 
\ref{dP3 model II} and underlines the role of the critical weight in making 
Conjecture \ref{conj:critical weight} reasonable for every dessin obtained from
$\cA$.

Transforming the above determinant as in the left-hand side
of Equation (\ref{eq:crit det}) yields
\begin{eqnarray}
\nonumber
&&\hspace{-3em}
4\vu^{[1,1,3,4,4,3]}\,+2\vu^{[3,3,3,2,2,3]}\,+\vu^{[3,1,2,3,3,4]}
\,+\vu^{[1,3,4,3,3,2]}\,-\vu^{[2,4,4,2,2,2]}\,-
\\ \nonumber &&\hspace{14em}
\,-\vu^{[4,2,2,2,2,4]}\,-6\vu^{[2,2,3,3,3,3]}
\\[1.5ex] \label{eq:factor det}
&&\hspace{-4em}=\:
u_1u_2u_3^2u_4^2u_5^2u_6^2\,(u_4u_5-u_1u_2)
(4u_3u_4u_5u_6+u_2^2u_3^2-2u_1u_2u_3u_6+u_1^2u_6^2).
\end{eqnarray}
From the matrix $B_8$ in Figure \ref{fig:dp3 models} one easily finds 
$\cA$ and the corresponding polynomial
$$
f_\cA\,=\, u_1x_1\,+\,u_2x_1x_2\,+\,u_3x_1x_3\,+\,u_4x_1x_4
\,+\,u_5x_1x_2x_4^{-1}\,+\,u_6x_1x_2x_3\,.
$$
The primary polytope $\mathsf{convex\;hull}(\cA)$ is $3$-dimensional
and has seven $2$-di\-mensional faces, six of which are triangles and 
one is a quadrangle. Its only integer points are its vertices and these
correspond to the monomials of $f_\cA$. 
One can compute the discriminants for the polynomials supported by the faces of the primary polytope and multiply these to get the principal $\cA$-determinant
$E_\cA(f_\cA)$ as in Equation (\ref{eq:Adetdisc}).
The result is exactly as in (\ref{eq:factor det}), with the $4$-term factor coming from the primary polytope itself, the $2$-term factor coming from 
the quadrangle face and the other factors with multiplicities
coming from the vertices.
 \emph{This shows that conjecture \ref{conj:critical weight} holds in this case.}

\

\ssnnl{ Further examples.}\label{further}
Exactly as in Example \ref{dP1 Kasteleyn} one can verify Conjecture 
\ref{conj:critical weight} for other ``four-nomials''
of degree $\leq 5$ (see Example \ref{Labc}) by using the formulas
in \cite{gkz4} for discriminants of quartic and quintic polynomials.

The method of Example \ref{dP3 determinants} can be successfully employed
to verify Conjecture \ref{conj:critical weight} for, for instance,
cases $B_5$, $B_6$, $B_{10}$ in Figure \ref{fig:dp3 models}.
It also works for
$B=
\left[\begin{array}{rrrrrrr}-3&1&1&1&0&0&0\\ -3&0&0&0&1&1&1\end{array}\right]$
and
$B=
\left[\begin{array}{rrrrr}1&1&0&-2&0\\ 0&0&1&1&-2\end{array}\right]$,
which are closely related to $B_{16}$ and $B_9$, respectively.

\

\end{document}